\documentclass[a4paper,10pt,fleqn]{article}
\usepackage{textpos}
\usepackage{amsmath,amssymb,amsfonts,graphicx,amsthm,colortbl}
\usepackage[a4paper, margin=3cm]{geometry}
\usepackage{amsthm}
\usepackage{dsfont}
\usepackage{booktabs}
\usepackage{caption,subcaption}
\usepackage{tikz}
\usetikzlibrary{patterns}
\usepackage{pgfplots}
\usepackage{standalone}
\usepackage{placeins}

\newif\ifdraft

\ifdraft
\usepackage{showlabels}
\usepackage[colorinlistoftodos]{todonotes}
\else

\fi

\allowdisplaybreaks

\newcommand{\beq}{\begin{eqnarray*}}
\newcommand{\eeq}{\end{eqnarray*}}
\newcommand{\beqn}{\begin{eqnarray}}
\newcommand{\eeqn}{\end{eqnarray}}
\newcommand{\notinclude}[1]{}

\usepackage[pdftex]{hyperref}
\hypersetup{
    colorlinks = true,
    citecolor=black,
    filecolor=black,
    linkcolor=black,
    urlcolor=black,
    pdftitle={Stochastic Dominance in Elastic Shape Optimization},
    bookmarksopen=true,
    a4paper,
    pdfdisplaydoctitle
}

\newcommand{\captionfonts}{\small \em}

\makeatletter  \long\def\@makecaption#1#2{%
  \vskip\abovecaptionskip
  \sbox\@tempboxa{{\captionfonts #1: #2}}%
  \ifdim \wd\@tempboxa >\hsize
    {\captionfonts #1: #2\par}
  \else
    \hbox to\hsize{\hfil\box\@tempboxa\hfil}%
  \fi
  \vskip\belowcaptionskip}
\makeatother   
\theoremstyle{definition}

\newcommand{\R}{\mathbb{R}}

\newcommand{\body}{\mathcal{O}}

\newcommand{\lObjective}{\ensuremath{\mathcal{J}}}

\newcommand{\prob}{\mathbb{P}}
\newcommand{\costThreshold}{\eta}
\newcommand{\expectedValue}{\mathbb{E}}
\newcommand{\I}{\mathcal{I}}

\newcommand{\V}{\mathcal{V}}
\newcommand{\FE}{\mathbf{FE}}
\newcommand{\Qev}{\mathbf{Q}_\expectedValue}
\newcommand{\Qee}{\mathbf{Q}_{\text{EE}_\costThreshold}}
\newcommand{\Qep}{\mathbf{Q}_{\text{EP}_\costThreshold}}

\newcommand{\domain}{D}
\newcommand{\object}{\mathcal{O}}

\newcommand{\epsu}{\epsilon[u]}

\newcommand{\E}{\mathcal{E}}
\renewcommand{\L}{\mathcal{L}}

\newcommand{\W}{\mathcal{W}}
\newcommand{\Comp}{\mathcal{C}}

\renewcommand{\d}{\,\mathrm{d}}

\newcommand{\hv}[1]{\mathcal{H} (#1)}
\newcommand{\hvs}[1]{\mathcal{H}_\gamma (#1)}
\newcommand{\maxs}{{\max}_\gamma}

\newcommand{\leqone}{\preceq_{\text{st}}}
\newcommand{\leqtwo}{\preceq_{\text{icx}}}

\newcommand{\heavreg}{\gamma}

\newcommand{\emphdef}[1]{\textit{#1}}

\DeclareMathOperator{\tr}{tr}

\newcommand{\charfunction}{\mbox{\Large $\chi$}}

\begin{document}

\newpage
\setcounter{page}{1}

\title{Stochastic Dominance Constraints 
in Elastic Shape Optimization}

\author{
Sergio Conti\footnotemark[1]
\and Martin Rumpf\footnotemark[3]
\and R\"{u}diger Schultz\footnotemark[2]
\and Sascha T{\"o}lkes \footnotemark[3]
}

\renewcommand{\thefootnote}{\fnsymbol{footnote}}
\footnotetext[1]{Institut f\"ur Angewandte Mathematik, Universit\"at Bonn, Endenicher Allee 60, 53115 Bonn,
  Germany, sergio.conti@uni-bonn.de}
\footnotetext[2]{Faculty of  Mathematics, 
University of  Duisburg-Essen, Thea-Leymann-Str.~9, 45127 Essen, Germany, schultz@math.uni-duisburg.de}
\footnotetext[3]{Institute for Numerical Simulation, Universit\"at Bonn, Endenicher Allee 60, 53115 Bonn, Germany, \{martin.rumpf\}\{sascha.toelkes\}@ins.uni-bonn.de}
\renewcommand{\thefootnote}{\arabic{footnote}}

\newcommand{\rcomment}[3]{\item[#1]\hypertarget{#1}{{\it #2}\\ {#3}}}
\newenvironment{rlist}{\begin{list}{$\cdot$}{
\renewcommand{\makelabel}[1]{{\bf {##1} }}
\setlength{\labelsep}{1pt}
\setlength{\labelwidth}{2.5cm}
\setlength{\leftmargin}{2.5cm} }}{\end{list}}

\newcommand{\link}[1]{\hyperlink{#1}{{\rm [{\color{blue} #1}]}}}
\newcommand{\change}[2]{\marginpar{{\parbox{\marginparwidth}{\hypertarget{#1}{\tiny {\bf #1}}\\#2 }}}}
\newcommand{\changedtext}{\color{red}}

\newpage

\maketitle

\begin{abstract}
This paper deals with shape optimization for elastic materials under stochastic loads.
It transfers the paradigm of stochastic dominance,
which allows for flexible risk aversion via comparison with benchmark random variables,
from finite-dimensional stochastic programming 
to shape optimization.  
Rather than handling risk aversion in the objective, this enables
risk aversion by including
dominance constraints that single out subsets of nonanticipative shapes which compare favorably to a chosen stochastic benchmark. 
This new class of stochastic shape optimization problems arises by optimizing over such feasible sets. The analytical description is built on risk--averse cost measures. The underlying cost functional is of compliance type plus 
a perimeter term,  in the implementation shapes are represented by a phase field which permits an easy estimate of a regularized perimeter.
The analytical description and the numerical implementation of dominance constraints are
built on risk-averse measures for the cost functional.
A suitable numerical discretization is obtained using finite elements both for the displacement and the phase field function.
Different numerical experiments demonstrate the potential of the proposed stochastic shape optimization model 
and in particular the impact of high variability of forces or probabilities in the different realizations.
\end{abstract}


\section{Introduction}\label{sec:introduction}
This paper picks up concepts of risk aversion from finite dimensional stochastic optimization.
Two major lines of research were pursued in this context.
At first, risk aversion was considered in the objective function via statistical parameters, so called risk measures 
(cf. the monograph \cite{PfRo07}).
Practically meaningful risk measures were identified, which enable multiperiod decision making and led 
to appropriate mathematical tools. This includes investigations of mean-risk models in a two-stage linear setting \cite{Ah06,Kr05}, in a two-stage mixed-integer linear setting \cite{ScTi06,RoVi08}, and in a multistage setting \cite{AsRu15,Sh12}.
The second line of research concerns risk aversion in the constraints using
concepts of stochastic dominance. Rather than heading for risk minimization,
this work aims at bounding risk with the help of benchmark random variables
and comparison via partial orders of random variables. In this way,
``acceptance'' of nonanticipative solutions can be made mathematically
rigorous, offering to optimize suitable objectives over such ``acceptable
sets''. 
The topic of introducing orders on families of random variables has  a long tradition in applied probability and statistics, see \cite{MuSt02}.
Basic motivation for these efforts is to refine the rating of random variables and go beyond the  procedure of evaluating random variables just by a single stochastic parameter, the mean, for instance. 
The latter may be not very informative, or may suffer from non-existence of the parameter. It also may happen that for specific applications at hand additional information is available which may be incorporated into the comparison criteria.
Incorporation of dominance constraints
into stochastic programs (with decision variables in finite as well as
infinite dimensional Banach spaces) was pioneered by Dentcheva and
Ruszczy\'nski \cite{DeHeRu07,DeRu03,DeRu04,DeRu08}. 
Numerical techniques in risk averse stochastic optimization as summarized
above are inherently finite dimensional and appealing to principles of linear,
mixed-integer linear, or nonlinear programming. To the best of our knowledge,
there is no previous work on numerical techniques for risk averse shape
optimization where nonanticipative domains arise as variables. 
Research on stochastic dominance in optimization under  uncertainty has seen a substantial increase during the last decade. Topics include basic analysis of models \cite{DeRu03, DeHeRu07, Lu08,GoNeSc08,GoGoSc11,
DeRo13,DeRu14}, algorithm development  \cite{RuRu08,DeRu10,DrSc10}, 
and industrial applications, see for instance \cite{CaGoSc09,FaMiRo11,GoGoNe11}.

Uncertainty in the context of shape optimization already attracted considerable attention.
Multiload approaches take into account a fixed and usually small number of 
loading configurations \cite{AlJo05,GuRoBe03}.
Applications of robust worst-case optimization to shape optimization can be
found in~\cite{BeNe02}.
Shape optimization with stochastic loading was discussed in the 
context of beam models in \cite{Me01a} and in aerodynamic design in \cite{ScSc15,ScSc09}.
A number of papers addressed worst-case optimization, e.g. \cite{BeKoNe99,BaNe07}.
An efficient optimization approach for the optimization of the expected value of compliance and tracking type cost functionals under stochastic loading,
which makes use  of the representation of 
realizations of surface and volume loads as linear combinations of a few basis modes,
was developed and used in \cite{CoHePa07}.
 In  \cite{CoHePa11}
it was shown how this approach can be used for the optimization of risk averse cost measures,
such as expected excess and excess probability.
Allaire and Dapogny \cite{AlDa14} used a linearization of the cost functional in a worst case optimization scenario and used adapted duality techniques to derive an efficient optimization algorithm. In \cite{AlDa15} they
investigated different types of uncertainties, e.g., in the loading and the material parameters, and minimized stochastic cost functionals such as the expected cost of the failure probabilities. Their approach is based on the Taylor expansion of the risk measure and also leads to deterministic algorithm with a cost depending on the number of realizations of the random configuration.
Recently, Dambrine, Dapogny and Harbrecht \cite{DaDaHa15} studied elastic shape optimization with stochastic surface loads.
They developed a deterministic algorithm to minimize the expected cost which is based on the observation the objective functional 
and its gradient are determined by the first order moments of the surface loads. Their method relies on the efficient approximation of 
integrals in $6$ dimensions.
 Pach studied the concept of stochastic dominance in the context of a parametric shape optimization approach, 
where certain thickness parameters of a given truss geometry have to be optimized \cite{Pa13}.

The paper is organized as follows.
In Section \ref{sec:shapeopt} we will revisit basic concepts in shape optimization under stochastic loading and prepare the later modeling
of stochastic dominance. 
Then, Section \ref{sec:dominance} will introduce the notion of stochastic dominance constraints in elastic shape optimization based on benchmark shapes.  The spatial discretization via finite elements and a suitable regularization of the constraints is proposed in Section 
\ref{sec:discrete}. Finally, in Section \ref{sec:numerics} we discuss properties and characteristics of the proposed approach for two different 
applications.

\section{Shape optimization with stochastic loading}\label{sec:shapeopt}
Shape optimization based on an explicit parametric description of the mechanical object is algorithmically quite demanding. Hence, we consider an approximating phase-field representation of the elastic object $\object \subset \domain$, where 
$\domain\subset \R^d$ ($d=2,3$) denotes the computational domain. 
Here, we follow \cite{PeRuWi12}, where a phase field appraoch for 
(nonlinear) elastic shape optimization was considered including existence results 
and $\Gamma$-convergence of the cost functional in the context of constraint optimization.
To this end, we take into account a phase-field function $v:\domain \to\R$ 
of Allen--Cahn or Modica--Mortola type with a phase field energy functional 
\begin{equation*}
\L^\varepsilon[v]:=\frac12\int_\Omega\varepsilon|\nabla v|^2+\frac1\varepsilon\Psi(v)\,\d x\,,
\end{equation*}
where the scale parameter $\varepsilon$ describes the width of the interfacial region.
In our context of shape modeling, we set $\Psi(v):=\frac{9}{16}(v^2-1)^2$,
which has two minima at $v=-1$ and $v=1$ representing
the two phases outside and inside of $\object$, respectively.
In the limit $\varepsilon\to0$, the phase field $v$ is forced towards the pure phases
$-1$ and $1$ and $\L^\varepsilon$ $\Gamma$-converges to the total interface area \cite{Br02}.
Furthermore, we approximate the total volume by the smooth functional $\V[v]:=\int_\domain\charfunction_\object(v)\,\d x$,
where $\charfunction_\object(v):=\frac14(v+1)^2$ is an approximation to the characteristic function $\charfunction_\object$.
Next, let recall the model of linearized elasticity taylored to the phase field description of the elastic object $\object$.
Instead of considering a void phase on $\domain \setminus \object$, we follow common practice in shape optimization assuming the presence of
a very soft material on the part of the domain 
 outside of $\object$.
The \emphdef{elastic energy} stored inside the material under a given displacement $u:\domain \to \R^d$
is then defined as
\begin{equation}
\W^\delta [v, u] := \frac{1}{2} \int_\domain ((1 - \delta)
\charfunction_\object(v) + \delta) C \epsu : \epsu \d x\,,
\end{equation}
where $\epsu = \frac12 (Du^T+Du)$ is the strain tensor. 
Here, $C$ is a fourth-order tensor satisfying the symmetry
relations
$C_{ijkl} = C_{jikl} = C_{ijlk} = C_{klij}$
and $\sum_{ijkl} C_{ijkl} \xi_{ij} \xi_{kl} \geq c \|\xi+\xi^T\|^2$ for all $\xi \in \mathbb{R}^{3 \times 3}$
with $c > 0$. In our implementation we consider only isotropic materials, with the
the Lam{\'e}-Navier elasticity tensor $C$ defined by $\sigma[u] = C
\epsu = \lambda \tr \epsu \mathds{1} + 2 \mu \epsu$.
Furthermore, we take into account a Dirichlet boundary $\Gamma_D\subset \partial \domain$ 
and an inhomogeneous Neumann boundary $\Gamma_N\subset \partial \domain$ with $\Gamma_D \cap \Gamma_N = \emptyset$. 
Since external forces can only be applied on the elastic object $\object$, we require
$v|_{\Gamma_N}=1$. 
The equilibrium displacement  $u\in H^1(\domain;\R^d)$ for fixed phase field $v$ 
is then given as the minimizer of the free energy
\begin{equation*}
\E[v,u]:=\W^\delta[v,u]-\Comp[u]
\end{equation*}
within a set of admissible displacements $u$ with trace $u|_{\Gamma_D}=0$, where
the compliance functional 
\begin{equation*}
\Comp[u]:=\int_{\Gamma_N}g\cdot u\,\d a
\end{equation*}
is the (negative) potential of the surface load for 
$g \in L^2(\Gamma_N;\R^d)$.

Finally, we define the cost functional
\begin{equation}\label{eqdefJ}
\mathcal{J}[v, u] := 2 \W^\delta[v,u] + \nu \V[v] + \eta \L^\varepsilon[v]\,,
\end{equation} where 
$\V$ measures approximately the volume of the elastic object $\object$ and $\mathcal{L}^\varepsilon$ the
object perimeter within the computational domain $\domain$.
The shape optimization problem is now to minimize $\mathcal{J}[v, u[v]]$ subject to the constraint that 
$u[v]$ is a minimizer of $u\mapsto \E [v,u]$.

The target functional $\lObjective$ 
defined in (\ref{eqdefJ}) depends deterministically on the load profiles.
If the surface load $g$ is uncertain, i.e., described by a random variable $g(\omega)$ on some  probability space $(\Omega,\mathcal{A},\prob)$, then so are the displacement $u(\omega)$ and the compliance  target functional $\lObjective[v,\omega]$. 
From our global perspective  the random variable $g(\omega)$ describes the surface loads which, together with the shape of the elastic object described by the phase field $v$, determine 
the displacement $u[v](\omega)$  and finally the compliance $\lObjective[v,u(\omega)]$. Optimizing or selecting shapes therefore amounts to optimizing the random variable representing compliance, 
\begin{equation}\label{eq:family}
\left\{\lObjective[v,u[v[(\cdot)]\;\;:\;\; v \in H^{1,2}_{\Gamma_N}(\domain)\right\},
\end{equation}
where
the Sobolev space $H^{1,2}_{\Gamma_N}(\domain)$ is the set of admissible phase fields which obey $v=1$ on $\Gamma_N$ in the sense of traces.

The members of the family being measurable functions further conceptual efforts have to be made to arrive at well-posed optimization or feasibility problems, on which we will briefly report here for the purpose of later comparison.
In mean-risk models objective functions are formulated  as weighted  sums of the  expected value and some risk measure, which is chosen as a quantity that formalizes the user's perception of risk.
For the  expected value one obtains
\[
\Qev (v)\;:= \int_\Omega \lObjective[v,u[v](\omega)]\,\prob(d\omega).
\]
Risk measures are nonlinear quantities, which may contain for example deviations from targets or probabilities of critical events,
and that can be used to replace  expected value as the quantity optimized.
 Typical examples of risk measures 
are the excess probability $\Qep$, which is the probability of  exceeding a preselected target value $\costThreshold\in\R$
\[
\Qep(v)\;:= \prob\left[ \lObjective[v,u[v](\omega)]>\costThreshold\right],
\]
and the expected excess $\Qee$ given by the mean-value of the outcomes  above a preselected $\costThreshold\in\R$
\[
\Qee(v)\;:= \int_\Omega \max\left\{\lObjective[v,u[v](\omega)]\,- \,\costThreshold,\,0\right\}\,\prob(d\omega)
\]
In \cite{CoHePa11,GeLeRu13} risk neutral and risk averse stochastic  shape optimization were discussed for mean-risk models with the above specifications. To this end the stochastic energies were used as stochastic cost functions for shape optimization directly.

\section{Stochastic dominance constraints}\label{sec:dominance}

Let us at first present the general concept of stochastic dominance.
Introducing orders on families of random variables has  a long tradition in applied probability and statistics, see \cite{MuSt02,ShSh07}, also for proofs of facts mentioned  in  the discussion below.
Since we are dealing with minimization problems, preference of small outcomes over big ones is postulated in the present paper.
We stress that, at variance with  applications to finance where a preference for big outcomes is standard,
in our context being dominant means is being small(er). \\
Then, the definitions of dominance of first order and second order, respectively, are given as follows:\\[1ex]
(i) With (real-valued) random variables ${\sf X}$ and ${\sf Y}$, on some probability space $(\Omega,\mathcal{A},\prob)$, it is said that ${\sf X}$ is stochastically smaller 
in first order than  ${\sf Y}$, denoted ${\sf X}\leqone{\sf Y}$, if and only if
$\expectedValue h({\sf X})\le \expectedValue h({\sf Y})$ for all nondecreasing disutility functions $h$ for which both expectations exist.\\[1ex]
(ii) Furthermore, ${\sf X}$ is said to be smaller than  ${\sf Y}$ with respect to the increasing convex order, denoted ${\sf X}\leqtwo{\sf Y}$, if and only if  $\expectedValue h({\sf X})\;\le\; \expectedValue h({\sf Y})$
for all nondecreasing convex disutility functions $h$ for which both expectations exist.
\smallskip

For a rational decision maker who prefers less to more $\leqone$ and $\leqtwo$  correspond to the traditional first- and second-order dominance rules which assume preference of more to less. In view of this analogy we will refer to  first- and second-order stochastic dominance in the subsequent text.

The above definitions immediately show that first-order dominance implies second-order dominance, the reverse implication is invalid.  
The definitions of dominance relations via disutility functions are, however, not particularly well-suited for computations. Therefore it is  instructive to  consider the following equivalences (see, for instance, Subsection~8.1.2 in  \cite{MuSt02})
\begin{align}
 \label{eq:FirstOrderEquiv}
&{\sf X}\leqone{\sf Y}  & \mbox{ iff}  &\; \;\prob\left[\left\{\omega\in \Omega : {\sf X}(\omega)\le\costThreshold\right\}\right]\ge \prob[\{\omega \in\Omega : {\sf Y}(\omega)\le\eta\}] \;\mbox{ for all } \costThreshold\in\R\,.\\
\label{eq:SecOrderEquiv}
&{\sf X}\leqtwo{\sf Y}  &\mbox{ iff} &   \int_\Omega \max\left\{{\sf X(\omega)}-\costThreshold,\,0\right\}\,\prob(d\omega)\le \int_\Omega \max\left\{{\sf Y(\omega)}-\costThreshold,\,0\right\}\,\prob(d\omega)\mbox{ for all }\costThreshold \in\R\,.
\end{align}
In our shape optimization context, first order dominance can be understood as turning
the excess probability approach based on the cost function $\Qep$ into a constraint.
Analogously, second order dominance is conceptually related to the expected excess approach 
with the cost function $\Qee$.

Denoting by 
\[
F_X(t) := \prob\left[ {\sf X} \le t\right] \;\;\mbox{ and  }\;\; \pi_X(t):= \int_t^{+\infty}\left(1\,-\,F_X(z)\right)\, dz
\]
the cumulative  distribution function and the integrated survival function of  ${\sf X}$, respectively, the following equivalences readily follow from \eqref{eq:FirstOrderEquiv} and, after integration by parts, from \eqref{eq:SecOrderEquiv},
\begin{equation}
\label{eq:cdf}
 {\sf X}\leqone{\sf Y}\quad \mbox{ iff } \quad F_X(t)\,\ge F_Y(t)\;\mbox{ for all } \; t\in\R
\end{equation}
and 
\begin{equation}
\label{eq:isf}
 {\sf X}\leqtwo {\sf Y}\quad \mbox{ iff } \quad \pi_X(t)\,\le \pi_Y(t)\;\mbox{ for all } \; t\in\R.
\end{equation}
These relations can be interpreted as follows: First-order dominance comes as a pointwise,  rather strict, requirement, while second-order dominance tolerates violation of the pointwise relation at some $\tilde t$, but requires compensation for $t>\tilde t$ in terms of the integral defining $\pi_X$.

Coming back to shape optimization under uncertainty, having fixed  a benchmark  phase field 
$v_{b}$ (which describes an elastic object $\object_{b}$)
and a  random benchmark surface load $g(\omega_{b})$ one can formulate the abstract optimization problem 
\[
\min\{ \mathcal{G}(v)\;:\; \lObjective[v,\omega]\;\preceq\; \lObjective[v_{b},\omega_{b}],\;v \in H^{1,2}(\domain),\,\, v=1 \text{ on } \Gamma_N\}.
\]
Here $\preceq$ is specified as either $\leqone$ or $\leqtwo$,
and for the sake of simplicity we  identify $v_b$ with $\object_{b}$ 
and neglect  the difference between soft material ($\delta >0$) and void ($\delta =0$).

The entity $\lObjective[v_{b},\omega_{b}]$ corresponds to the compliance-valued random variable that arises when exposing the reference (or benchmark) shape 
$\object_b$ 
to the reference (or benchmark) random surface load $g(\omega_{b})$. In this model $\mathcal{G}$ stands for an objective function which, in the present context may involve the volume or the elastic energy of $\body$.
In the applications considered below we choose 
$$\mathcal{G}(v) := \V(v) + \epsilon \L(v)\,.$$

In what follows, we will assume a finite set of realizations of the stochastic loading,
each with different boundary data, in which every scenario is represented by a random variable $\omega_k$ and 
comes with a certain probability $\pi_k$ and a force $g_k$ 
for $k=1,\ldots, K$, with $\pi_k\in[0,1]$ and 
$\sum_{k=1}^K \pi_k =1$. 

Since the  set of realizations $\{\omega_k\}_{k=1,\ldots, K}$ 
of the stochastic load is finite, it suffices to consider a finite set of constraints.
For the random variable $X$ and random benchmark variable $Y$ 
we obtain in case of first order dominance: \\[2ex]
\centerline{$X \leqone Y$ iff $\prob[X \leq \eta] \geq \prob[Y \leq \eta]$ for 
$\eta =Y(\omega_j)$ with $j=1,\ldots, K$}  \\[2ex]
and in case of second order dominance:\\[2ex]
\centerline{$X \leqtwo Y$ iff $\mathbb{E}[\max\{X-\eta,0\}] \leq \mathbb{E}[\max\{Y-\eta,0\}]$ 
for all $\eta = Y(\omega_j)$ with $j=1,\ldots, K$.}\\[2ex]
We apply this now in our space discretized shape optimization context to
$X=\lObjective[v,u[v](\omega)]$ and $Y=\lObjective[v_{b},u[v_{b}](\omega)]$.
For notational uniformity we write
$$\prob (\lObjective[v,u[v](\omega)]\leq\eta)=\sum_{k=1}^K 
\pi_k \hv{\eta-\lObjective[v,u[v](\omega_k)]}$$ 
where $\hv{\cdot}$ is the Heaviside
function, defined by $\hv{x}= 1$ if $x \geq 0$ and $\hv{x}= 0$ otherwise.
We obtain the following sets of constraints:\\[2ex]
\emph{First order dominance},
$\lObjective[v,u[v](\omega)] \leqone \lObjective[v_{b},u[v](\omega)]$, is equivalent to
$$
\sum_{k=1}^K \pi_k \hv{\lObjective[v_{b},u[v_{b}](\omega_j)])-\lObjective[v,u[v](\omega_k)]}
\geq 
\sum_{k=1}^K \pi_k \hv{\lObjective[v_{b},u[v_{b}](\omega_j)])-
\lObjective[v_{b},u[v_{b}](\omega_k)]}
$$
for all $j=1,\ldots, K$.\\[2ex]
 \emph{Second order dominance},
$\lObjective[v,u[v](\omega)] \leqtwo \lObjective[v_{b},u[v](\omega)]$, is equivalent to
$$\sum_{k=1}^K \pi_k 
\max\{\lObjective[v,u[v](\omega_k)]-\lObjective[v_{b},u[v_{b}](\omega_j)],0\} \leq
\sum_{k=1}^K \pi_k 
\max\{\lObjective[v_{b},u[v_{b}(\omega_k)]-\lObjective[v_{b},u[v_{b}](\omega_j)],0\}
$$
for all $j=1,\ldots, K$.\\[2ex]

\section{Spatial discretization}\label{sec:discrete}
We use adaptive Finite Elements to discretize the optimization problem.
To this end, we consider an adaptive mesh of rectangular elements covering the computational domain $\domain$. 
The grid is handled using a quad-tree as the underlying hierarchical data structure. We suppose that there is
a transition of at most one level at element faces, so that there is at most one hanging node on each edge.
On this mesh we define the space $\FE_h$ of piecewise bilinear, continuous functions,
where $h$ denotes the piecewise constant mesh size function.
Values of function in $\FE_h$ at hanging nodes are obtained by interpolation.
For a discrete phase field function $V \in \FE_h$ we define the discrete counterpart 
of the Modica--Mortola type energy functional 
\begin{equation*}
\mathbf{L}^\varepsilon[V]:=\frac12\int_\Omega\varepsilon|\nabla V|^2+\frac1\varepsilon \I_h(\Psi(V))\,\d x\,,
\end{equation*}
where $\I_h$ denotes the bilinear Lagrangian interpolation on the mesh.
The discrete volume functional is given as  $\mathbf{V}(V):=\int_\domain\charfunction_\object(V)\,\d x$ for any $V\in \FE_h$
and the discrete elastic energy for a phase field $V\in \FE_h$ and a discrete displacement $U\in \FE_h^2$ is
\begin{equation}
\mathbf{W}^\delta [V, U] := \frac{1}{2} \int_\domain ((1 - \delta)
\I_h(\charfunction_\object(V)) + \delta) \, \frac{C}{4} (DU^T+DU):   (DU^T+DU) \d x\,.
\end{equation}
Furthermore, we suppose that the Dirichlet boundary $\Gamma_D$ and the inhomogeneous Neumann boundary  $\Gamma_N$
are resolved on the adaptive rectangular mesh. We define the space $\FE_{h,0}$ as the subset of $\FE_{h}$ with vanishing trace of $\Gamma_D$ and require that $V|_{\Gamma_N}=1$ for discrete phase field functions $V$.
The discrete compliance functional is given by
\begin{equation*}
\mathbf{C}[U]:=\int_{\Gamma_N} \I_h(g\cdot U) \,\d a\,.
\end{equation*}
Thus we obtain the discrete total free energy
\begin{equation*}
\mathbf{E}[V,U]:=\mathbf{W}[V,U]-\mathbf{C}[U]\,.
\end{equation*}
Finally the discrete cost functional is given by
\begin{equation}
\label{eq:discrJ}
\mathbf{J}[V, U] := 2 \mathbf{W}^\delta[V,U] + \nu \mathbf{V}[V] + \eta \mathbf{L}^\varepsilon[V]\,.
\end{equation} 
The deterministic, discrete shape optimization problem is to minimize $\mathbf{J}[V, U[V]]$ subject to the constraint that 
$U[V]$ is a minimizer of $U\mapsto \mathbf{E}[V,U]$.
Due to the finite set of realizations $\{\omega_i\}_{i=1,\ldots, K}$ 
of the stochastic load it suffices to consider a finite set of constraints,
and the stochastic dominance  conditions reduce to the finitely many conditions stated in the last
two equations of the previous section.
We apply them in our space discretized shape optimization context to
$X=\lObjective[V,U[V](\omega)]$ and $Y=\lObjective[V_{b},U[V_{b}](\omega)]$,
where $V_{b}\in \FE_h$ is a given discrete benchmark phase field describing the benchmark shape.
To ensure differentiability of the cost in the numerical descent method
we consider smooth approximations of the Heavyside function and the $\max$ function. Specifically,
fixing a small regularization parameter $\gamma>0$ we replace $\hv{x}$ by
\begin{equation*}
  \hvs{x} := \frac{1}{1 + \exp{(-2 \gamma x)}}
\end{equation*}
and $\max\{x,0\}$ by
\begin{equation*}
  \maxs \left\lbrace x, 0 \right\rbrace := \frac{\sqrt{x^2 + \gamma} + x}{2}
\end{equation*}
After this regularization we obtain the following sets of constraints:\\[2ex]
The \emph{first order dominance} condition 
$\lObjective[v,u[v](\omega)] \leqone \lObjective[v_{b},u[v](\omega)]$
is numerically approximated by
\beq
&&\sum_{k=1}^K \pi_k \hvs{\mathbf{J}[V_{b},U[V_{b}](\omega_j)])-\mathbf{J}[V,U[V](\omega_k)]} \\[-2ex]
&& \qquad \qquad \qquad \qquad \geq 
\sum_{k=1}^K \pi_k \hvs{\mathbf{J}[V_{b},U[V_{b}](\omega_j)])-\mathbf{J}[V_{b},U[V_{b}](\omega_k)]}
\eeq
for all $j=1,\ldots, K$.\\[2ex]
The \emph{second order dominance} condition 
$\lObjective[v,u[v](\omega)] \leqtwo \lObjective[v_{b},u[v](\omega)]$
is numerically approximated by
\beq
&&\sum_{k=1}^K \pi_k 
\maxs\{\mathbf{J}[V,U[V](\omega_k)]-\mathbf{J}[V_{b},U[V_{b}](\omega_j)],0\}\\[-2ex] 
&& \qquad \qquad \qquad \qquad \leq 
\sum_{k=1}^K \pi_k 
\maxs\{\mathbf{J}[V_{b},U[V_{b}](\omega_k)]-\mathbf{J}[V_{b},U[V_{b}](\omega_j)],0\}
\eeq
for all $j=1,\ldots, K$.\\[2ex]
Finally, as the discrete counterpart of  $\mathcal{G}(\cdot)$
we choose
$$
\mathbf{G}(V) := \mathbf{V}(V) + \epsilon \mathbf{L}(V)\,.
$$
The constraint minimization problem will be solved using the {\tt IPOPT}
\cite{WaBi06} package.
To this end, we have to provide derivatives of $\mathbf{G}(V)$ and
the constraints. Let $\delta$ denote the Gateaux-derivative.
The derivative of $\mathbf{G}(V)$ can be computed directly as
$\delta_V \mathbf{G}(V)(\Theta) = \delta_V \mathbf{V}(V)(\Theta) 
+ \epsilon \delta_V \mathbf{L}(V)(\Theta)$ with
\begin{align}
  \label{eq:DV}
  \delta_V \mathbf{V}(V)(\Theta) &= \int_\Omega \partial_V \charfunction_\object(V) \Theta\quad
  \text{ and}\\
  \label{eq:DL}
  \delta_V \mathbf{L}(V)(\Theta) &=
  \frac{\epsilon}{2} \int_\Omega \varepsilon \nabla V \nabla \Theta
  + \frac1\varepsilon \I_h(\partial_V \Psi(V) \Theta)\,\d x\,.
\end{align}
For the derivative $\delta_V (\mathbf{J}[V, U[V]]) = (\partial_V \mathbf{J})[V, U[V]]
+ (\partial_U \mathbf{J})[V, U[V]](\partial_V U[V])$ it is convenient to use the dual formulation. We define
$P := (\delta^2_{UU} \mathbf{E}[V, U[V]]^{-1})(\delta_U \mathbf{J}[V, U[V]])$ and consider
\begin{align*}
  \delta_V &(\mathbf{J}[V, U[V]](S)) = (\partial_V \mathbf{J})[V, U[V]](S)
  + \partial^2_{UU} \mathbf{E}[V, U[V]](P)(\partial_V U[V](S))\\
  &= (\partial_V \mathbf{J})[V, U[V]](S)
  - \partial^2_{UU} \mathbf{E}[V, U[V]](P)((\partial^2_{UU} \mathbf{E}^{-1}[V, U[V]])
  (\partial^2_{UV} \mathbf{E}[V, U[V]])(S))\\
  &= (\partial_V \mathbf{J})[V, U[V]](S) - \partial^2_{UV} \mathbf{E}[V, U[V]])(P)(S)\,.
\end{align*}
The second equation makes uses of the implicit function theorem which states that
$$\partial_V U[V] = -(\partial^2_{UU} \mathbf{E}^{-1}[V, U[V]])
  (\partial^2_{UV} \mathbf{E}[V, U[V]])$$ in a neighbourhood of $V$.

For the energies $\mathbf{J}[V, U[V]]$ and $\mathbf{E}[V, U]$ described here, the dual
problem can be solved explicitly as $P = 2U$. Thus, the derivative of
$\mathbf{J}[V, U[V]]$ can be written as
\begin{equation}
  \delta_V (\mathbf{J}[V, U[V]]) = -\delta_V \mathbf{W}^\delta[V,U] +
  \nu  \delta_V \mathbf{V}[V] + \eta  \delta_V \mathbf{L}^\varepsilon[V]
\end{equation}
with $\delta_V \mathbf{W}^\delta[V,U](\Theta) = \int_\domain ((1 - \gamma)
\I_h(\charfunction_\object(V)) + \delta) \, \frac{C}{4} (DU^T+DU):(D\Theta^T+D\Theta) \d x$
and $\delta_V \mathbf{V}[V]$ and $\delta_V \mathbf{L}^\varepsilon[V]$ as given in (\ref{eq:DV})
and (\ref{eq:DL}).

We use adaptive grid refinement 
after each run of the optimization algorithm.
Specifically, after each optimization run elements where $\nabla V > t$ for a $t > 0$ and their
top, bottom, left and right neighbors are marked for grid refinement.
The solution is then prolongated to the new grid and
used to restart the optimization algorithm on the finer grid with a smaller parameter
$\varepsilon$ for the interfacial region. 
Each run of {\tt IPOPT} stops when the error estimate $E_0(x, \lambda, z)$ as given in 
\cite[eq. (5)]{WaBi06} becomes smaller than a given tolerance $\epsilon_\text{tol}$, which
makes the complete algorithm stop when $E_0(x, \lambda, z) < \epsilon_\text{tol}$ at
a grid level supporting a sufficiently small $\epsilon$.
Furthermore, it turned out to be advantageous to start with a
relatively large parameter $\heavreg$ (the heaviside/max regularization
parameter) and to successively decrease that parameter during the course of
the optimization and refinement algorithm described above. 

\section{Numerical experiments}\label{sec:numerics}
We investigated two different optimization tasks. In both applications the working domain is $D=[0,1]^2$.
All results presented in this chapter have been computed with paramters $\lambda = \mu = 80$
for the elaticity tensor, a factor $\delta = 10^{-4}$ for the soft material and weights  $\nu = 0.04096$
and $\eta = 0.00064$ for the cost functional $\mathbf{J}$ given in (\ref{eq:discrJ}).

\paragraph{Cantilever.}
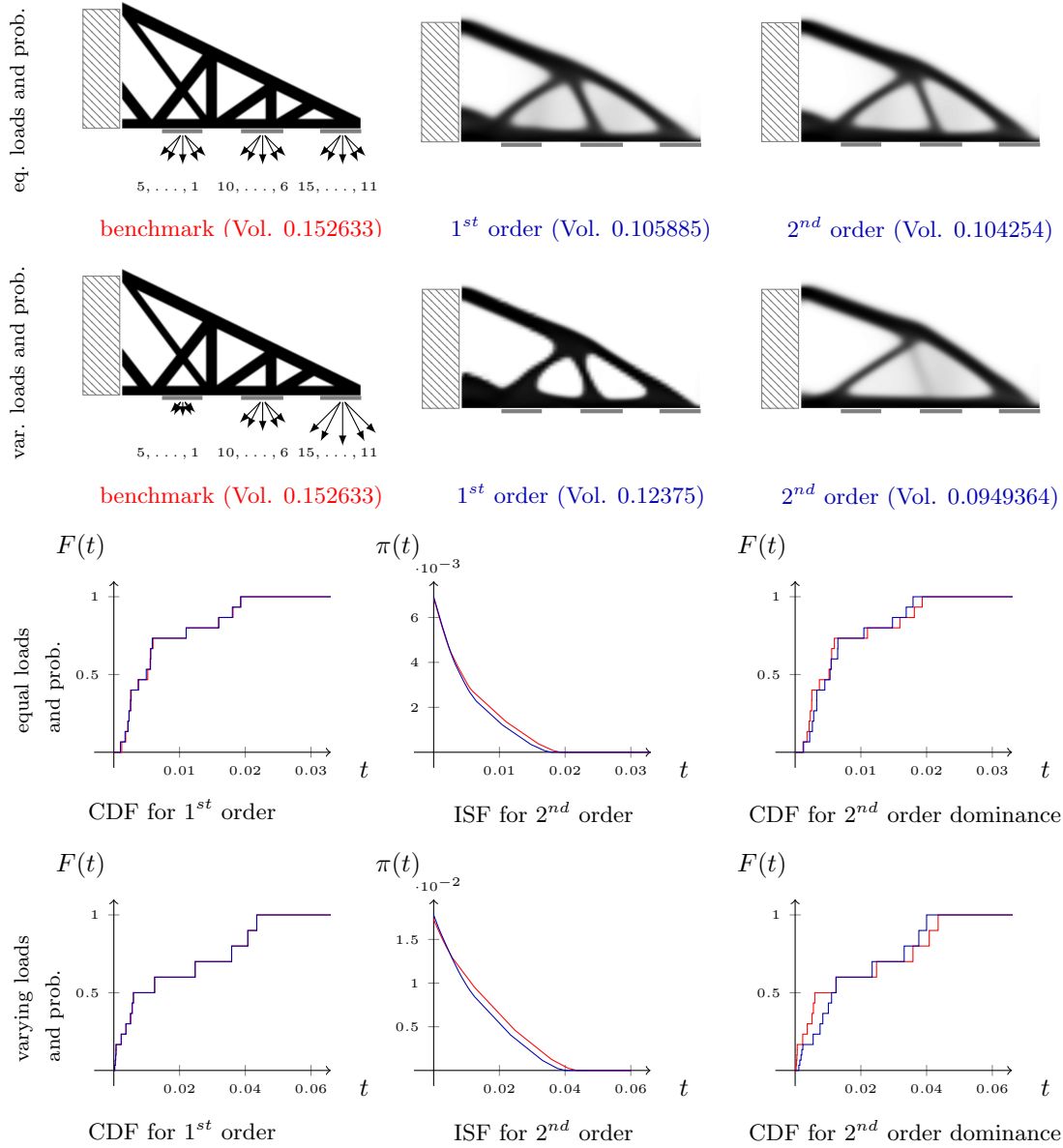
\begin{figure}[h!]
  \centering
  \begin{subfigure}[t]{0.06\textwidth}
    \begin{tikzpicture}
      \node[rotate=90,align=left] at (0.0,0.0) {\footnotesize eq. loads and prob.};
    \end{tikzpicture}
  \end{subfigure}
  \begin{subfigure}[t]{0.3\textwidth}
    \includegraphics[width=\textwidth]{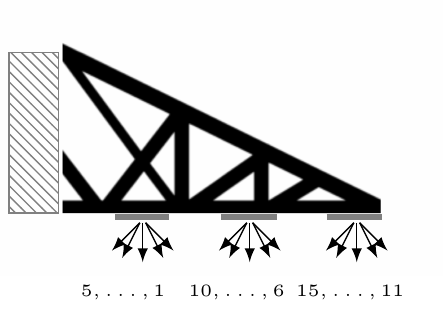}
    \caption{\color{red} benchmark (Vol. $0.152633$)}
    \label{fig:benchcanti3e}
  \end{subfigure}
  \begin{subfigure}[t]{0.3\textwidth}
    \includegraphics[width=\textwidth]{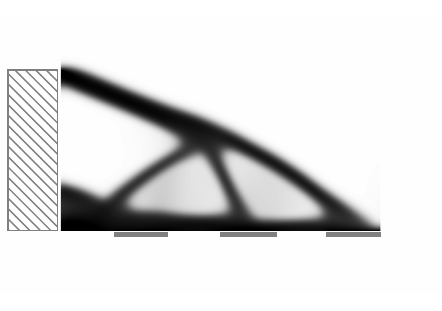}
    \caption{\color{blue!70!black} $1^{st}$ order (Vol. $0.105885$)}
    \label{fig:res1canti3e}
  \end{subfigure}
  \begin{subfigure}[t]{0.3\textwidth}
    \includegraphics[width=\textwidth]{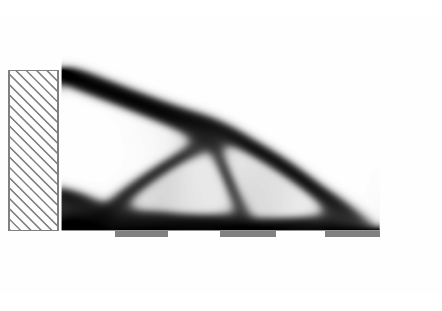}
    \caption{\color{blue!70!black} $2^{nd}$ order (Vol. $0.104254$)}
    \label{fig:res2canti3e}
  \end{subfigure}\\[-1ex]
   \begin{subfigure}[t]{0.06\textwidth}
    \begin{tikzpicture}
      \node[rotate=90,align=left] at (0.0,0.0) {\footnotesize var. loads and prob.};
    \end{tikzpicture}
  \end{subfigure}
  \begin{subfigure}[t]{0.3\textwidth}
    \includegraphics[width=\textwidth]{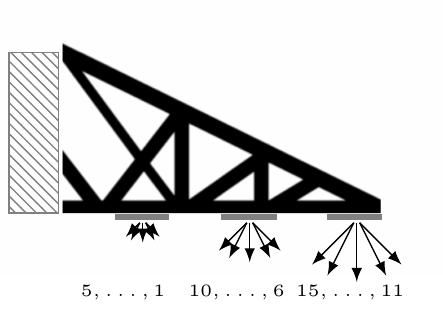}
    \caption{\color{red} benchmark (Vol. $0.152633$)}
    \label{fig:benchcanti3w}
  \end{subfigure}
  \begin{subfigure}[t]{0.3\textwidth}
    \includegraphics[width=\textwidth]{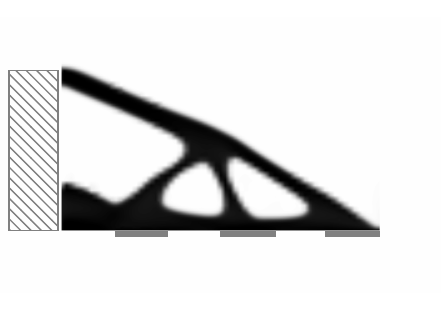}
    \caption{\color{blue!70!black} $1^{st}$ order (Vol. $0.12375$)}
    \label{fig:res1canti3w}
  \end{subfigure}
  \begin{subfigure}[t]{0.3\textwidth}
    \includegraphics[width=\textwidth]{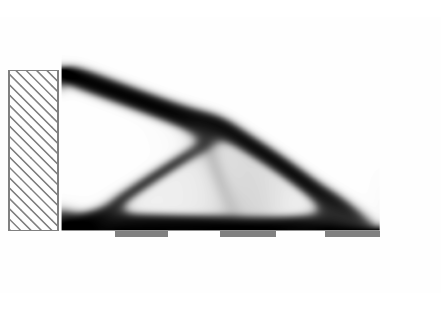}
    \caption{\color{blue!70!black} $2^{nd}$ order (Vol. $0.0949364$)}
    \label{fig:res2canti3w}
  \end{subfigure}\\[1ex]
  \begin{subfigure}[t]{0.32\textwidth}
    \begin{tikzpicture}
      \begin{axis}[axis lines=middle,
        axis line style={->},
        enlargelimits=true,
        width=\textwidth,
        xmin=0, xmax=0.03,
        xlabel={$t$},
        x label style={at={(axis description cs:1.2, -0.1)}},
        xticklabel style={
          /pgf/number format/fixed,
          /pgf/number format/precision=4,
          font=\tiny },
        scaled x ticks=false,
        ylabel={$F(t)$},
        y label style={
          at={(axis description cs:-0.2, 1.3)},
          rotate=0 },
        yticklabel style={
          font=\tiny },
        ]
        \addplot[red]   table[x={x}, y={F}, mark=none] {images/res/canti3e/bench_cdf1};
        \addplot[blue!70!black] table[x={x}, y={F}, mark=none] {images/res/canti3e/pf_cdf1};
      \end{axis}
      \node[rotate=90,align=left] at (-0.75,1.2) {\footnotesize equal loads\\ \footnotesize and prob.};
    \end{tikzpicture}
    \vspace*{-3ex}
    \caption{CDF for $1^{st}$ order }
    \label{fig:cdf1canti3e} 
  \end{subfigure}
  \begin{subfigure}[t]{0.32\textwidth}
    \begin{tikzpicture}
      \begin{axis}[axis lines=middle,
        axis line style={->},
        enlargelimits=true,
        width=\textwidth,
        xmin=0, xmax=0.03,
        xlabel={$t$},
        x label style={at={(axis description cs:1.2, -0.1)}},
        xticklabel style={
          /pgf/number format/fixed,
          /pgf/number format/precision=4,
          font=\tiny },
        scaled x ticks=false,
        ylabel={$\pi(t)$},
        y label style={
          at={(axis description cs:-0.2, 1.3)},
          rotate=0 },
        yticklabel style={
          font=\tiny },
        ]
        \addplot[red]   table[x={x}, y={isf}, mark=none] {images/res/canti3e/bench_isf2};
        \addplot[blue!70!black] table[x={x}, y={isf}, mark=none] {images/res/canti3e/pf_isf2};
      \end{axis}
    \end{tikzpicture}
    \caption{ISF for $2^{nd}$ order }
    \label{fig:isf2canti3e}
  \end{subfigure}
  \begin{subfigure}[t]{0.32\textwidth}
    \begin{tikzpicture}
      \begin{axis}[axis lines=middle,
        axis line style={->},
        enlargelimits=true,
        width=\textwidth,
        xmin=0, xmax=0.03,
        xlabel={$t$},
        x label style={at={(axis description cs:1.2, -0.1)}},
        xticklabel style={
          /pgf/number format/fixed,
          /pgf/number format/precision=4,
          font=\tiny },
        scaled x ticks=false,
        ylabel={$F(t)$},
        y label style={
          at={(axis description cs:-0.2, 1.3)},
          rotate=0 },
        yticklabel style={
          font=\tiny },
        ]
        \addplot[red]   table[x={x}, y={F}, mark=none] {images/res/canti3e/bench_cdf2};
                \addplot[blue!70!black] table[x={x}, y={F}, mark=none] {images/res/canti3e/pf_cdf2};
              \end{axis}
    \end{tikzpicture}
    \caption{CDF for $2^{nd}$ order dominance}
    \label{fig:cdf2canti3w}
  \end{subfigure}\\[1ex]
  \begin{subfigure}[t]{0.32\textwidth}
    \begin{tikzpicture}
      \begin{axis}[axis lines=middle,
        axis line style={->},
        enlargelimits=true,
        width=\textwidth,
        xmin=0, xmax=0.06,
        xlabel={$t$},
        x label style={at={(axis description cs:1.2, -0.1)}},
        xticklabel style={
          /pgf/number format/fixed,
          /pgf/number format/precision=4,
          font=\tiny },
        scaled x ticks=false,
        ylabel={$F(t)$},
        y label style={
          at={(axis description cs:-0.2, 1.3)},
          rotate=0 },
        yticklabel style={
          font=\tiny },
        ]
        \addplot[red]   table[x={x}, y={F}, mark=none] {images/res/canti3w/bench_cdf1};
        \addplot[blue!70!black] table[x={x}, y={F}, mark=none] {images/res/canti3w/pf_cdf1};
      \end{axis}
      \node[rotate=90,align=left] at (-0.75,1.2) {\footnotesize varying loads \\
        \footnotesize and prob.};
    \end{tikzpicture}
    \caption{CDF for $1^{st}$ order }
    \label{fig:cdf1canti3w}
  \end{subfigure}
  \begin{subfigure}[t]{0.32\textwidth}
    \begin{tikzpicture}
      \begin{axis}[axis lines=middle,
        axis line style={->},
        enlargelimits=true,
        width=\textwidth,
        xmin=0, xmax=0.06,
        xlabel={$t$},
        x label style={at={(axis description cs:1.2, -0.1)}},
        xticklabel style={
          /pgf/number format/fixed,
          /pgf/number format/precision=4,
          font=\tiny },
        scaled x ticks=false,
        ylabel={$\pi(t)$},
        y label style={
          at={(axis description cs:-0.2, 1.3)},
          rotate=0 },
        yticklabel style={
          font=\tiny },
        ]
        \addplot[red]   table[x={x}, y={isf}, mark=none] {images/res/canti3w/bench_isf2};
        \addplot[blue!70!black] table[x={x}, y={isf}, mark=none] {images/res/canti3w/pf_isf2};
      \end{axis}
    \end{tikzpicture}
    \caption{ISF for $2^{nd}$ order }
    \label{fig:isf2canti3w}
  \end{subfigure}
  \begin{subfigure}[t]{0.32\textwidth}
    \begin{tikzpicture}
      \begin{axis}[axis lines=middle,
        axis line style={->},
        enlargelimits=true,
        width=\textwidth,
        xmin=0, xmax=0.06,
        xlabel={$t$},
        x label style={at={(axis description cs:1.2, -0.1)}},
        xticklabel style={
          /pgf/number format/fixed,
          /pgf/number format/precision=4,
          font=\tiny },
        scaled x ticks=false,
        ylabel={$F(t)$},
        y label style={
          at={(axis description cs:-0.2, 1.3)},
          rotate=0 },
        yticklabel style={
          font=\tiny },
        ]
        \addplot[red]   table[x={x}, y={F}, mark=none] {images/res/canti3w/bench_cdf2};
        \addplot[blue!70!black] table[x={x}, y={F}, mark=none] {images/res/canti3w/pf_cdf2};
      \end{axis}
    \end{tikzpicture}
    \caption{CDF for $2^{nd}$ order dominance}
    \label{fig:cdf2canti3w}
  \end{subfigure}
  \caption{Benchmark with load configuration plot, optimal shapes cumulative distance funcions (CDF)
    and integrated survival functions (ISF) for the cantilever setup.}
  \label{fig:rescanti}
\end{figure}

As a first application, we discuss the shape optimization of a 2D cantilever. 
The geometry is sketched in the first panel of Figure \ref{fig:rescanti}.
The cantilever is fixed on the left-hand side, where  the phase field $v$ is required to be 1
and the elastic displacement $u$ obeys  homogeneous Dirichlet boundary conditions.
We select three segments in the lower boundary on which piecewise constant surface loads are applied, 
around each of them we prescribe $v=1$ in a thin strip.
In each scenario a single load is applied on one of these segments. We compare different directions of the loading, different absolute values of the load 
and different probabilities.
Two different stochastic loading configurations are depicted in the first and second
row of Fig. \ref{fig:rescanti}, respectively. 

For the initial configuration, a phase field regularization parameter of $\epsilon = 0.025$ was used. After each
grid refinement, $\epsilon$ was multiplied by $0.75$. The smallest $\epsilon$ was set to $5.93 \cdot 10^{-3}$ for the first order varying load and probability configuration and to $7.91 \cdot 10^{-3}$ for all others.
The regularization parameter
$\heavreg$ was $1$ initially and decreased during the simulation,  
its smallest value was $1.95 \cdot 10^{-3}$.
Fig. \ref{fig:rescanti} shows results for first order dominance and second order dominance calculations based on the scheme derived in the previous section.

\begin{figure}[h!]
  \begin{subfigure}[t]{0.05\textwidth}
    \begin{tikzpicture}
      \node[rotate=90] at (0.0,0.0) {\footnotesize $1^{st}$ order};
    \end{tikzpicture}
  \end{subfigure}
  \begin{subfigure}[t]{0.11\textwidth}
    \includegraphics[width=\textwidth]{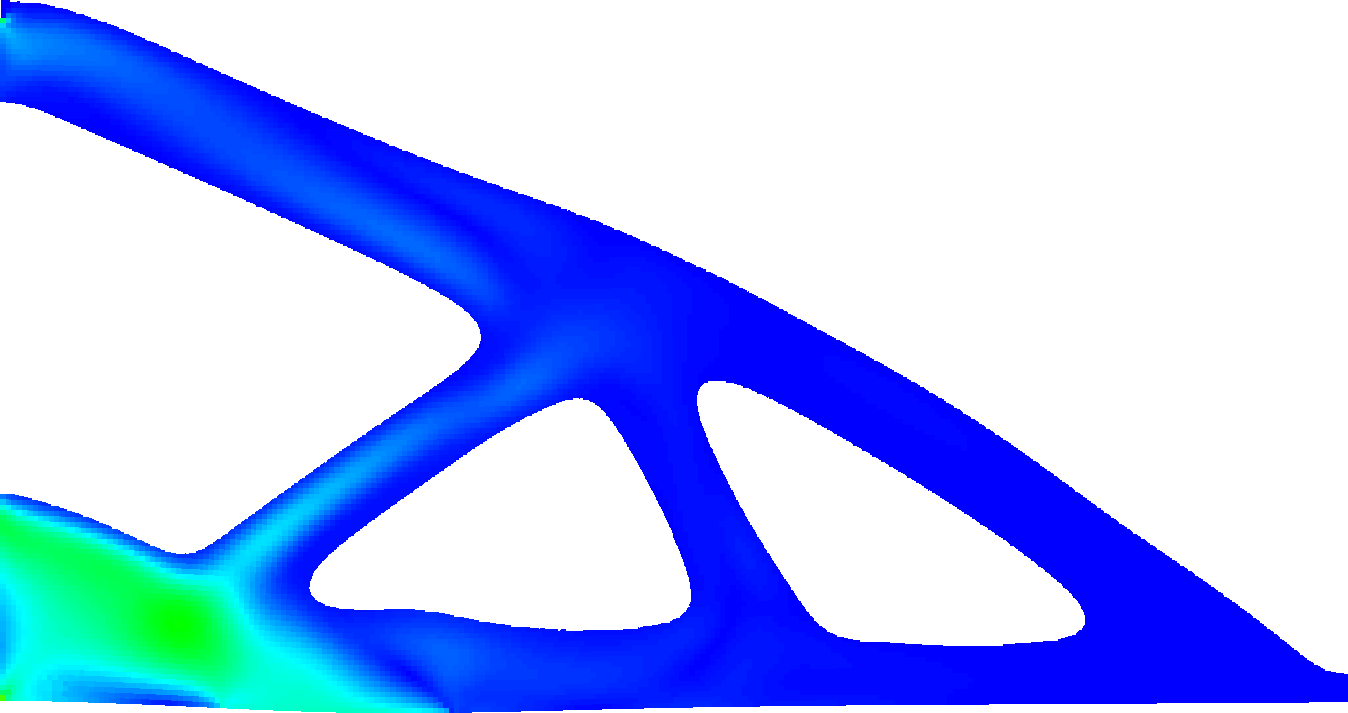}
  \end{subfigure}
  \begin{subfigure}[t]{0.11\textwidth}
    \includegraphics[width=\textwidth]{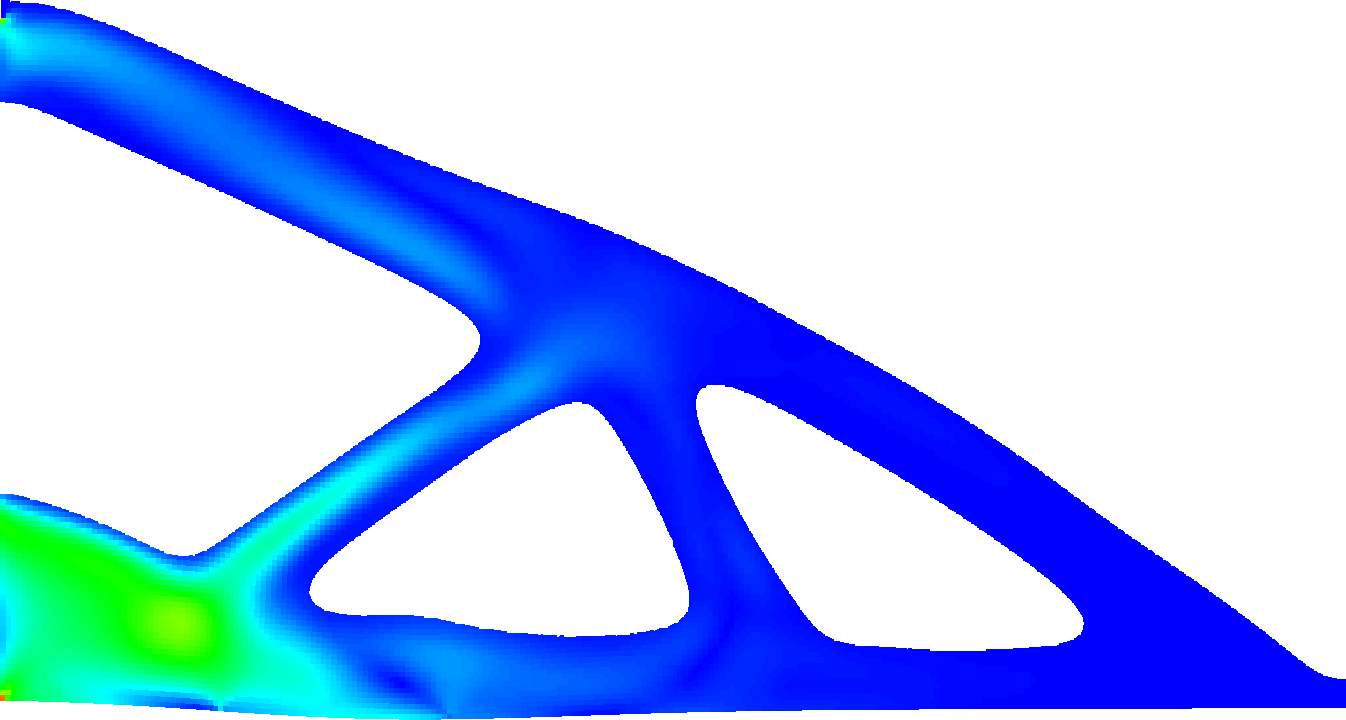}
  \end{subfigure}
  \begin{subfigure}[t]{0.11\textwidth}
    \includegraphics[width=\textwidth]{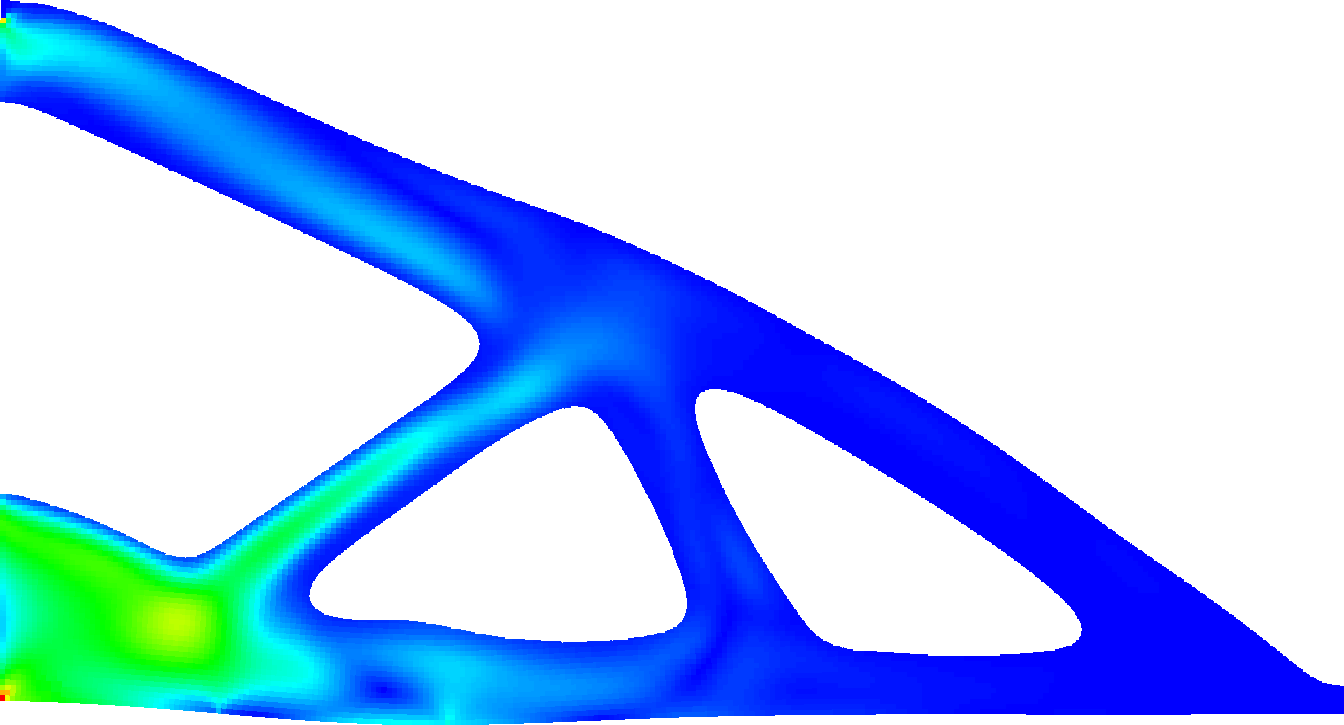}
  \end{subfigure}
  \begin{subfigure}[t]{0.11\textwidth}
    \includegraphics[width=\textwidth]{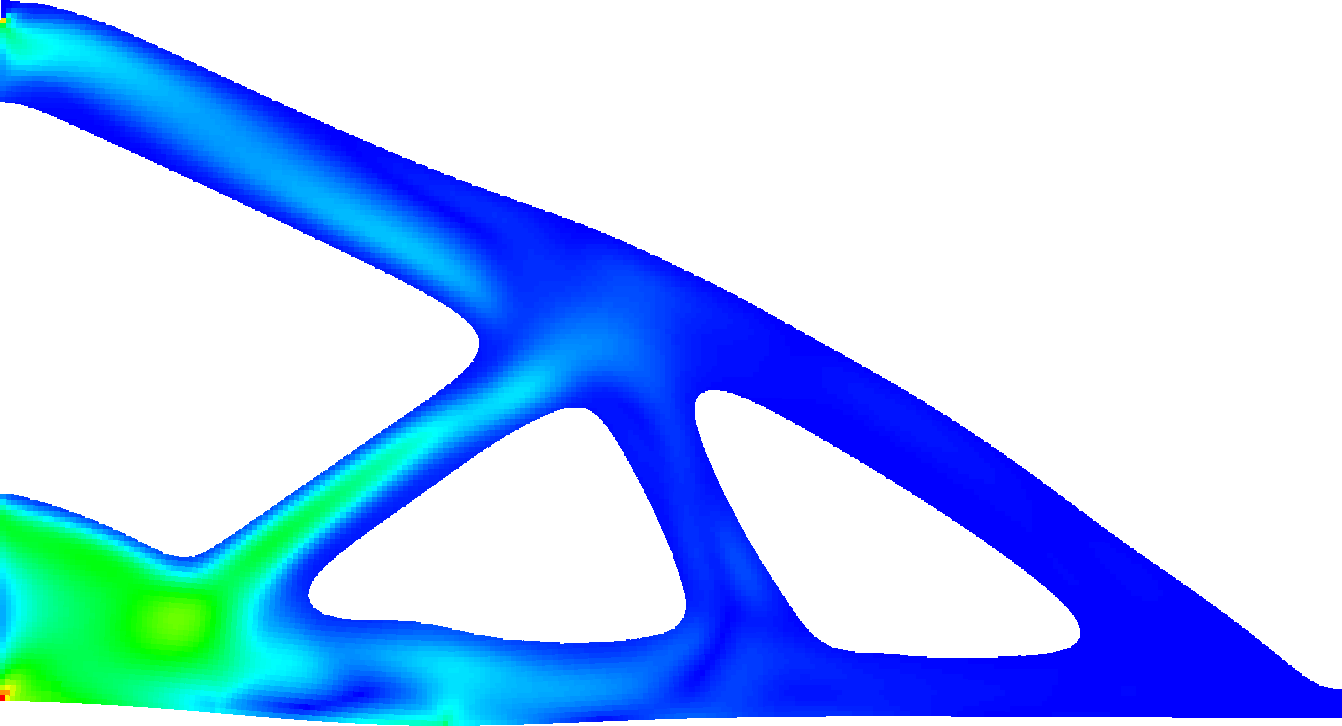}
  \end{subfigure}
  \begin{subfigure}[t]{0.11\textwidth}
    \includegraphics[width=\textwidth]{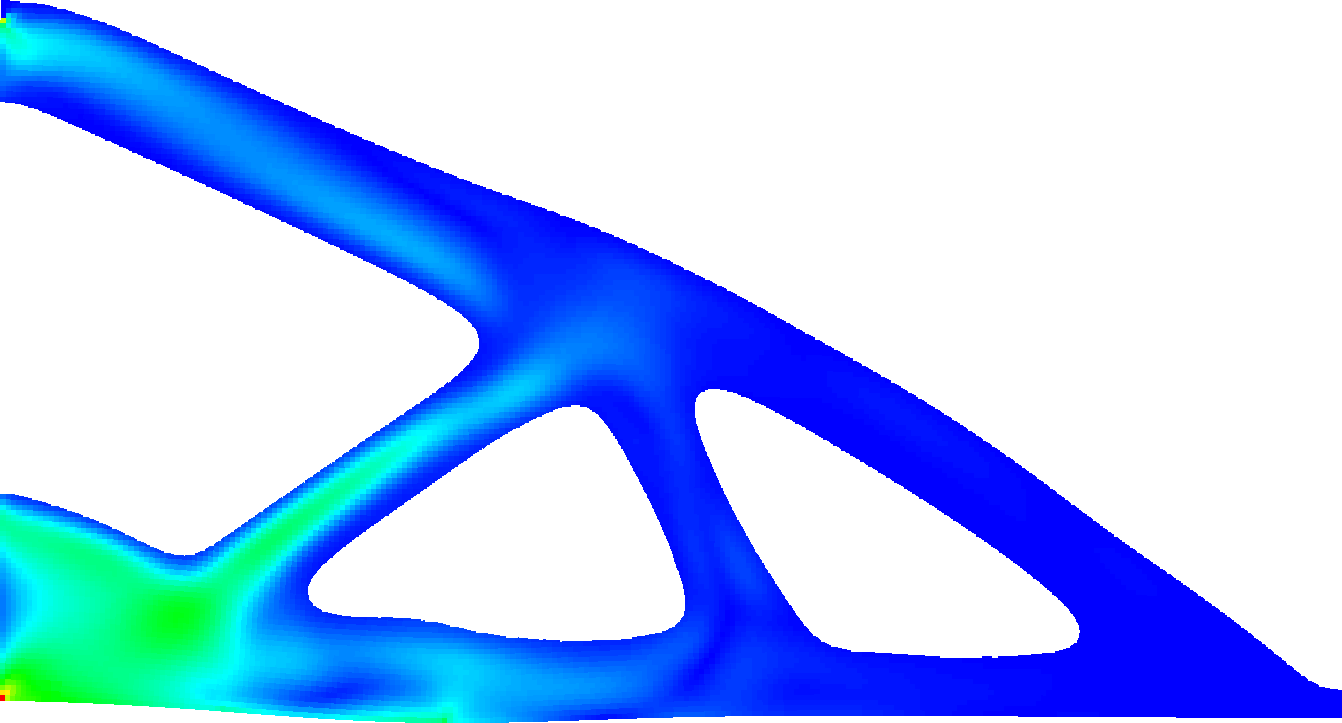}
  \end{subfigure}
  \begin{subfigure}[t]{0.11\textwidth}
    \includegraphics[width=\textwidth]{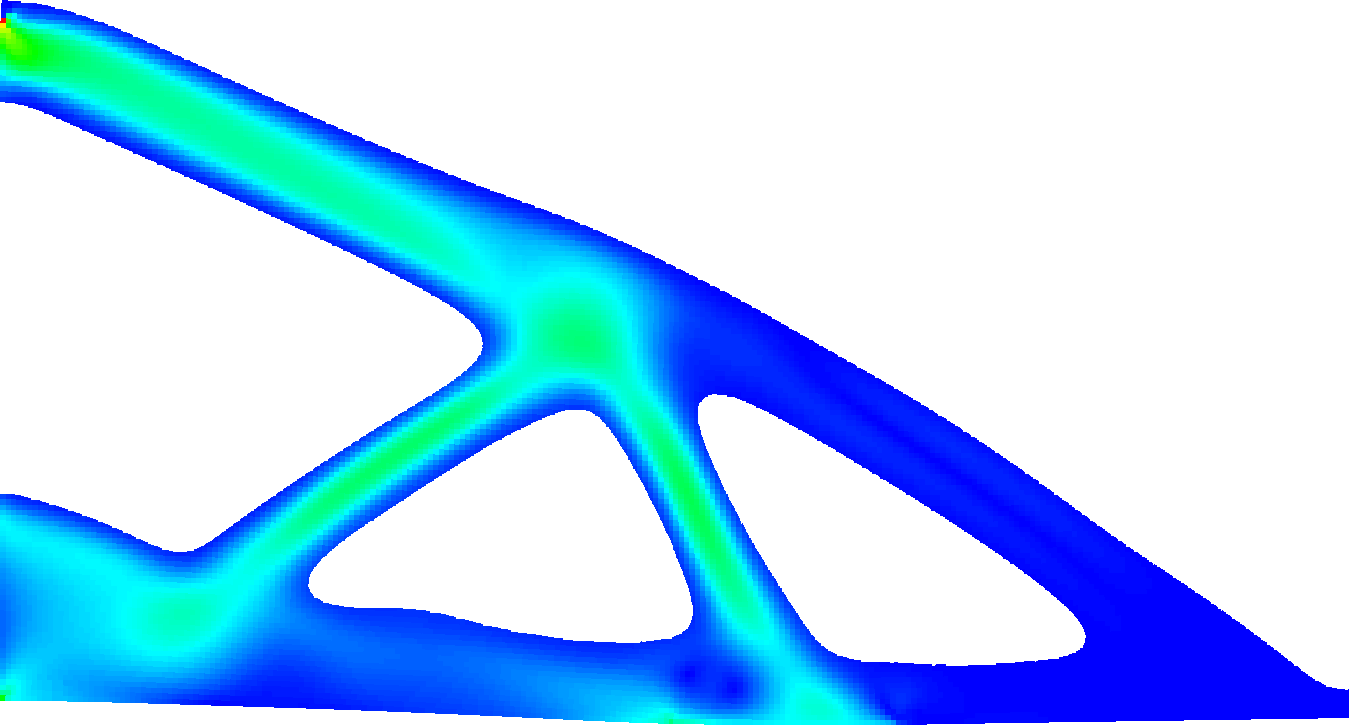}
  \end{subfigure}
  \begin{subfigure}[t]{0.11\textwidth}
    \includegraphics[width=\textwidth]{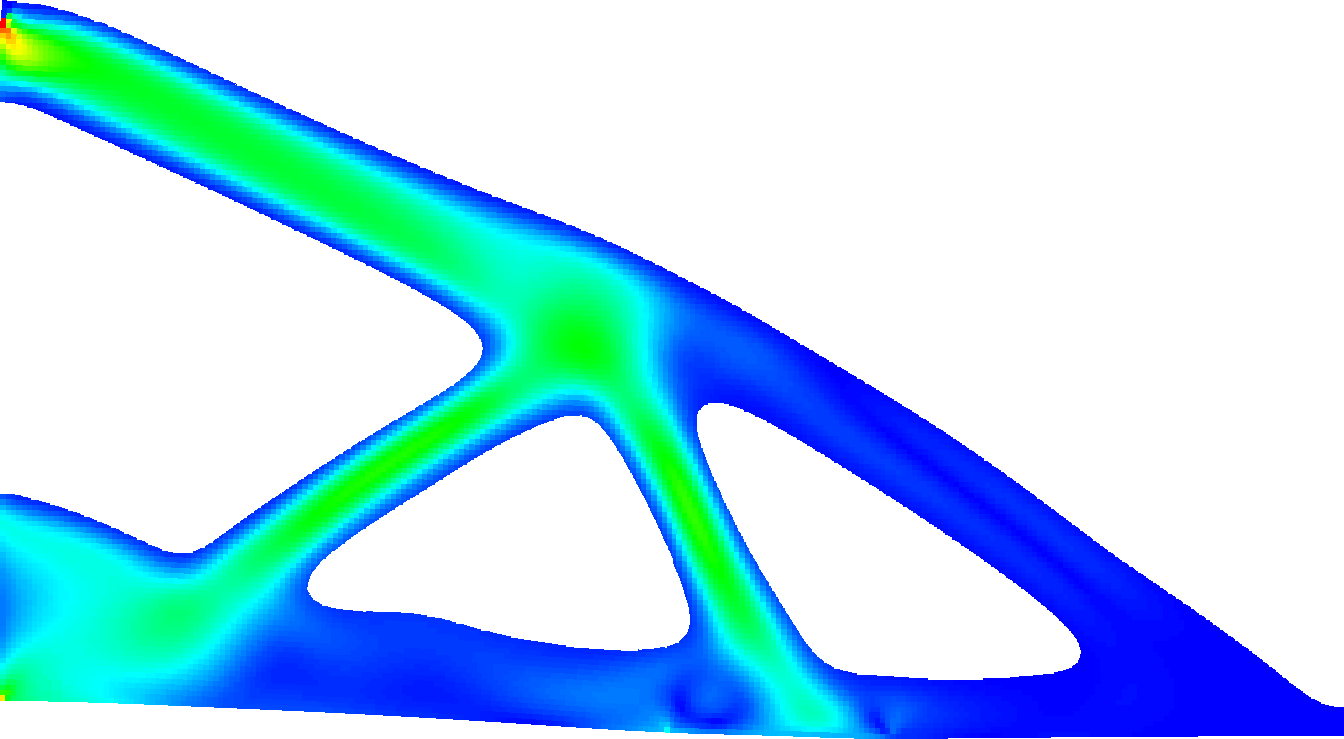}
  \end{subfigure}
  \begin{subfigure}[t]{0.11\textwidth}
    \includegraphics[width=\textwidth]{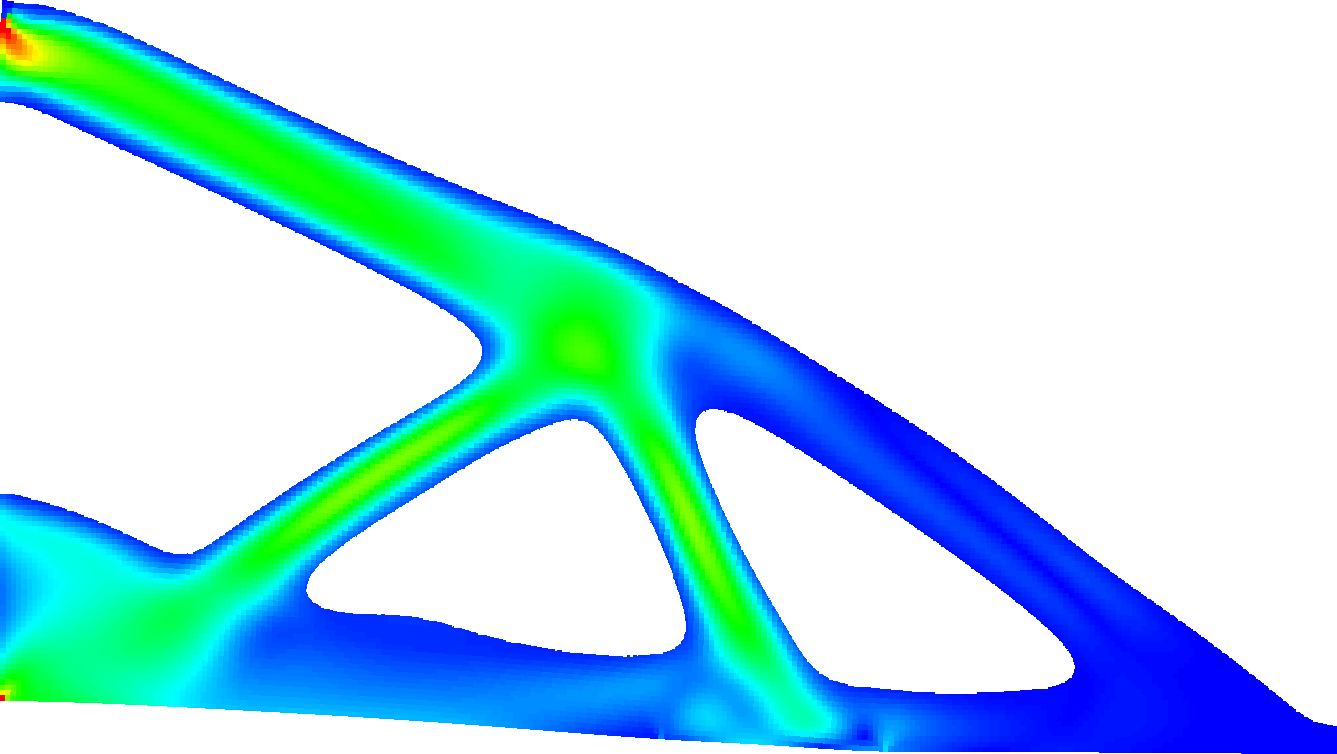}
  \end{subfigure}\\[1ex]
  \begin{subfigure}[t]{0.06\textwidth}
    \begin{tikzpicture}
      \node[rotate=90,align=left] at (0.0,0.0) {\footnotesize eq. loads\\
        \footnotesize and prob.};
    \end{tikzpicture}
  \end{subfigure}
  \begin{subfigure}[t]{0.11\textwidth}
    \includegraphics[width=\textwidth]{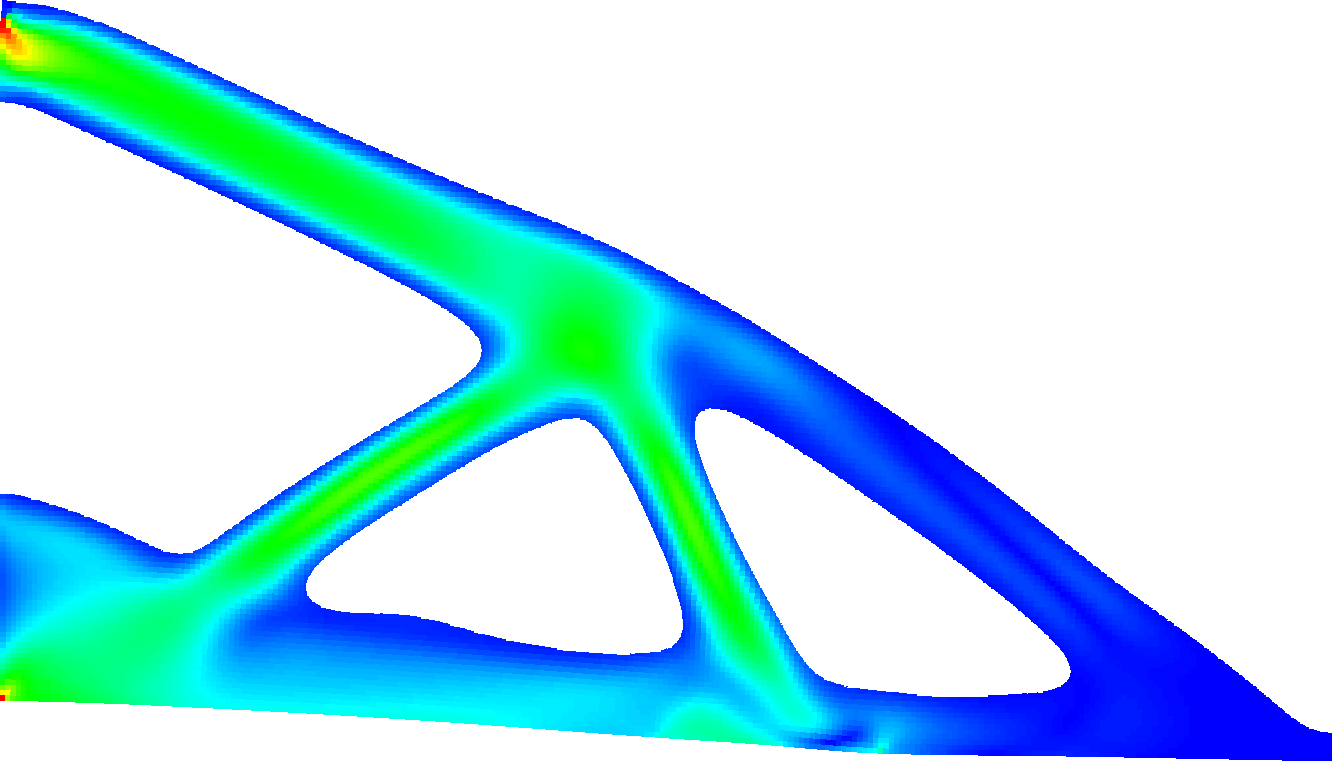}
  \end{subfigure}
  \begin{subfigure}[t]{0.11\textwidth}
    \includegraphics[width=\textwidth]{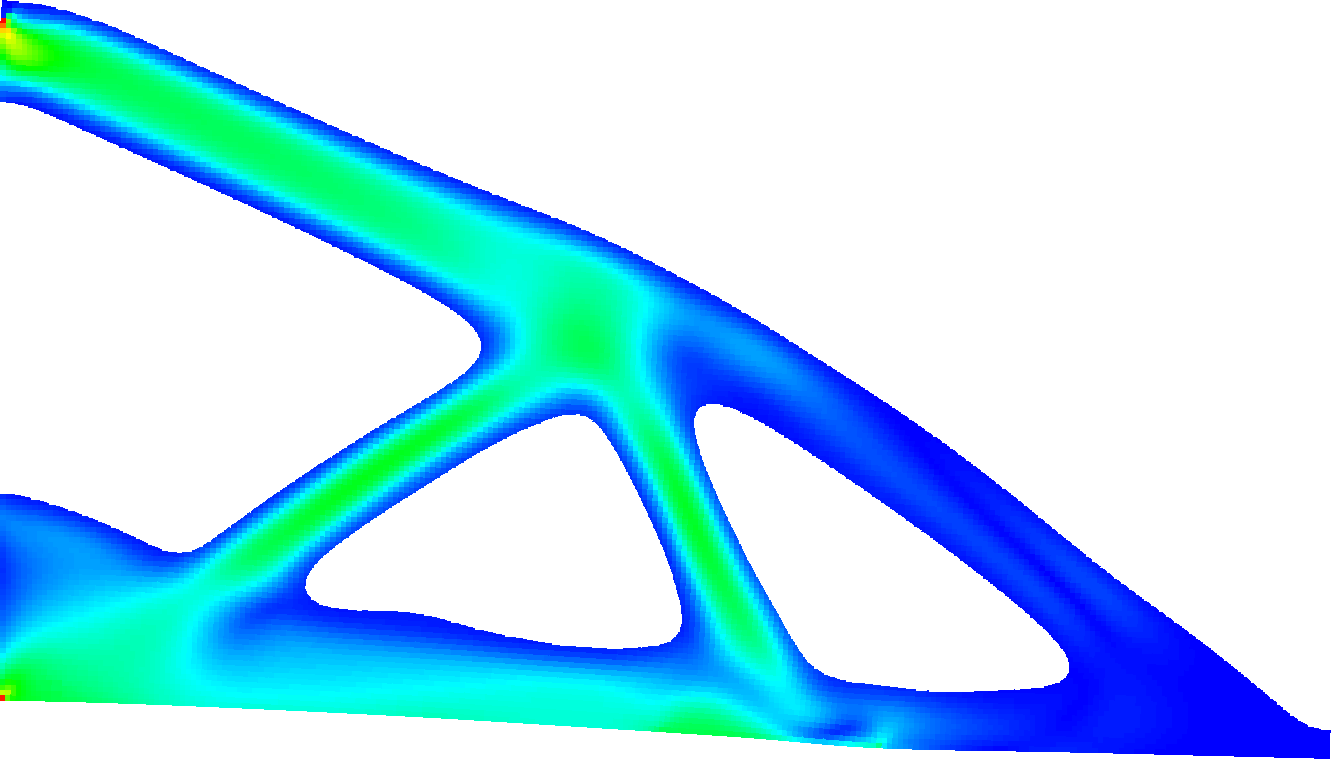}
  \end{subfigure}
  \begin{subfigure}[t]{0.11\textwidth}
    \includegraphics[width=\textwidth]{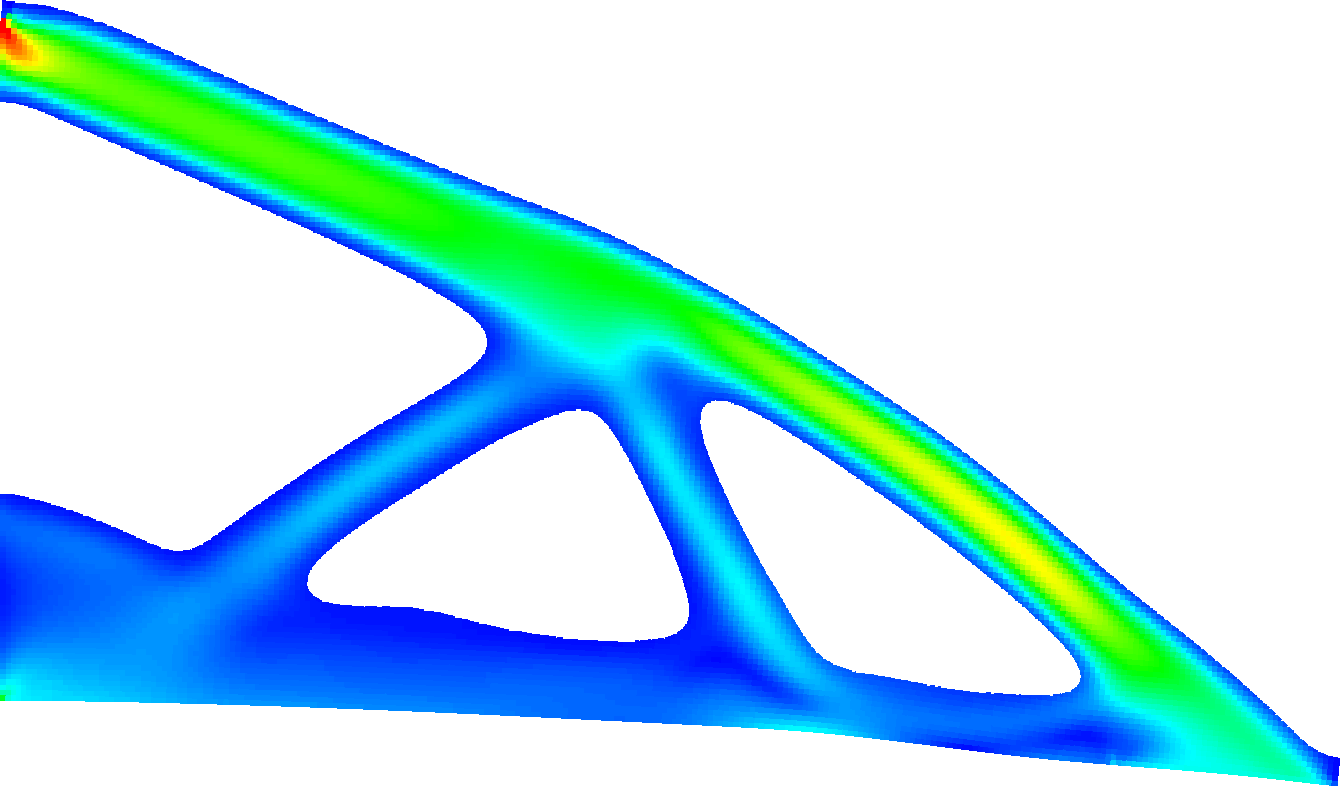}
  \end{subfigure}
  \begin{subfigure}[t]{0.11\textwidth}
    \includegraphics[width=\textwidth]{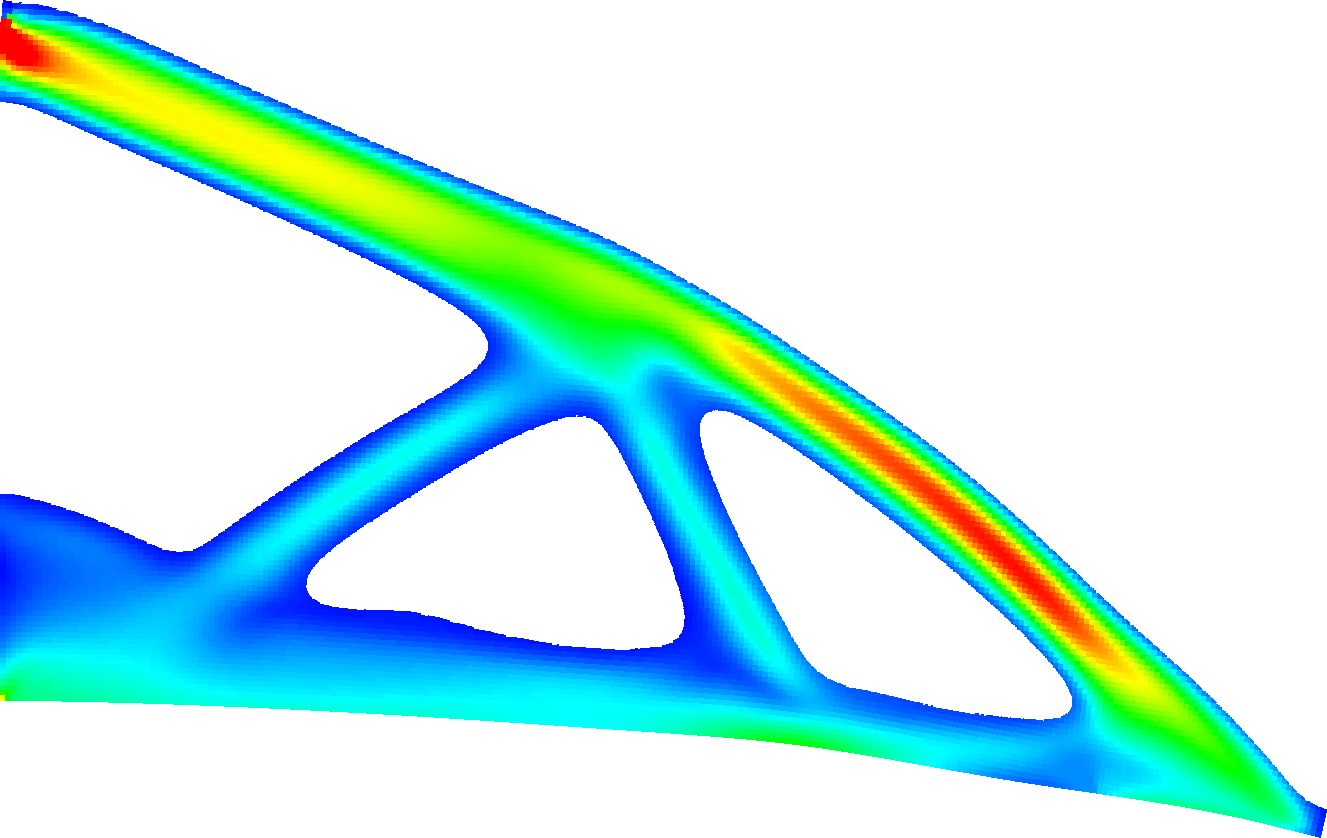}
  \end{subfigure}
  \begin{subfigure}[t]{0.11\textwidth}
    \includegraphics[width=\textwidth]{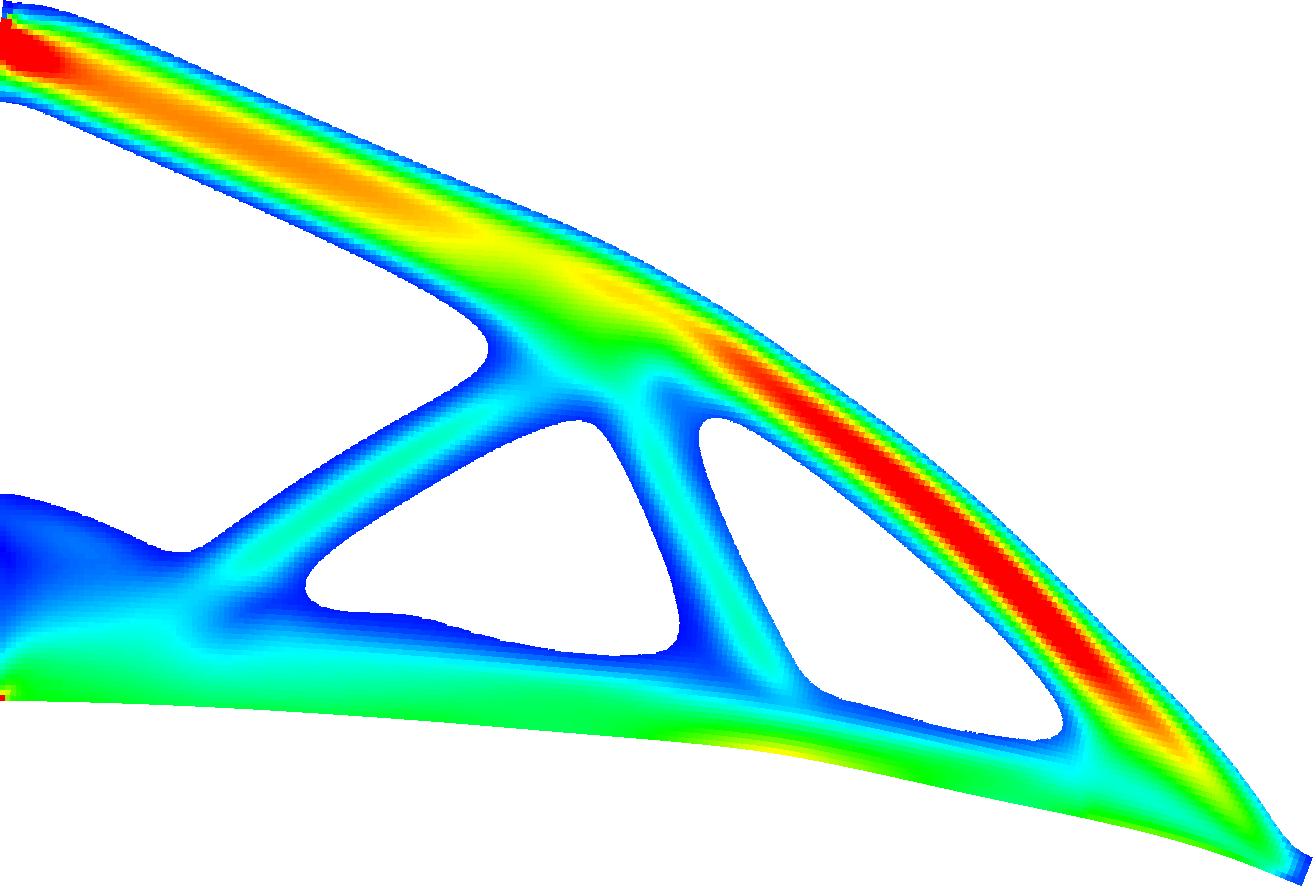}
  \end{subfigure}
  \begin{subfigure}[t]{0.11\textwidth}
    \includegraphics[width=\textwidth]{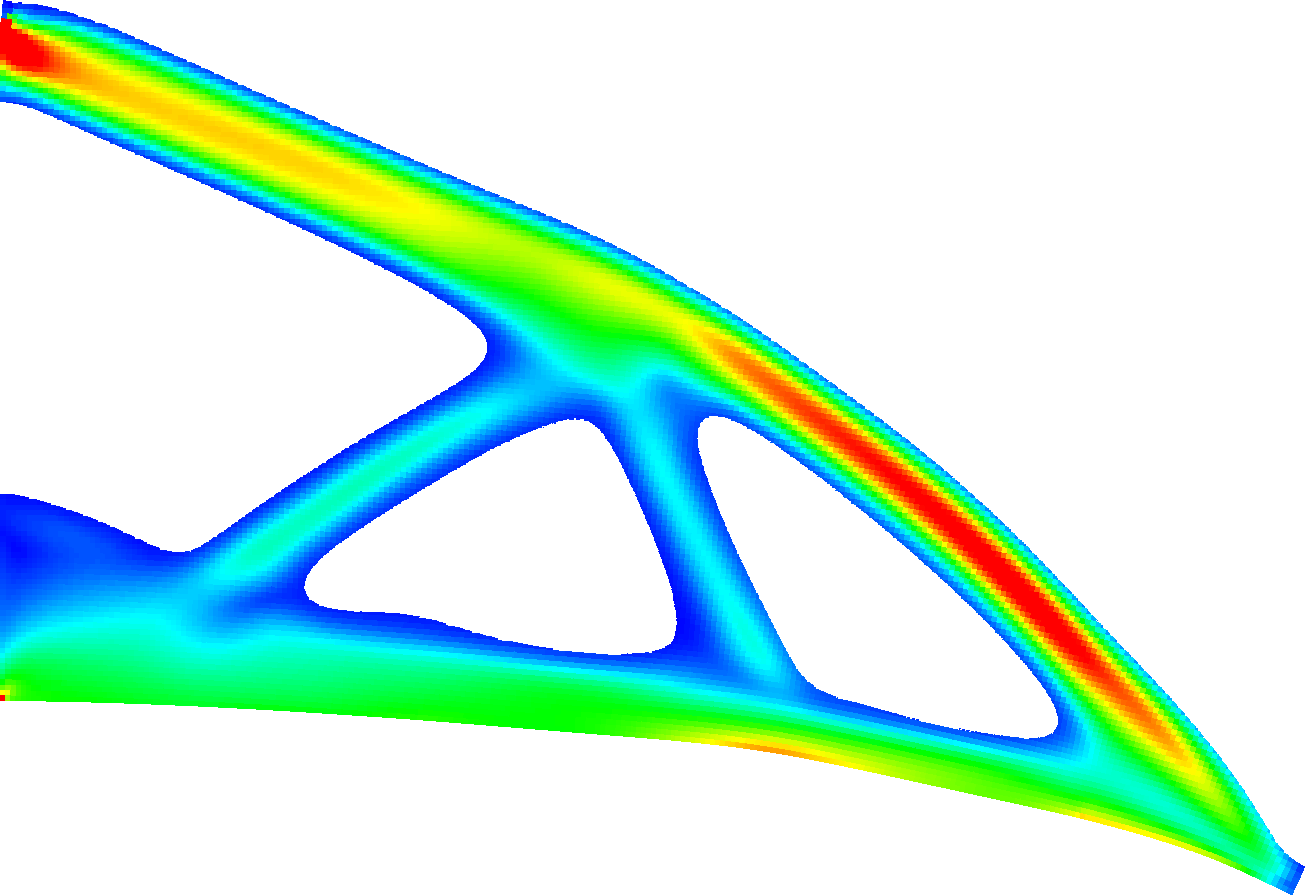}
  \end{subfigure}
  \begin{subfigure}[t]{0.11\textwidth}
    \includegraphics[width=\textwidth]{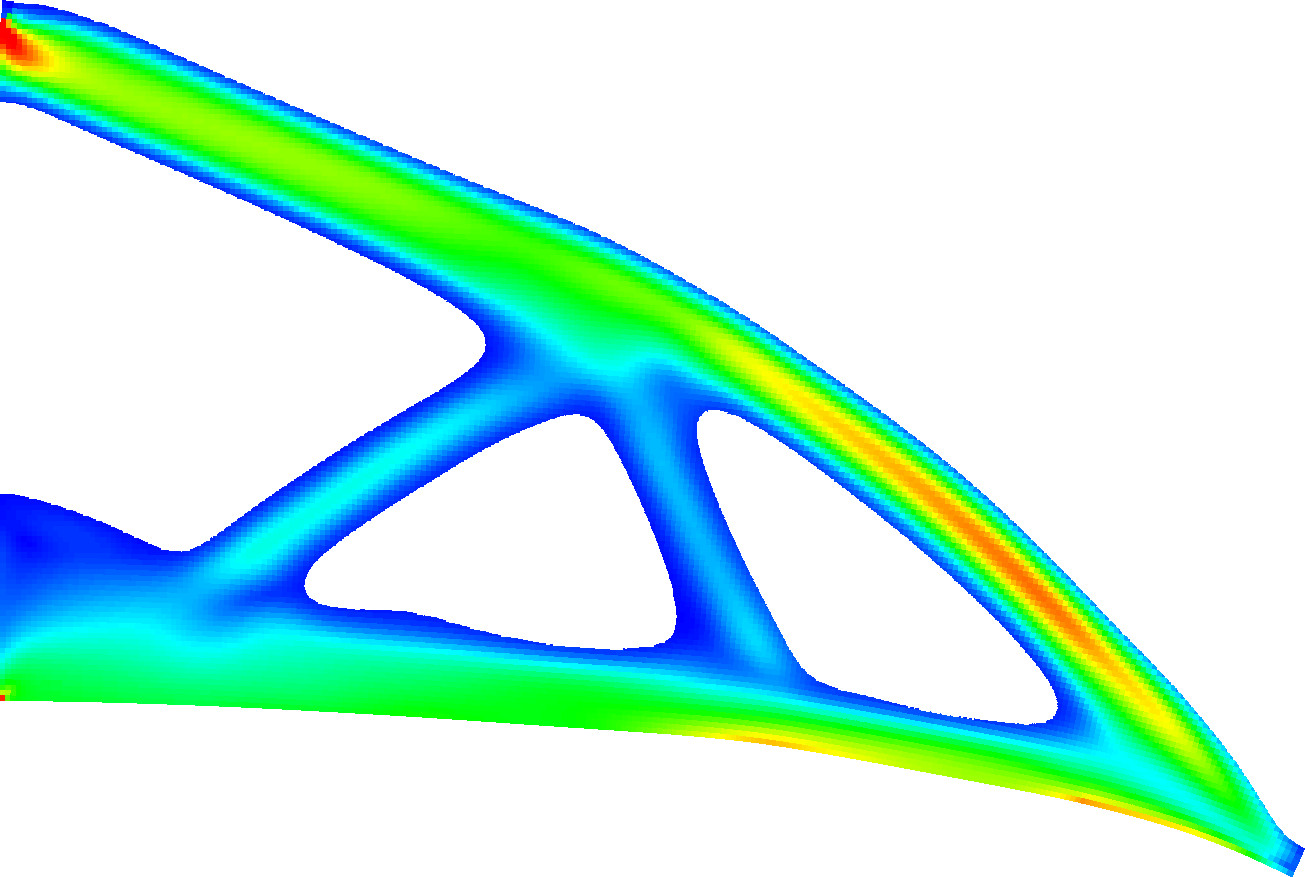}
  \end{subfigure}
  \begin{subfigure}[t]{0.11\textwidth}
    \hspace{\textwidth}
  \end{subfigure}\\[1ex]
   \begin{subfigure}[t]{0.05\textwidth}
    \begin{tikzpicture}
      \node[rotate=90] at (0.0,0.0) {\footnotesize $2^{nd}$ order};
    \end{tikzpicture}
  \end{subfigure}
  \begin{subfigure}[t]{0.11\textwidth}
    \includegraphics[width=\textwidth]{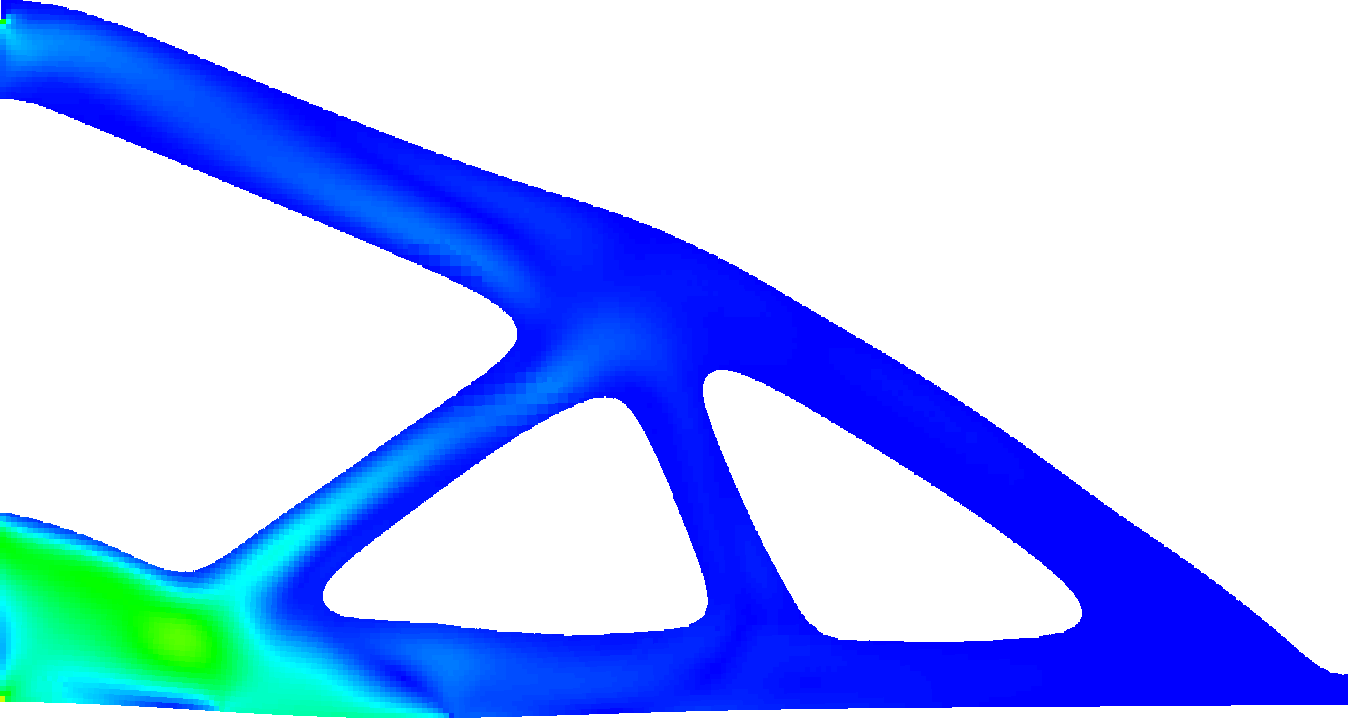}
  \end{subfigure}
  \begin{subfigure}[t]{0.11\textwidth}
    \includegraphics[width=\textwidth]{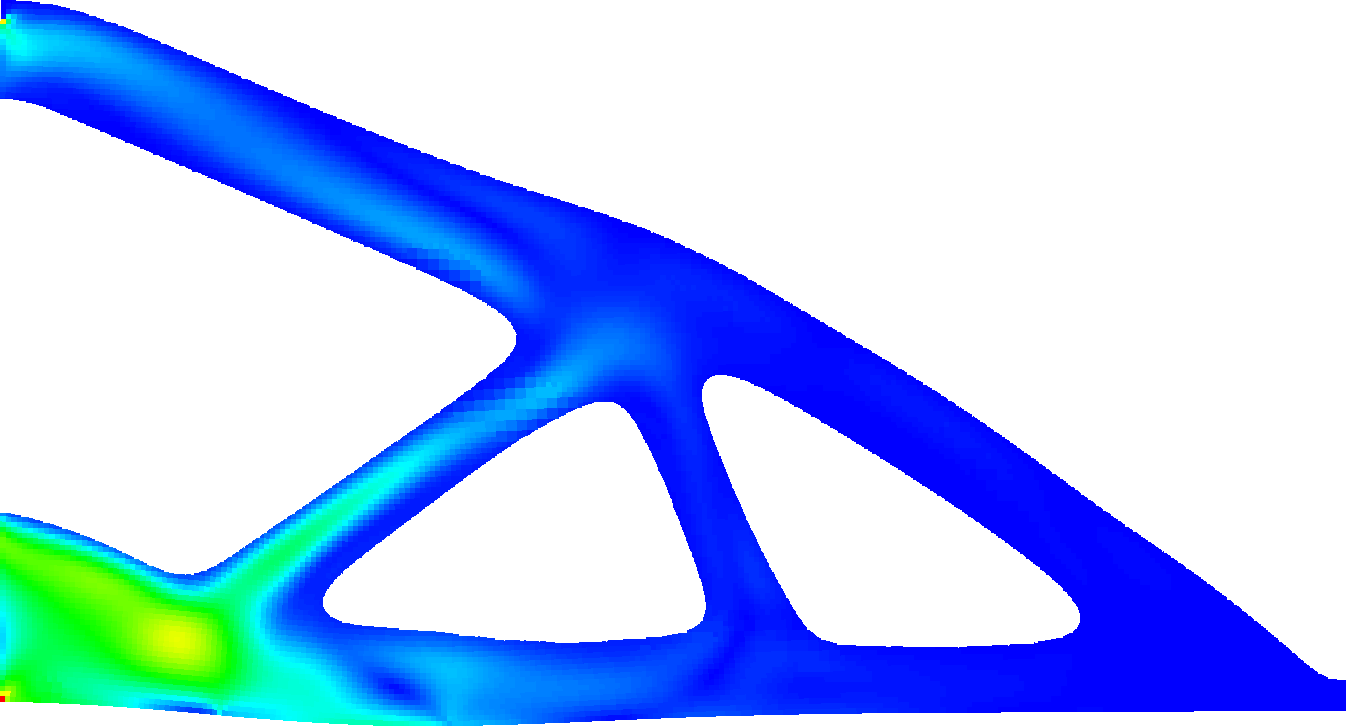}
  \end{subfigure}
  \begin{subfigure}[t]{0.11\textwidth}
    \includegraphics[width=\textwidth]{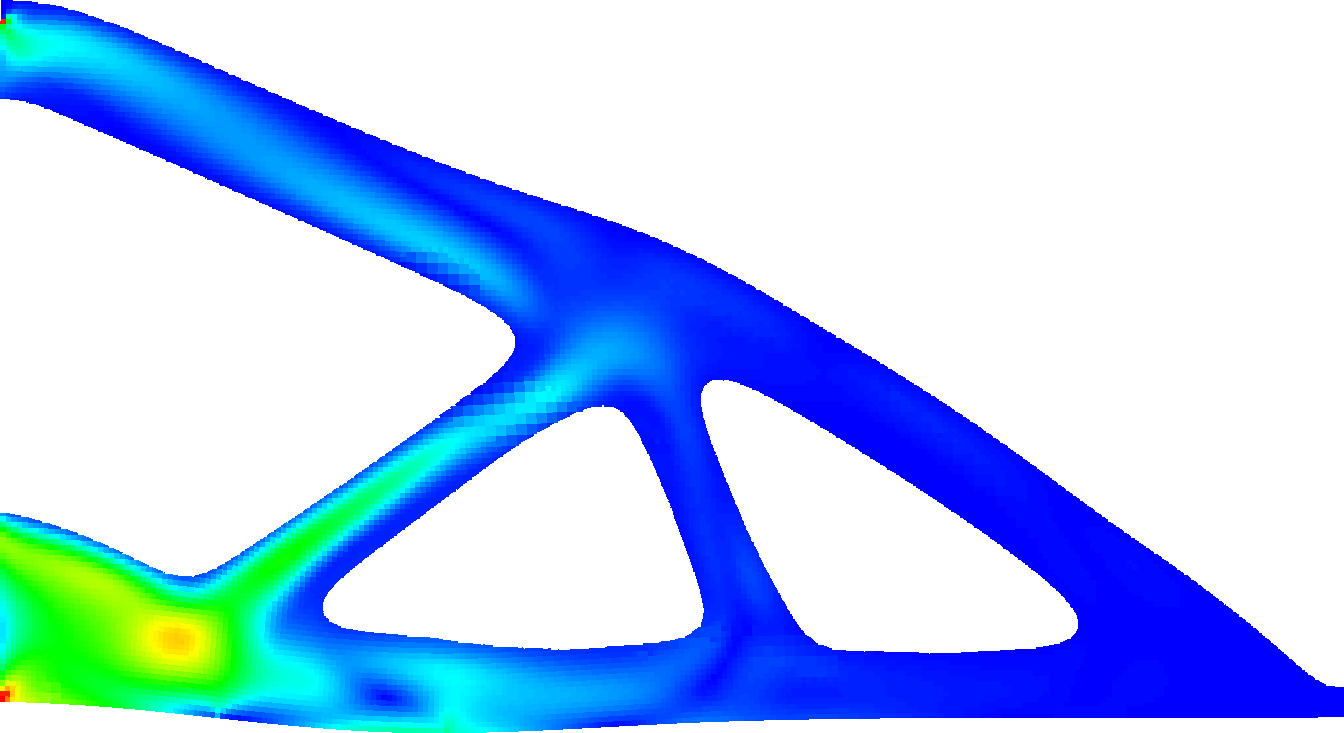}
  \end{subfigure}
  \begin{subfigure}[t]{0.11\textwidth}
    \includegraphics[width=\textwidth]{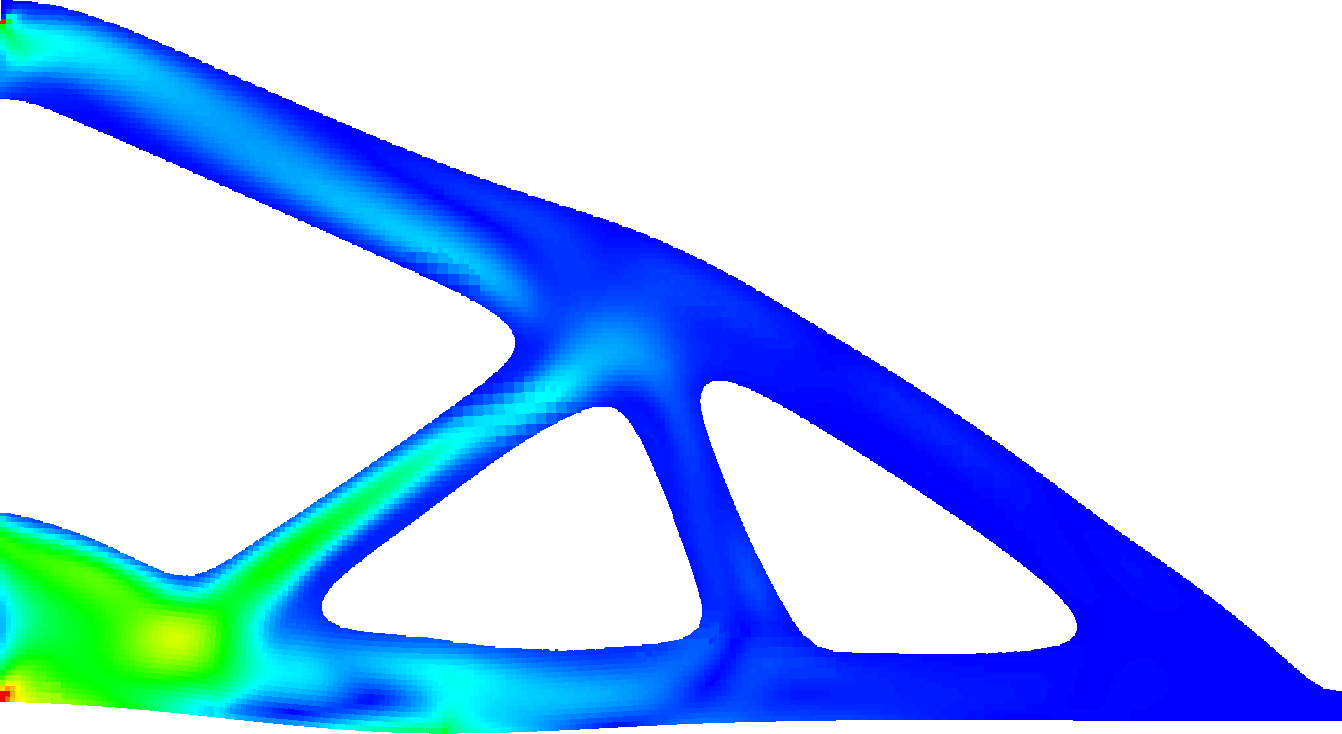}
  \end{subfigure}
  \begin{subfigure}[t]{0.11\textwidth}
    \includegraphics[width=\textwidth]{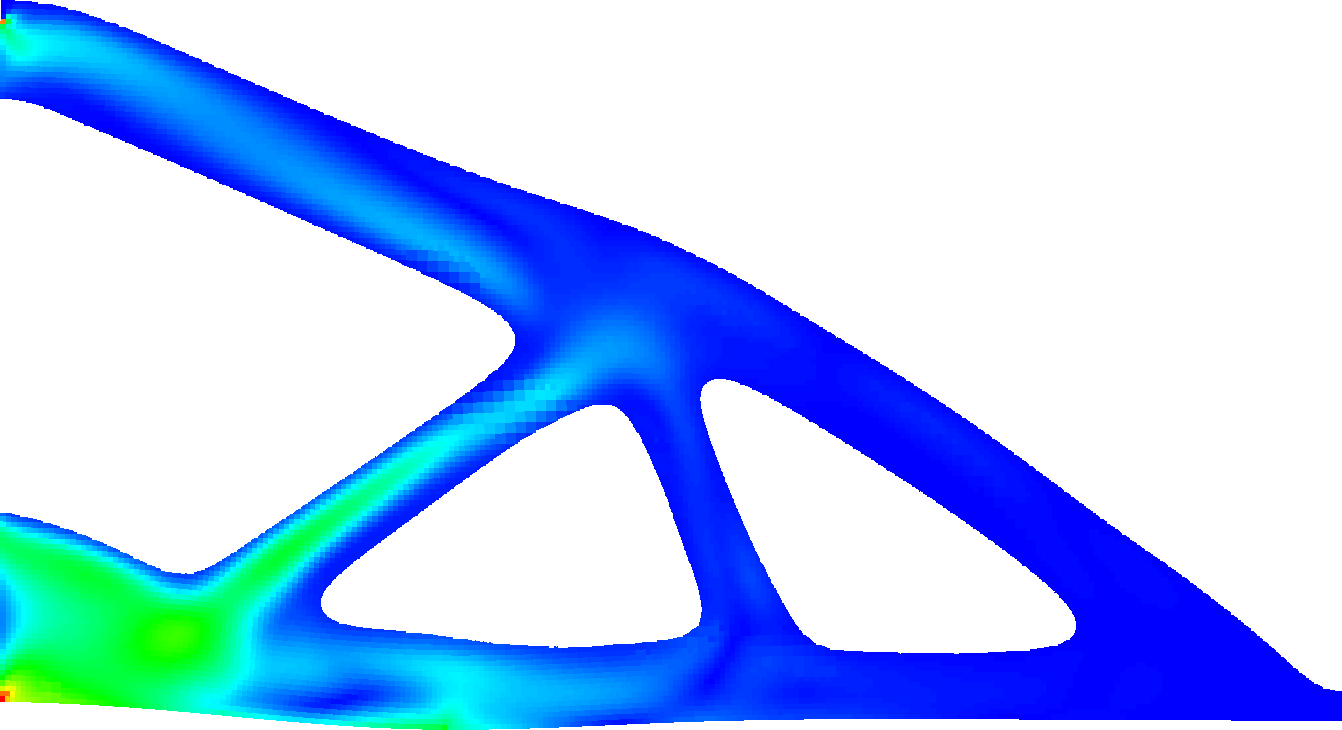}
  \end{subfigure}
  \begin{subfigure}[t]{0.11\textwidth}
    \includegraphics[width=\textwidth]{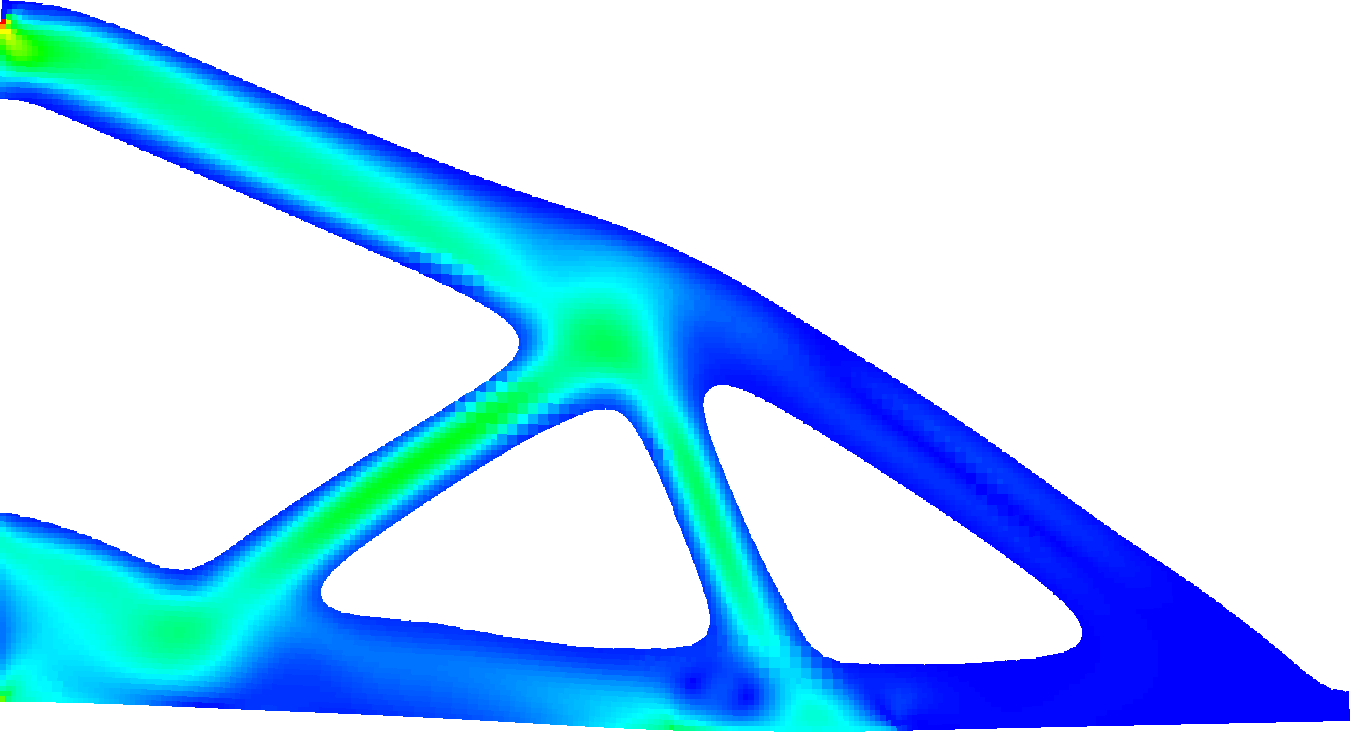}
  \end{subfigure}
  \begin{subfigure}[t]{0.11\textwidth}
    \includegraphics[width=\textwidth]{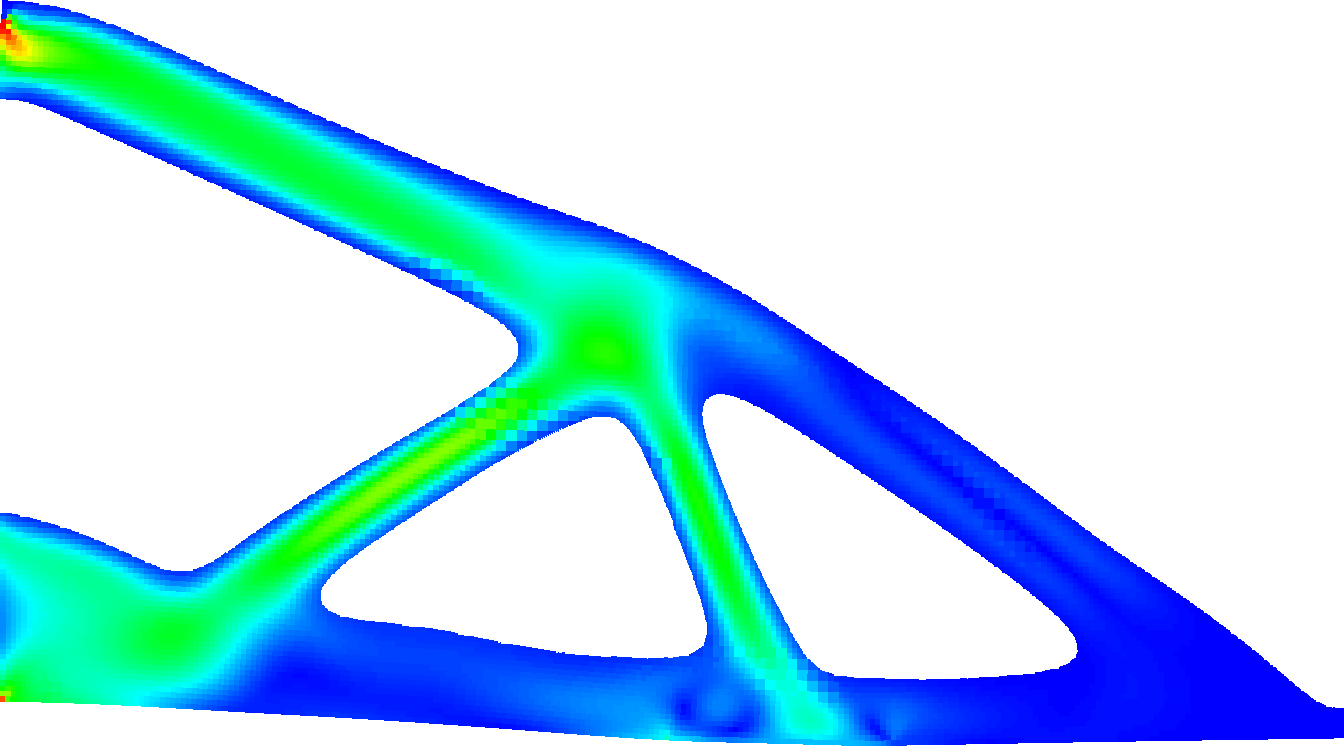}
  \end{subfigure}
  \begin{subfigure}[t]{0.11\textwidth}
    \includegraphics[width=\textwidth]{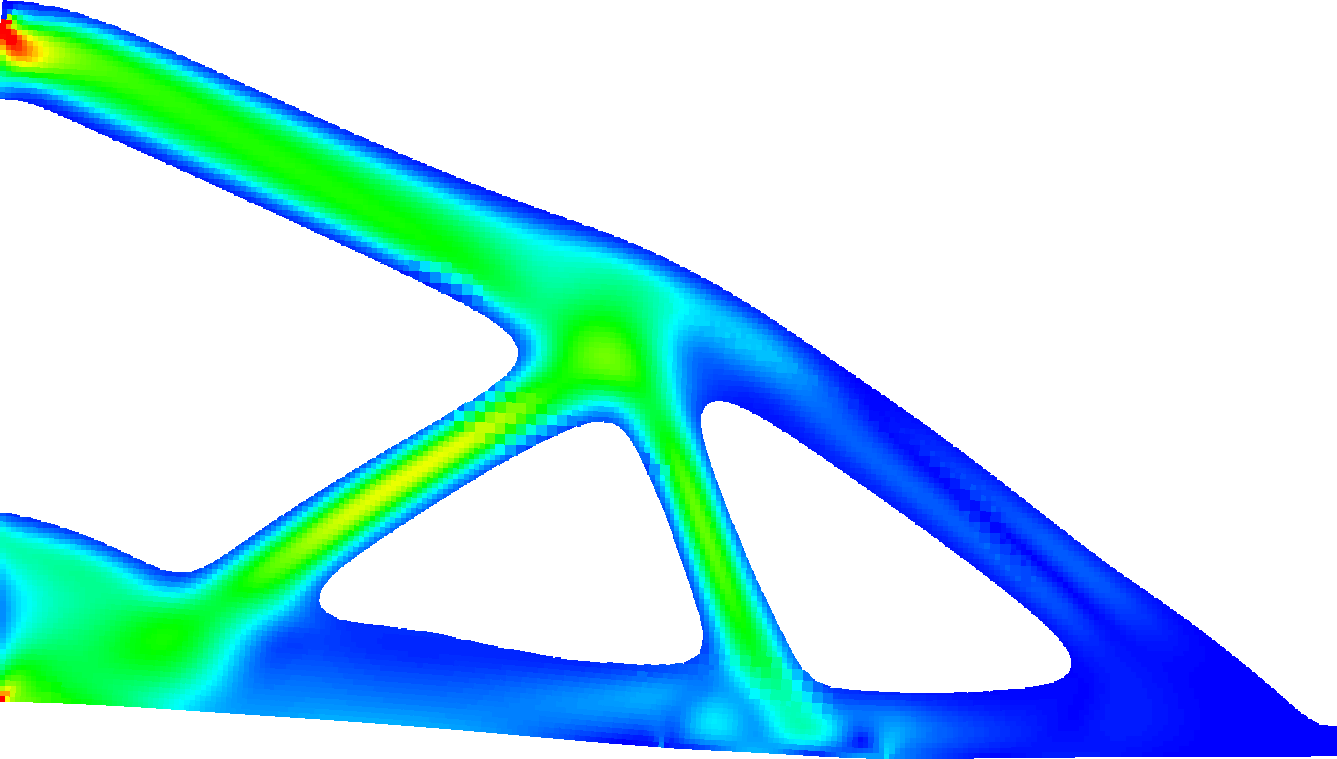}
  \end{subfigure}\\[1ex]
  \begin{subfigure}[t]{0.06\textwidth}
    \begin{tikzpicture}
      \node[rotate=90,align=left] at (0.0,0.0) {\footnotesize eq. loads\\
      \footnotesize and prob.};
    \end{tikzpicture}
  \end{subfigure}
  \begin{subfigure}[t]{0.11\textwidth}
    \includegraphics[width=\textwidth]{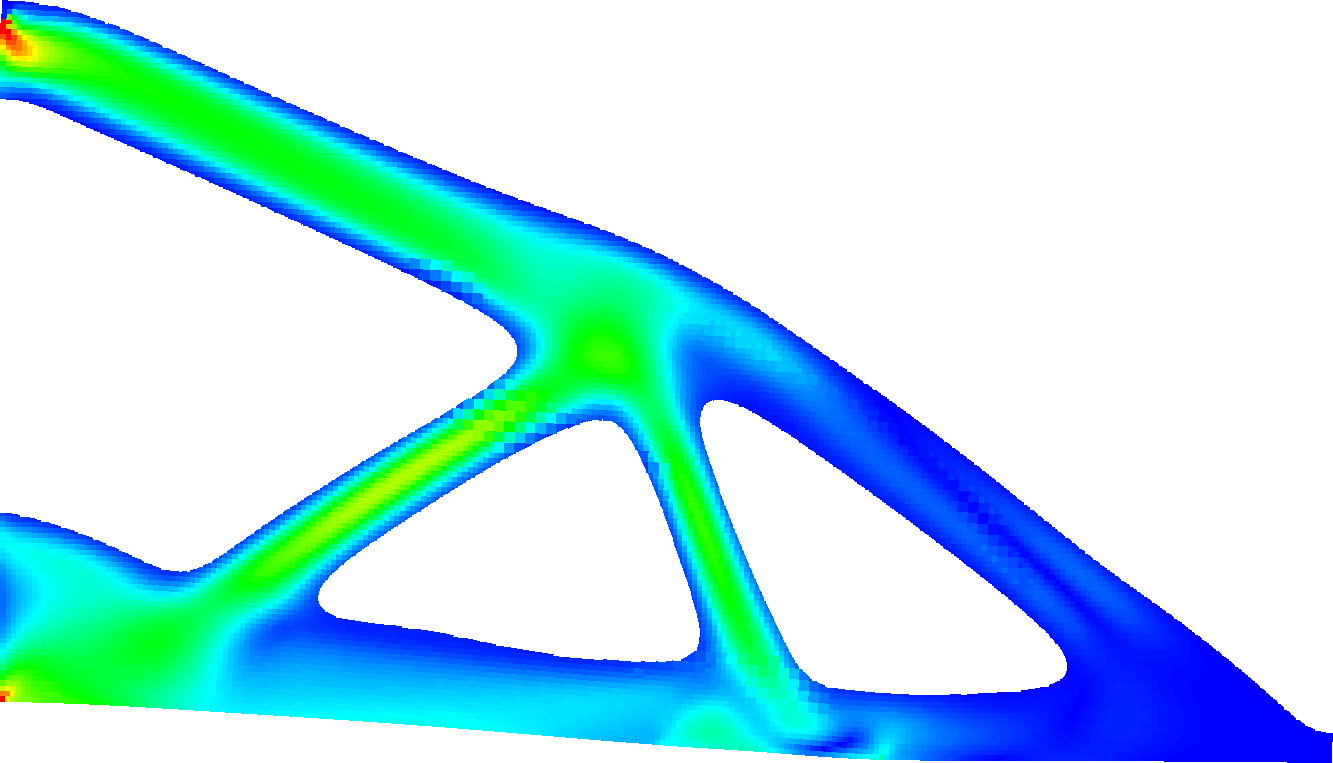}
  \end{subfigure}
  \begin{subfigure}[t]{0.11\textwidth}
    \includegraphics[width=\textwidth]{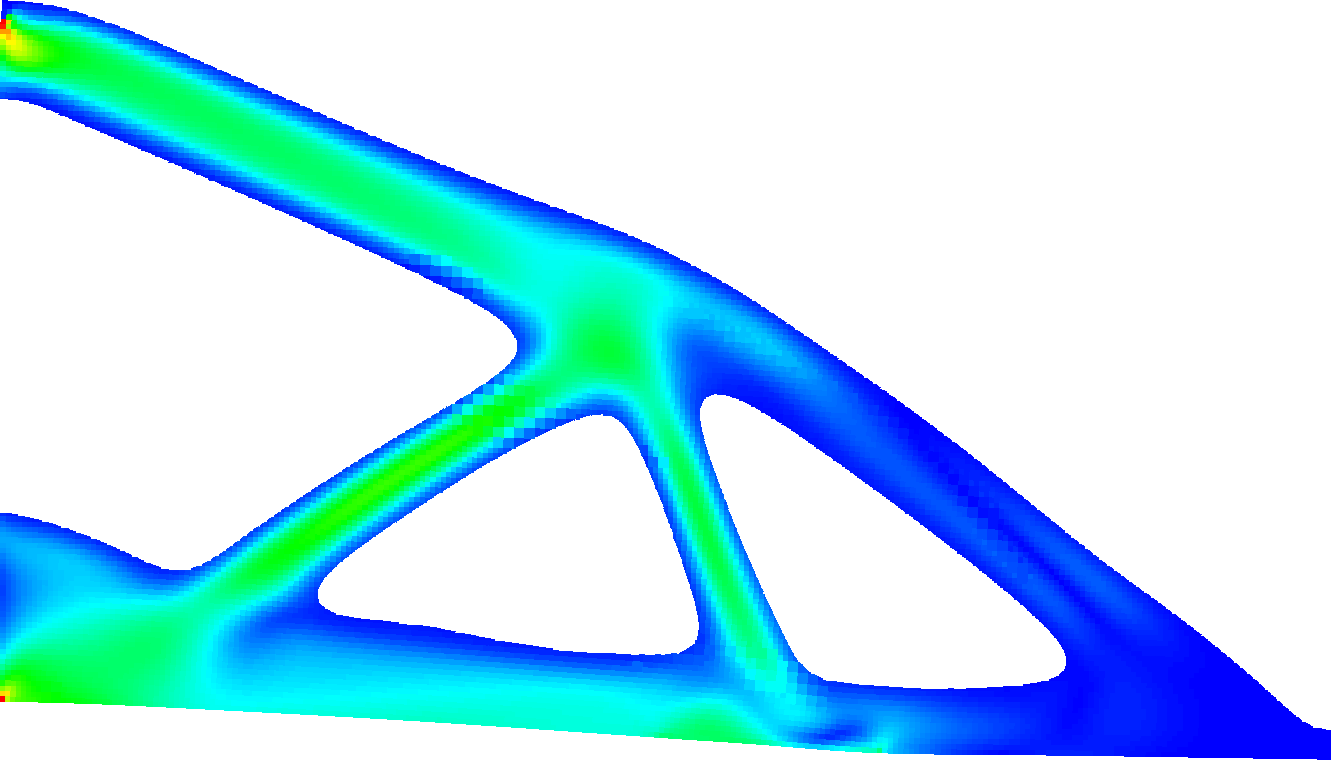}
  \end{subfigure}
  \begin{subfigure}[t]{0.11\textwidth}
    \includegraphics[width=\textwidth]{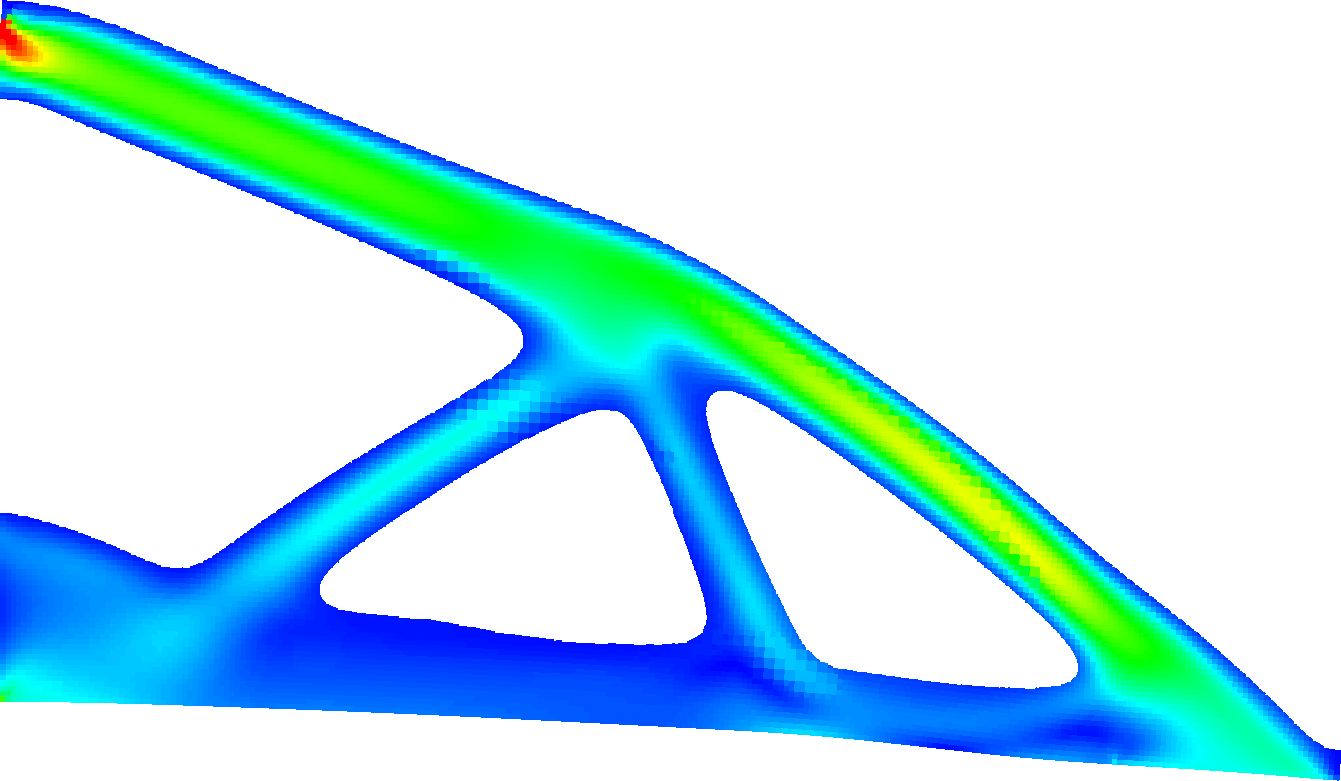}
  \end{subfigure}
  \begin{subfigure}[t]{0.11\textwidth}
    \includegraphics[width=\textwidth]{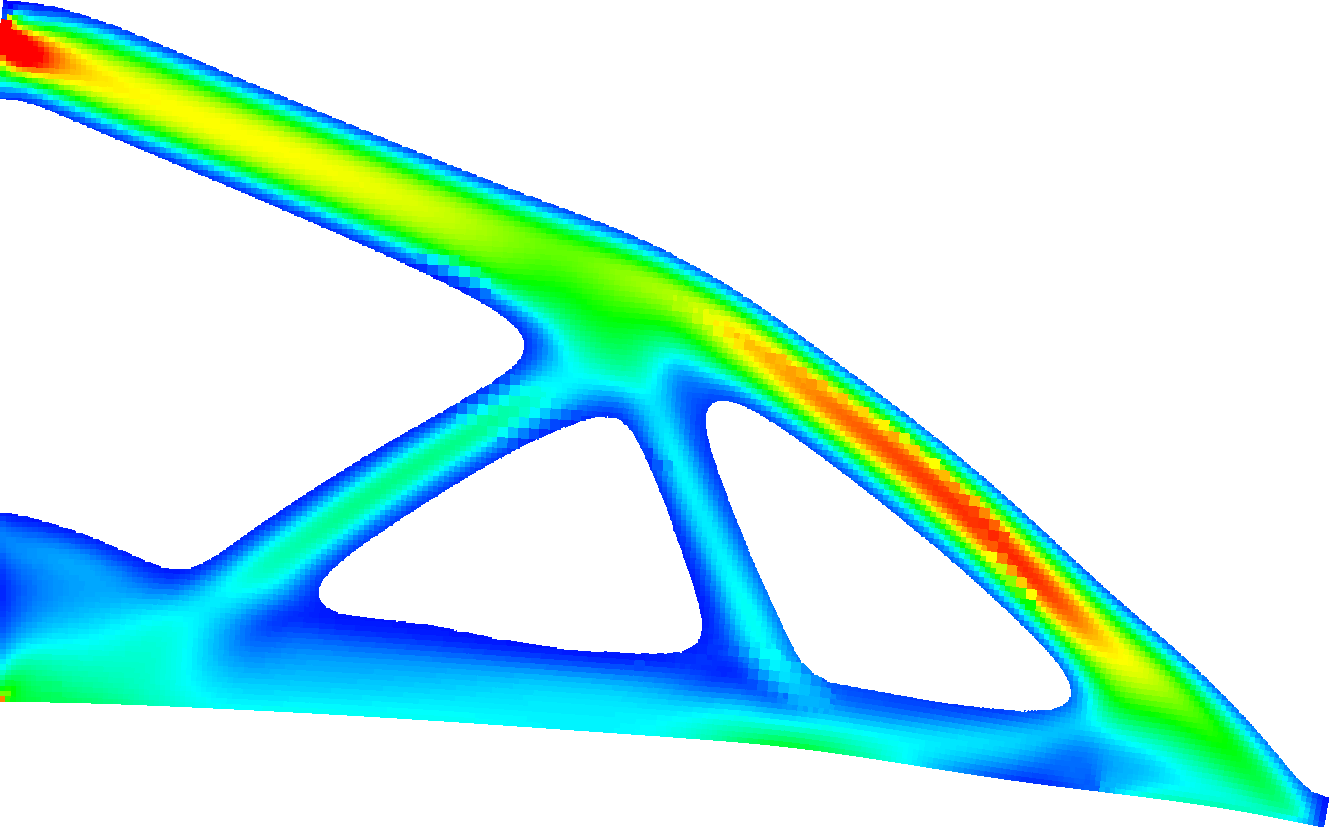}
  \end{subfigure}
  \begin{subfigure}[t]{0.11\textwidth}
    \includegraphics[width=\textwidth]{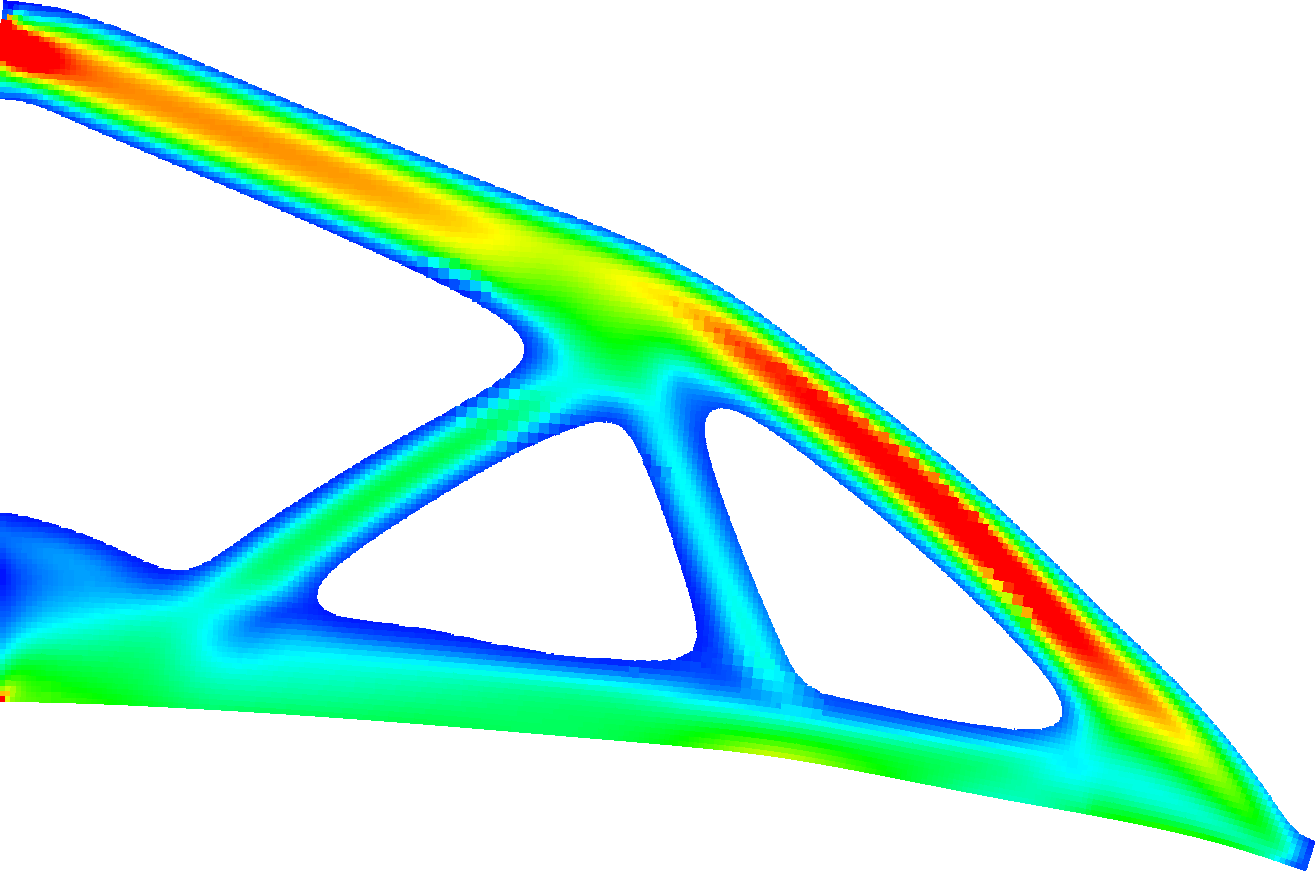}
  \end{subfigure}
  \begin{subfigure}[t]{0.11\textwidth}
    \includegraphics[width=\textwidth]{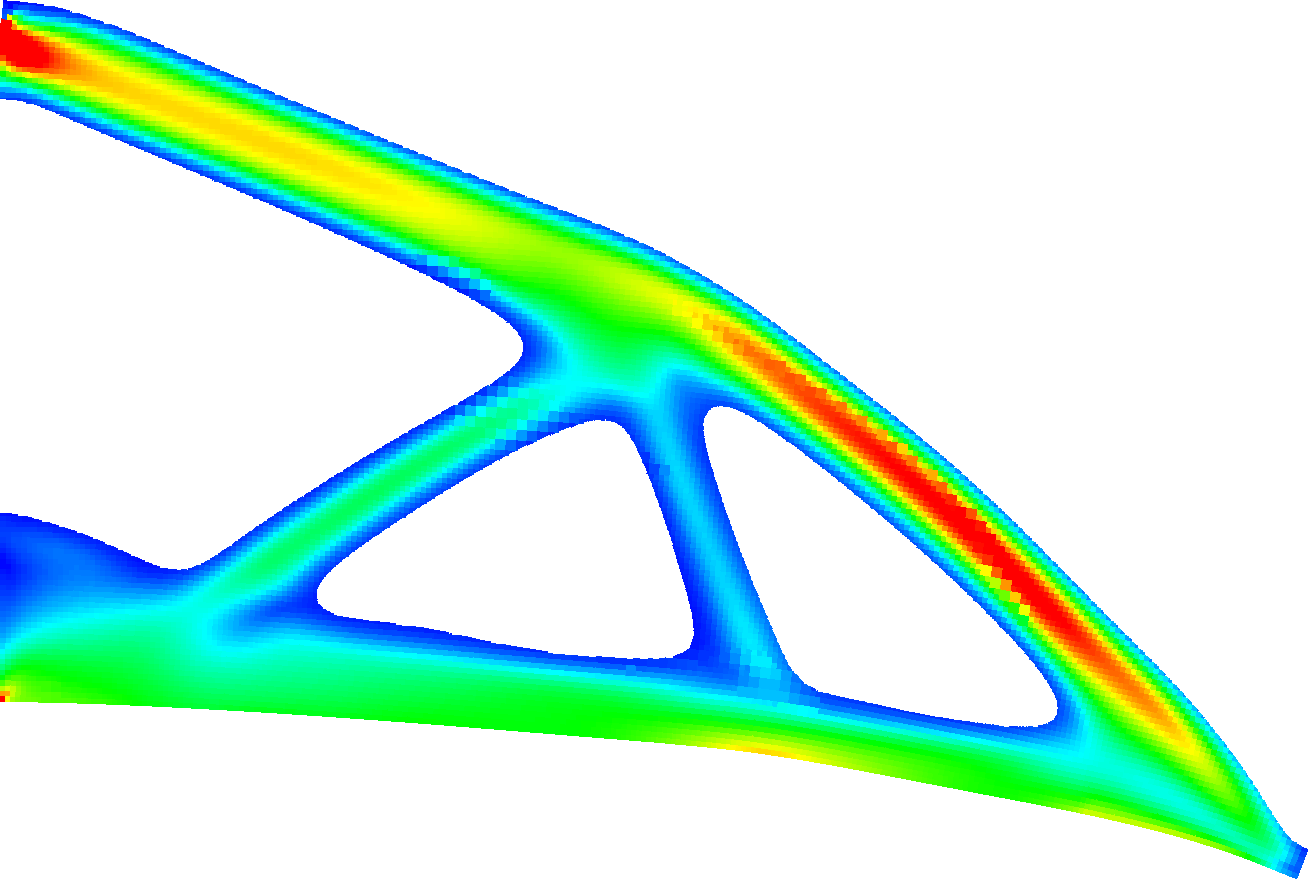}
  \end{subfigure}
  \begin{subfigure}[t]{0.11\textwidth}
    \includegraphics[width=\textwidth]{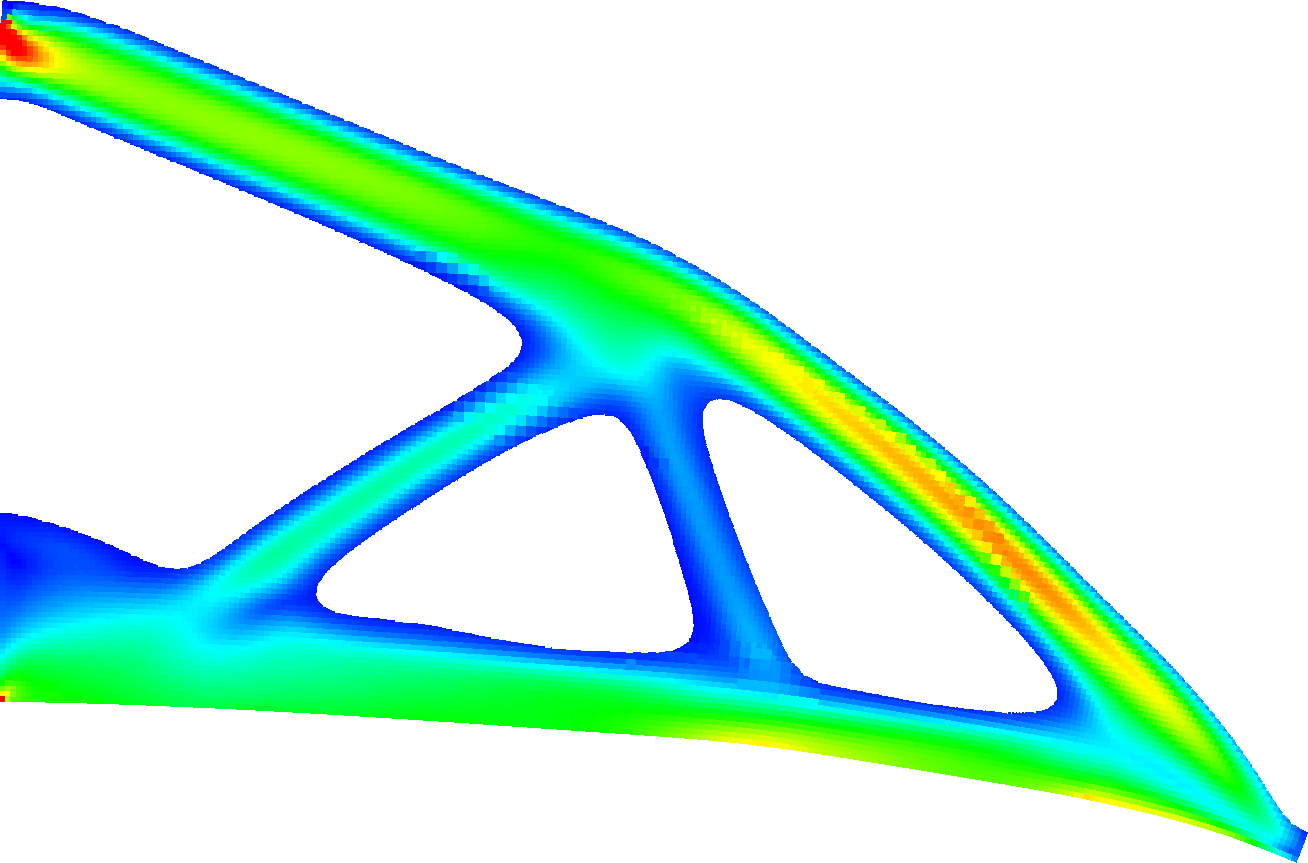}
  \end{subfigure}
  \begin{subfigure}[t]{0.11\textwidth}
    \hspace{\textwidth}
  \end{subfigure}\\[2ex]
  \hrule \vspace{2ex}
                    
  \begin{subfigure}[t]{0.05\textwidth}
    \begin{tikzpicture}
      \node[rotate=90] at (0.0,0.0) {\footnotesize $1^{st}$ order};
    \end{tikzpicture}
  \end{subfigure}
  \begin{subfigure}[t]{0.11\textwidth}
    \includegraphics[width=\textwidth]{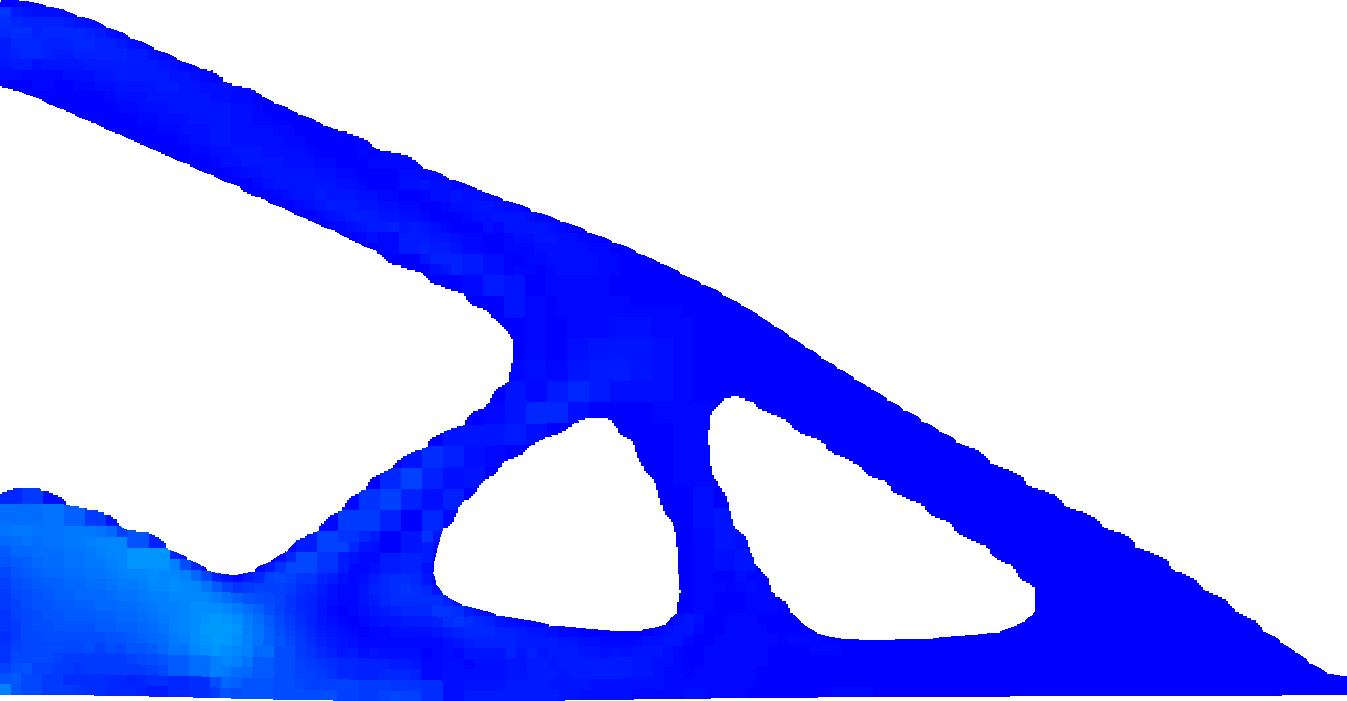}
  \end{subfigure}
  \begin{subfigure}[t]{0.11\textwidth}
    \includegraphics[width=\textwidth]{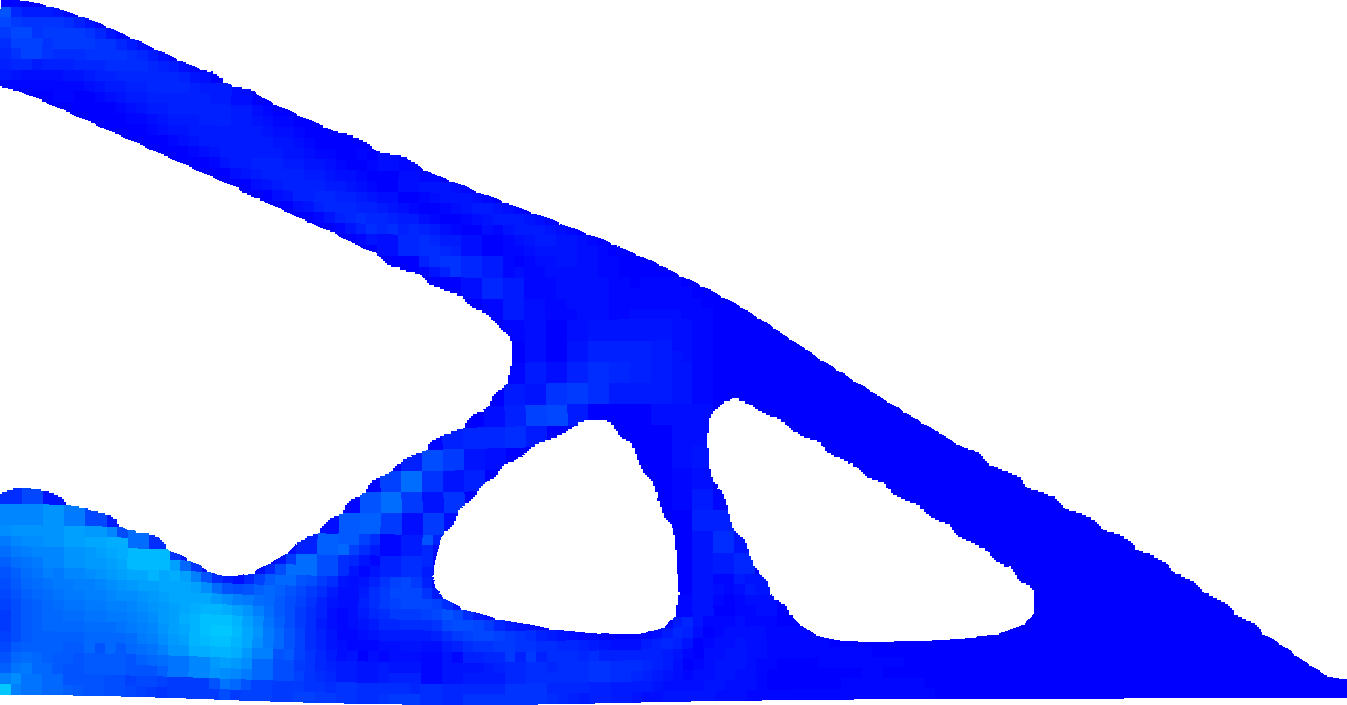}
  \end{subfigure}
  \begin{subfigure}[t]{0.11\textwidth}
    \includegraphics[width=\textwidth]{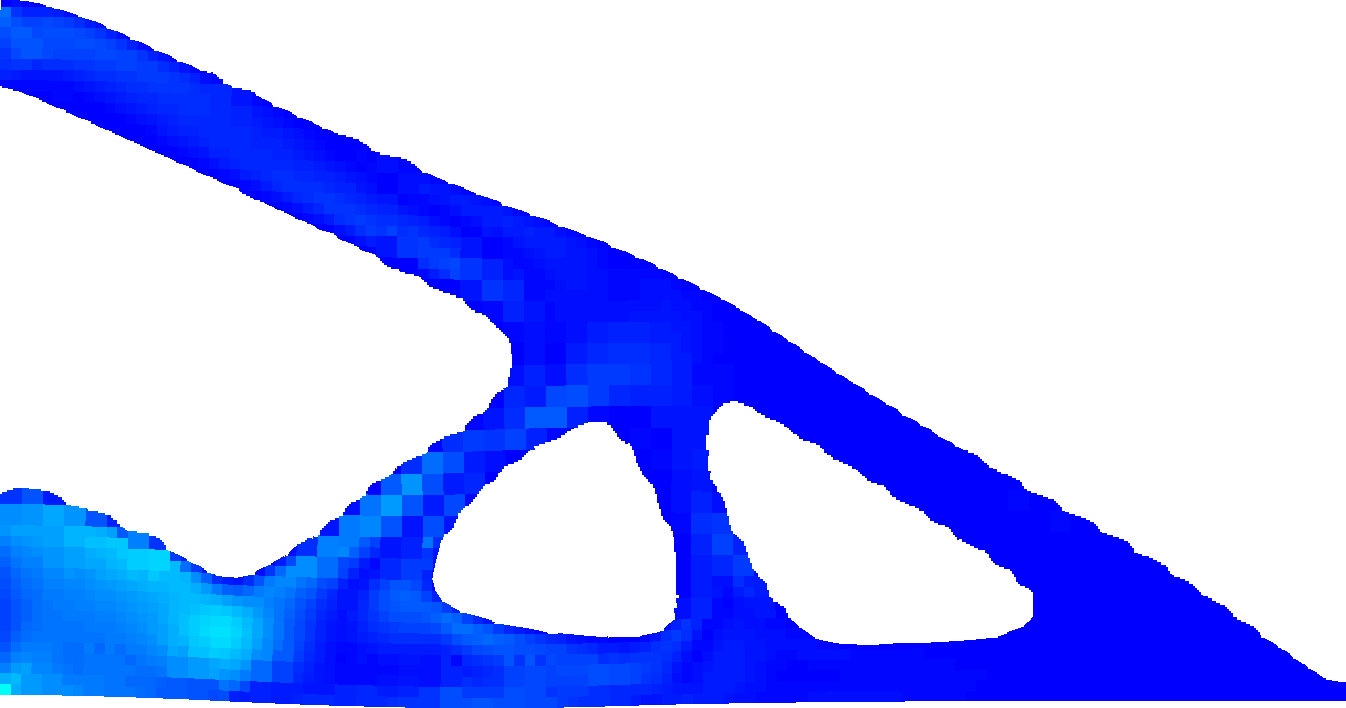}
  \end{subfigure}
  \begin{subfigure}[t]{0.11\textwidth}
    \includegraphics[width=\textwidth]{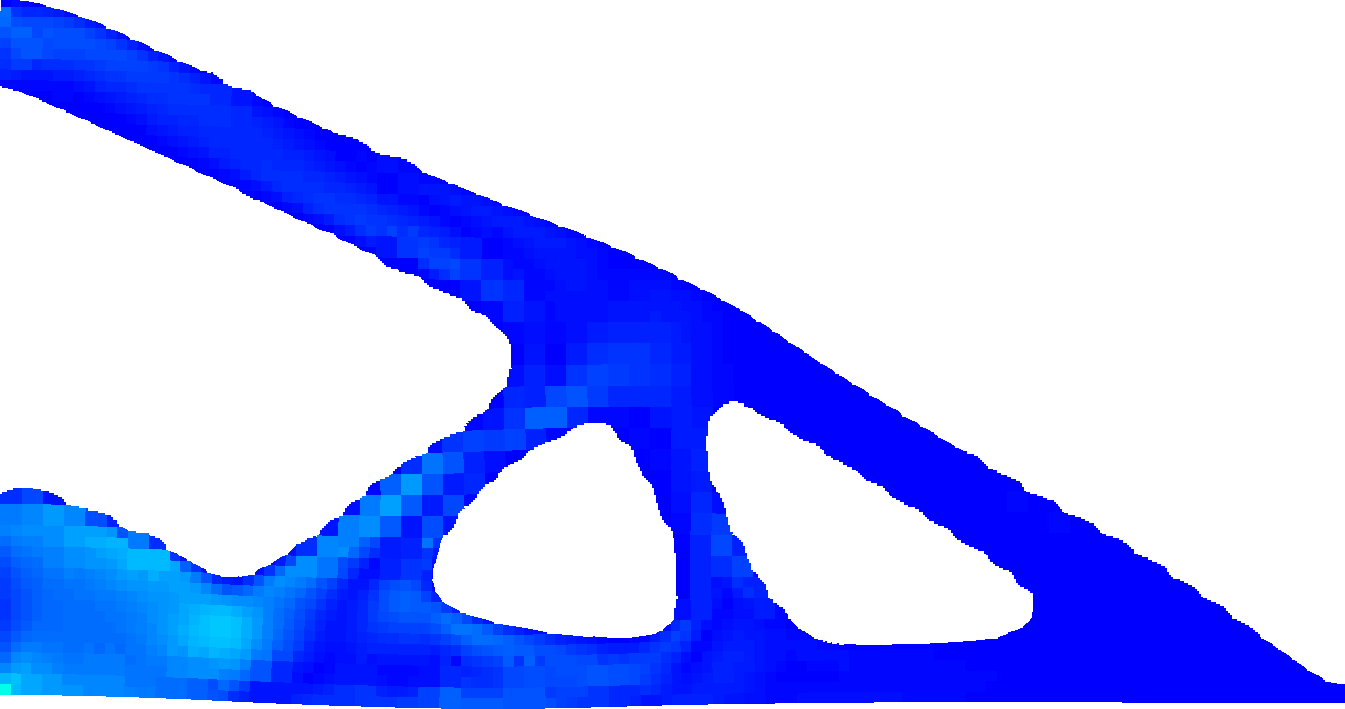}
  \end{subfigure}
  \begin{subfigure}[t]{0.11\textwidth}
    \includegraphics[width=\textwidth]{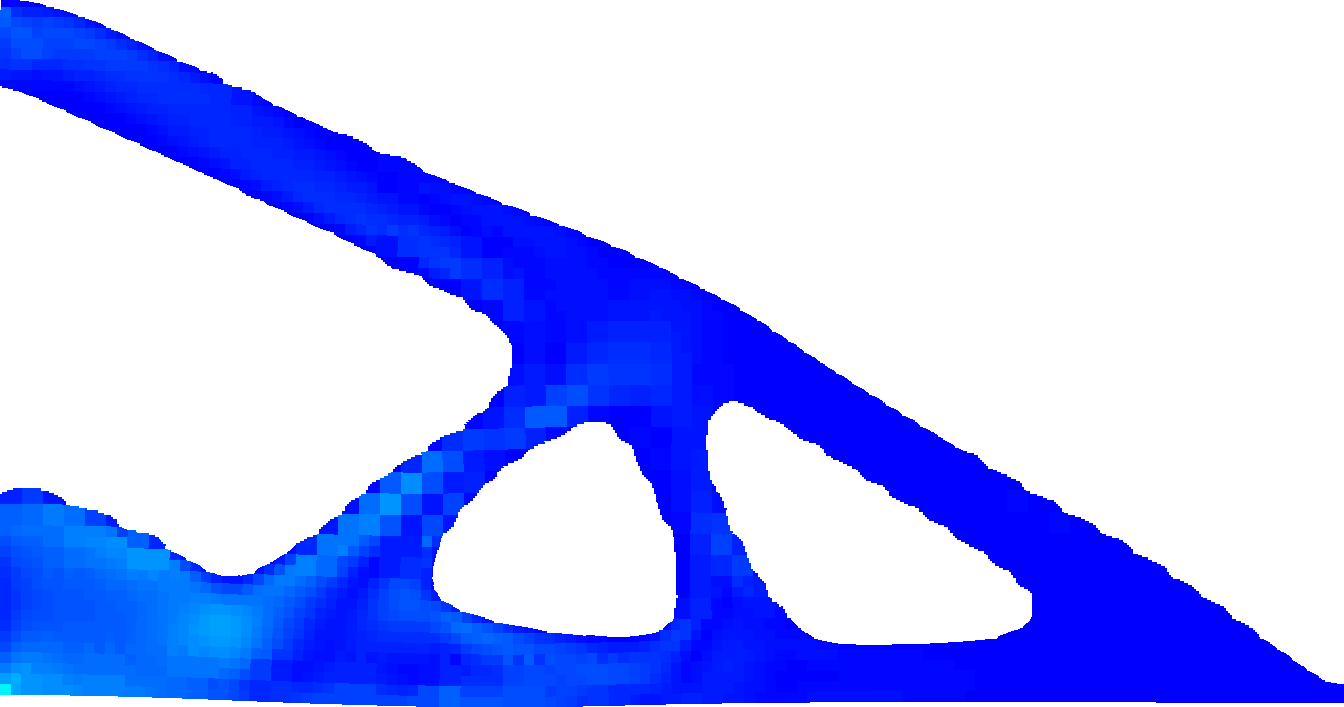}
  \end{subfigure}
  \begin{subfigure}[t]{0.11\textwidth}
    \includegraphics[width=\textwidth]{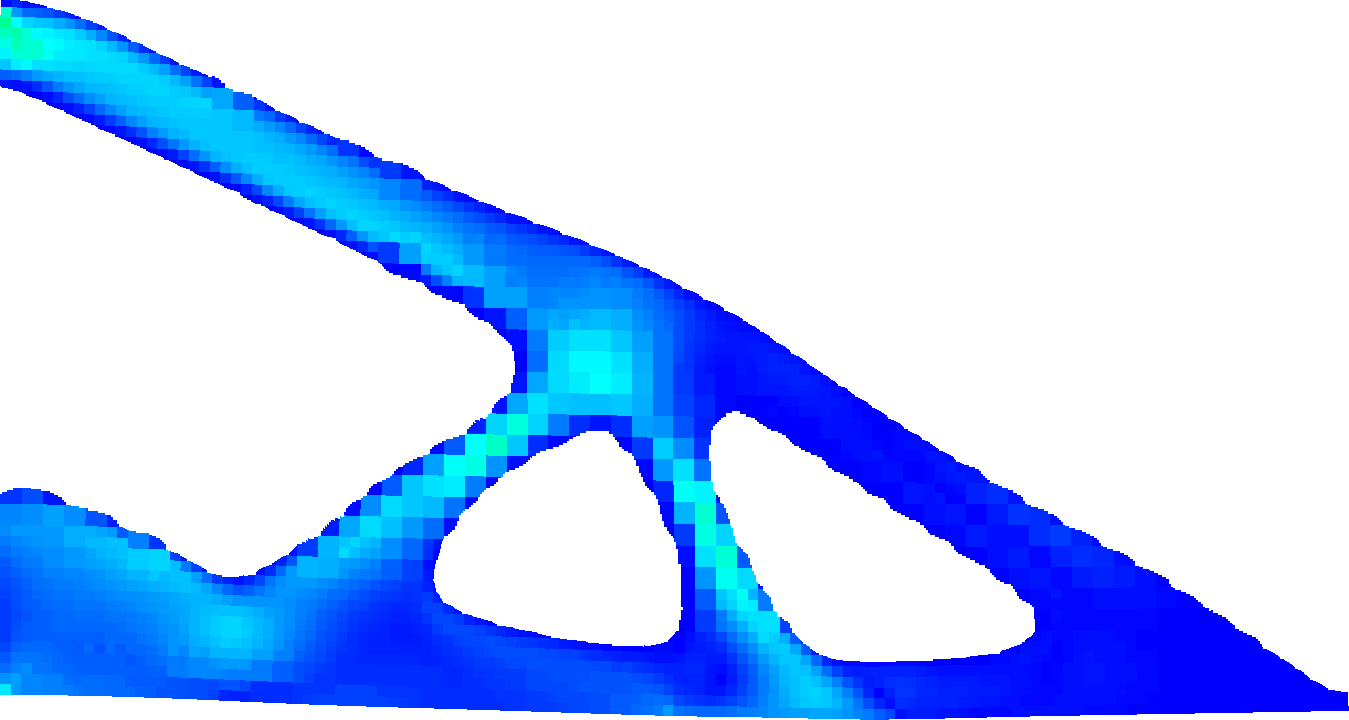}
  \end{subfigure}
  \begin{subfigure}[t]{0.11\textwidth}
    \includegraphics[width=\textwidth]{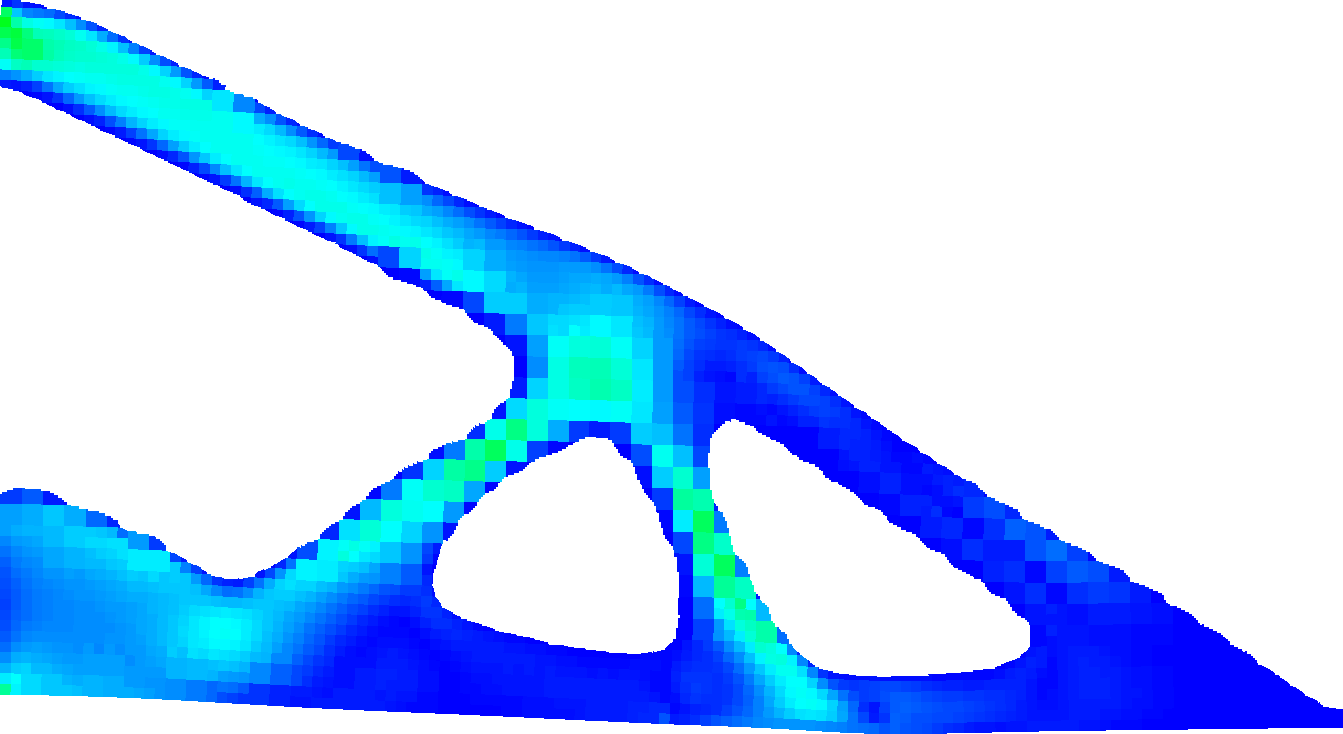}
  \end{subfigure}
  \begin{subfigure}[t]{0.11\textwidth}
    \includegraphics[width=\textwidth]{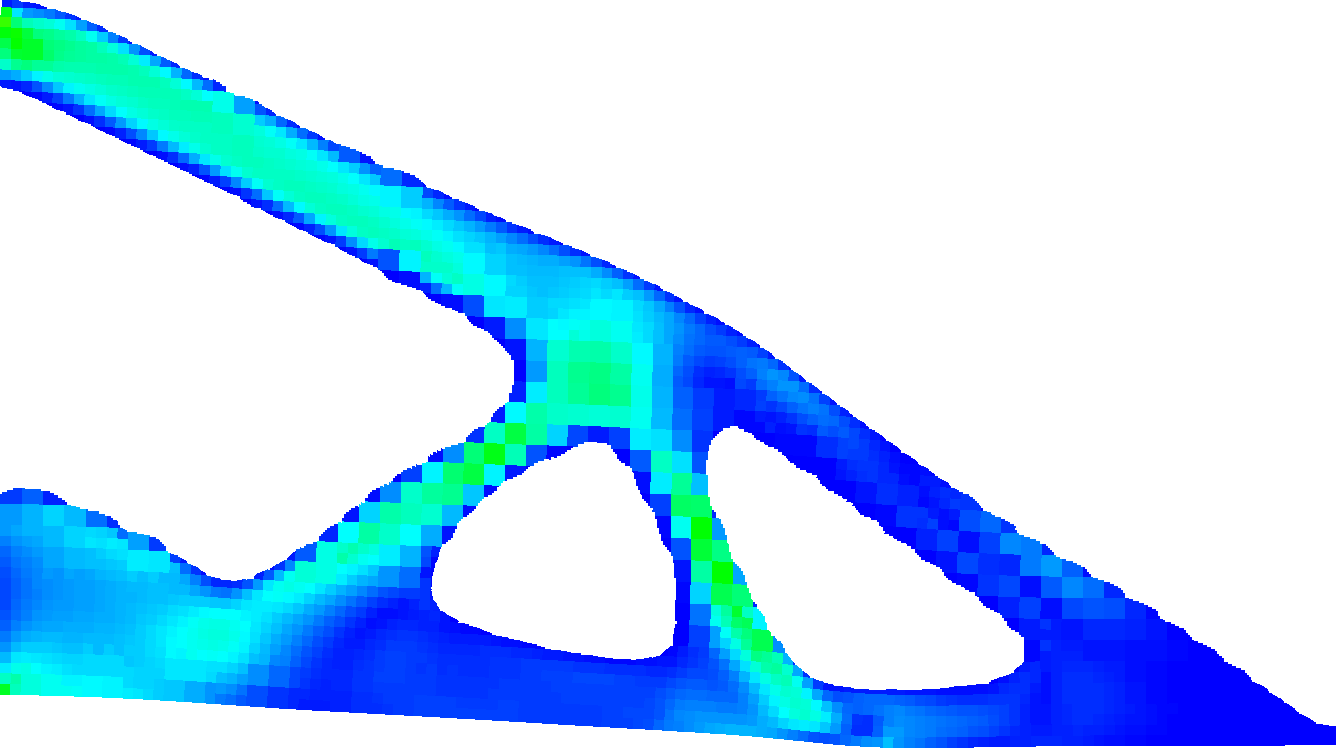}
  \end{subfigure}\\[1ex]
  \begin{subfigure}[t]{0.06\textwidth}
    \begin{tikzpicture}
      \node[rotate=90,align=left] at (0.0,0.0) {\footnotesize var. loads\\
      \footnotesize and prob.};
    \end{tikzpicture}
  \end{subfigure}
  \begin{subfigure}[t]{0.11\textwidth}
    \includegraphics[width=\textwidth]{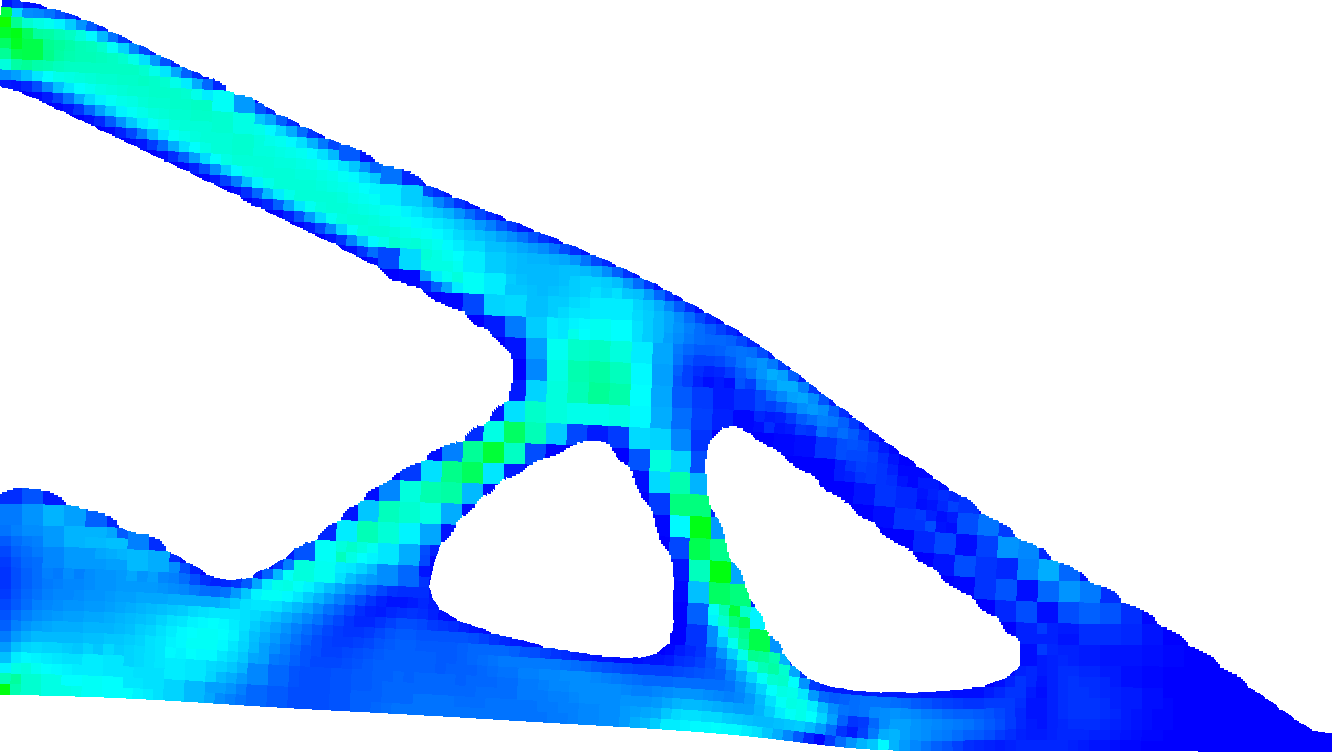}
  \end{subfigure}
  \begin{subfigure}[t]{0.11\textwidth}
    \includegraphics[width=\textwidth]{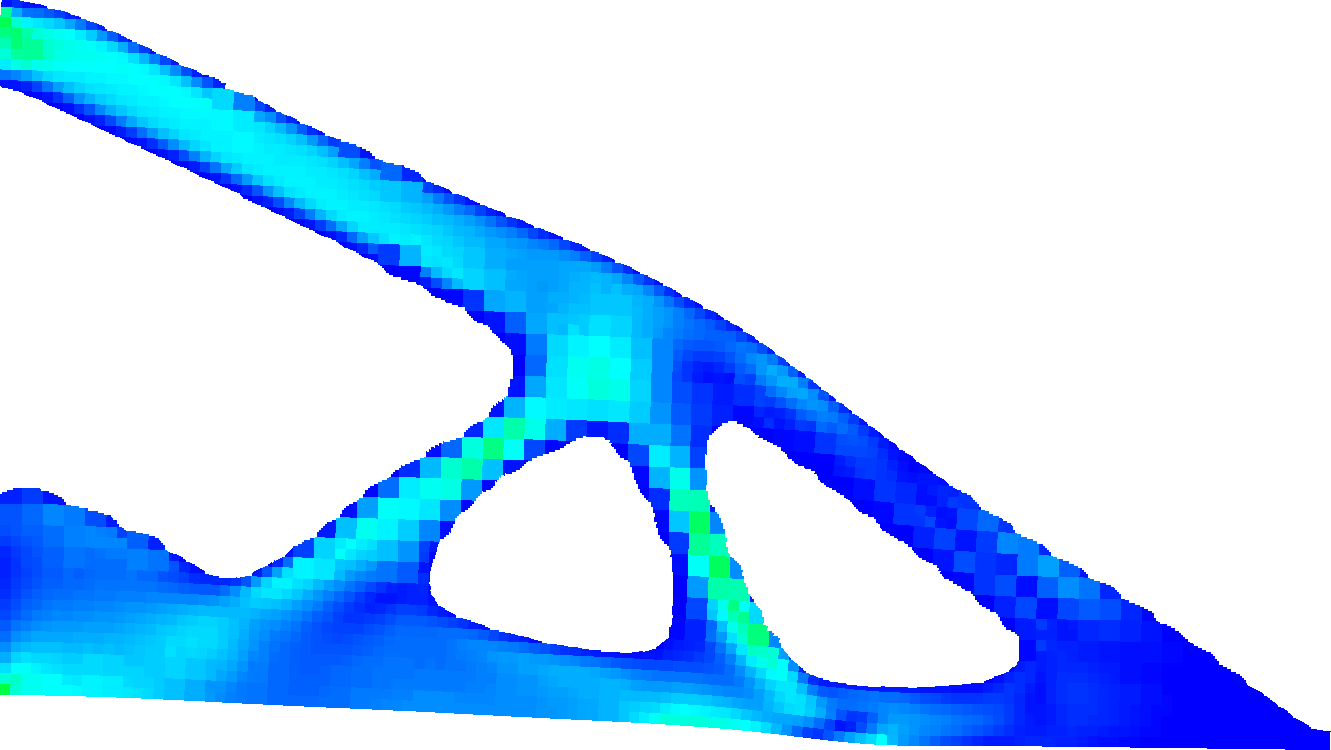}
  \end{subfigure}
  \begin{subfigure}[t]{0.11\textwidth}
    \includegraphics[width=\textwidth]{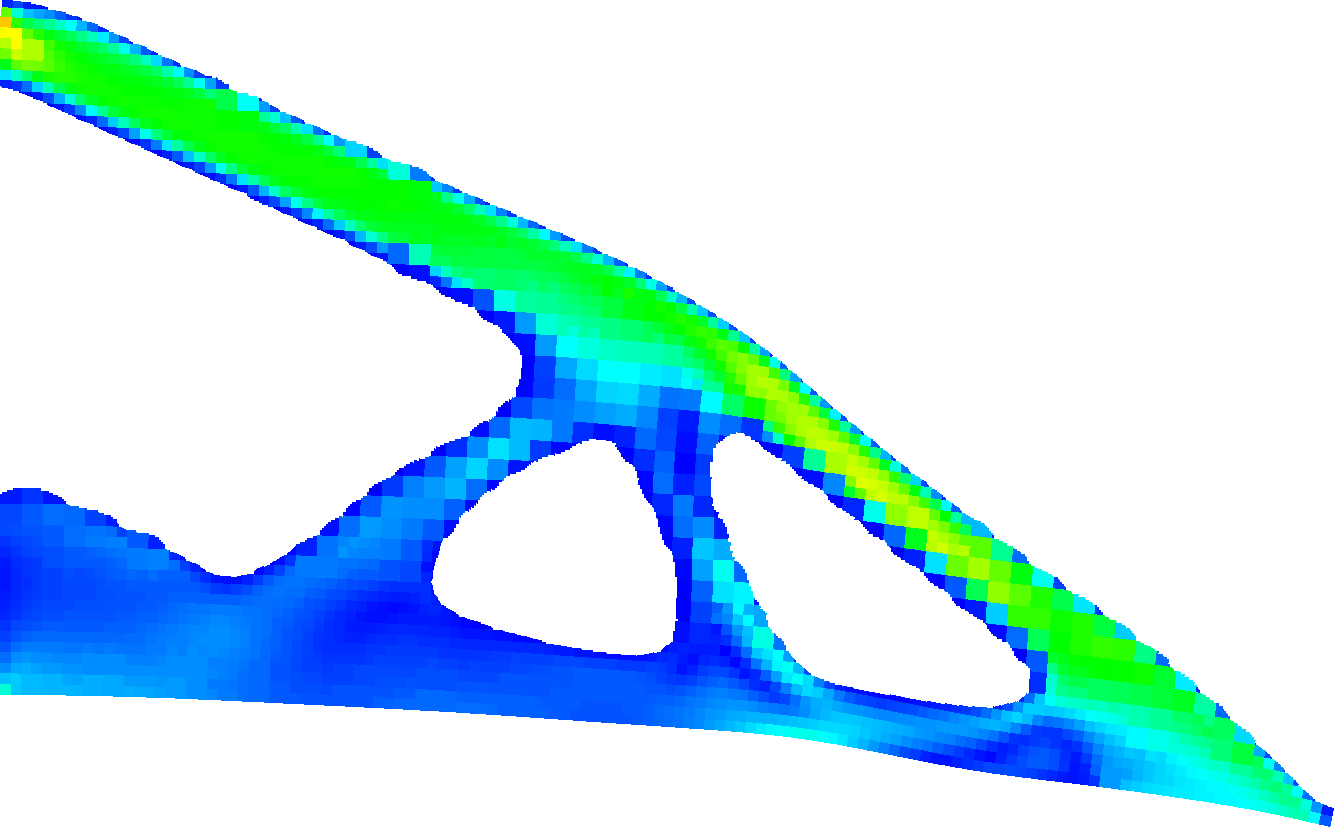}
  \end{subfigure}
  \begin{subfigure}[t]{0.11\textwidth}
    \includegraphics[width=\textwidth]{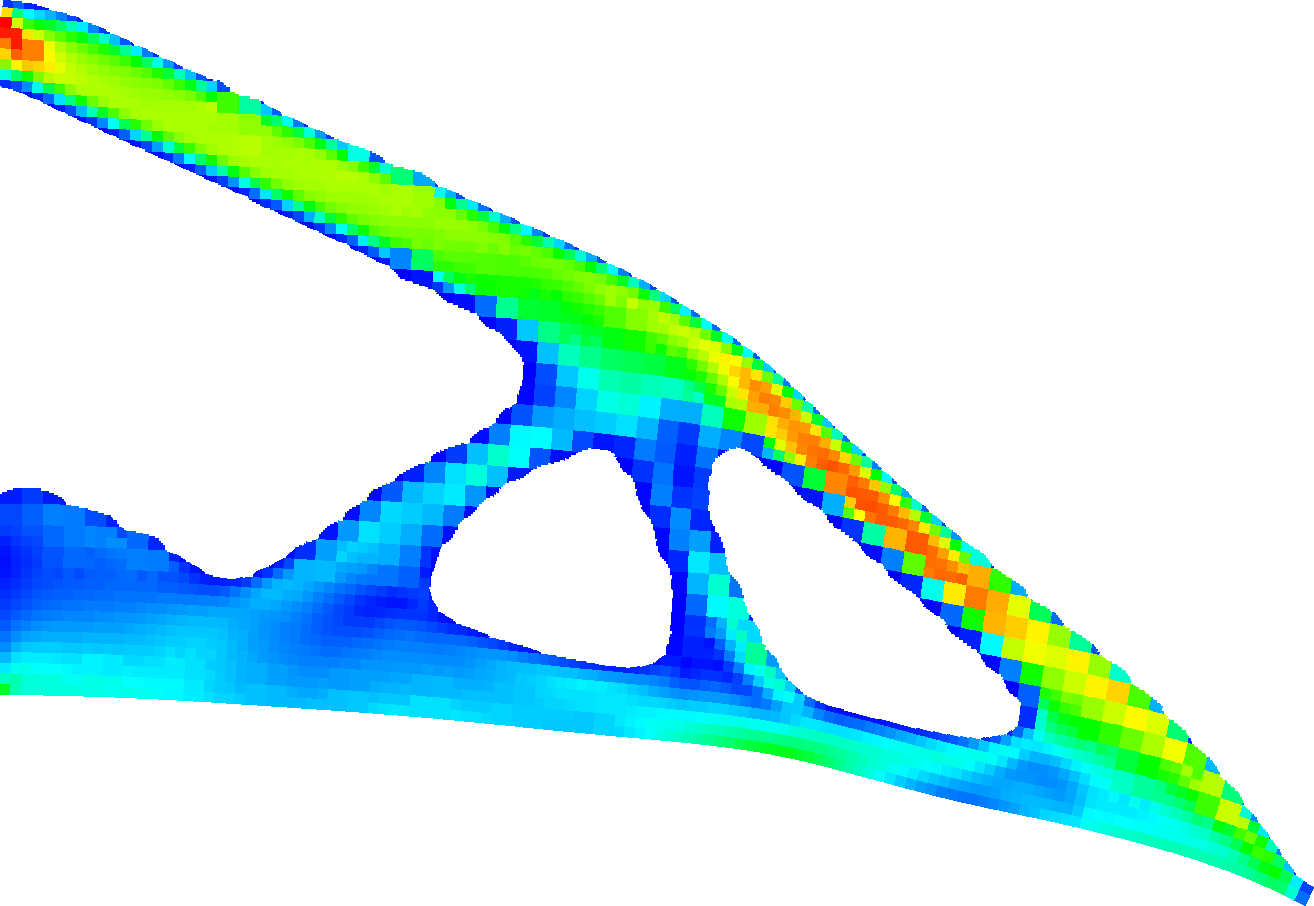}
  \end{subfigure}
  \begin{subfigure}[t]{0.11\textwidth}
    \includegraphics[width=\textwidth]{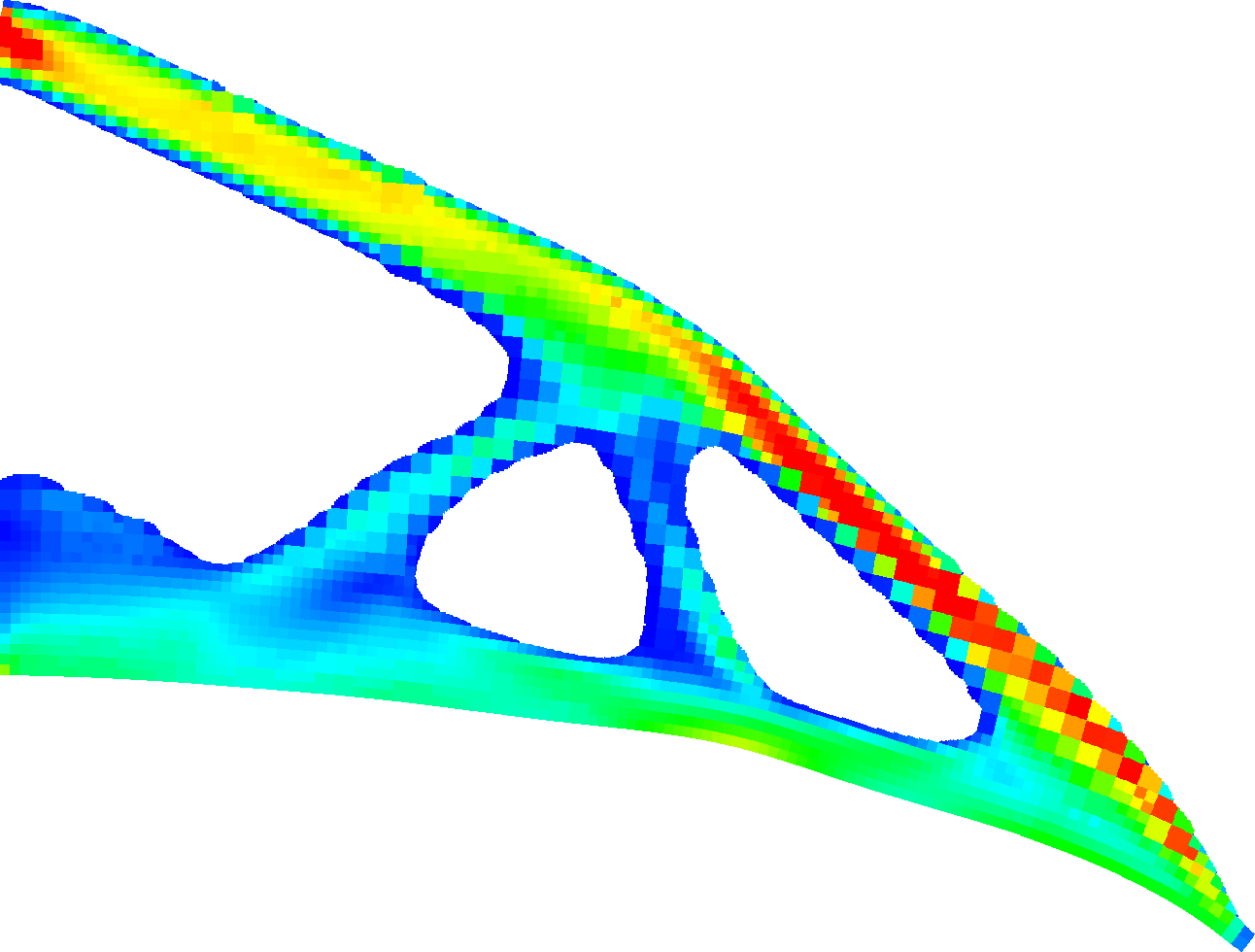}
  \end{subfigure}
  \begin{subfigure}[t]{0.11\textwidth}
    \includegraphics[width=\textwidth]{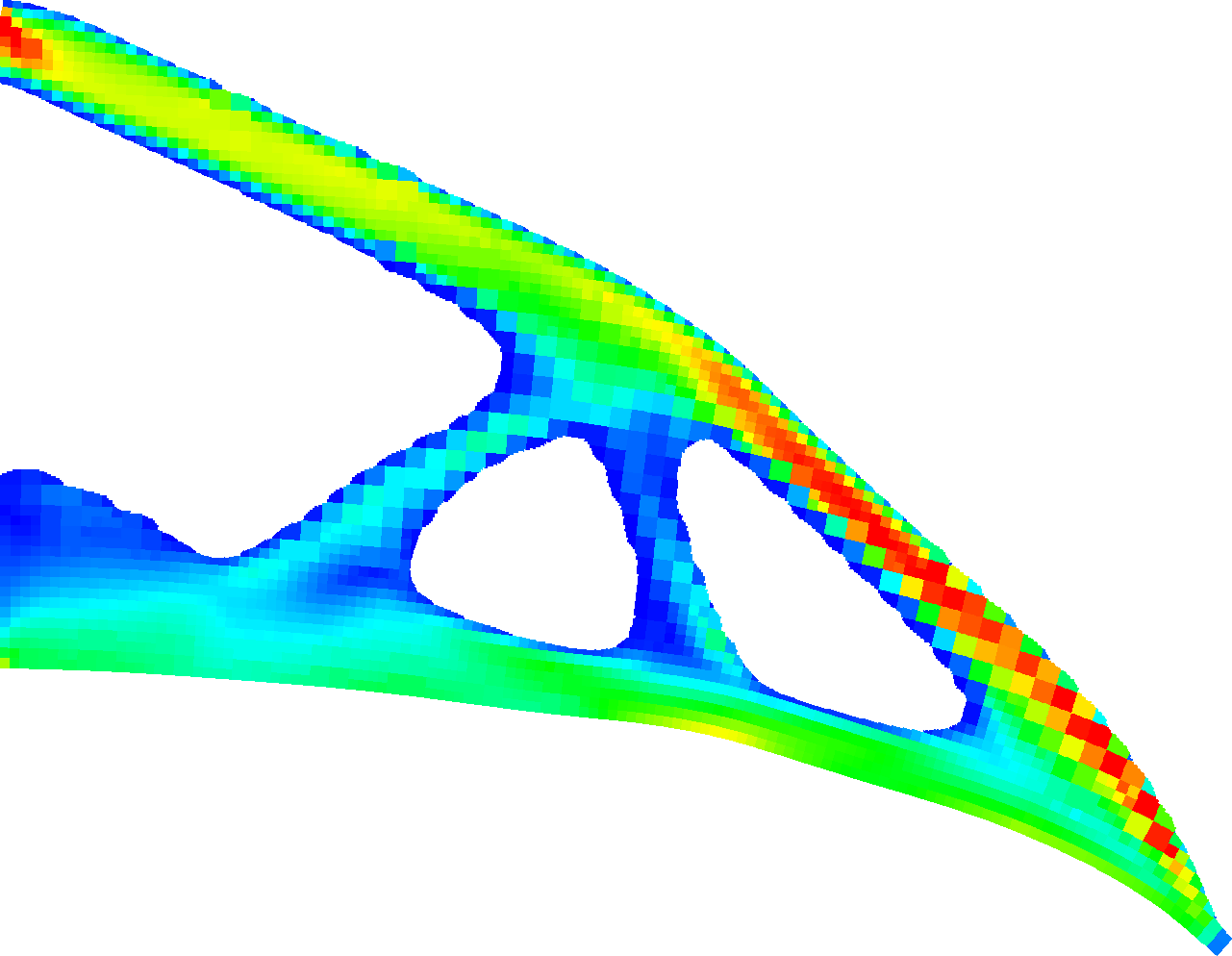}
  \end{subfigure}
  \begin{subfigure}[t]{0.11\textwidth}
    \includegraphics[width=\textwidth]{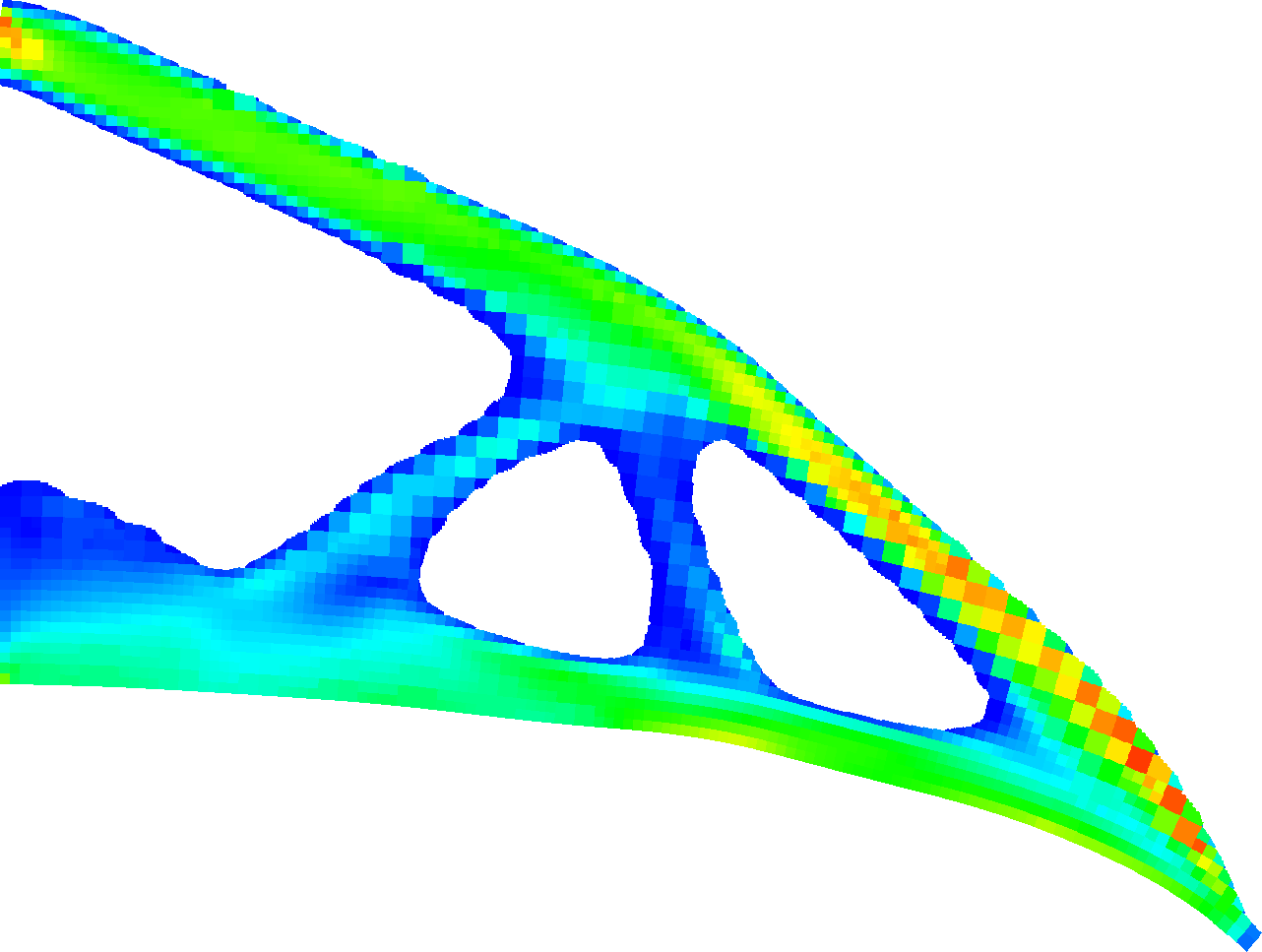}
  \end{subfigure}
  \begin{subfigure}[t]{0.11\textwidth}
    \hspace{\textwidth}
  \end{subfigure}\\[1ex]
   \begin{subfigure}[t]{0.05\textwidth}
    \begin{tikzpicture}
      \node[rotate=90] at (0.0,0.0) {\footnotesize $2^{nd}$ order};
    \end{tikzpicture}
  \end{subfigure}
  \begin{subfigure}[t]{0.11\textwidth}
    \includegraphics[width=\textwidth]{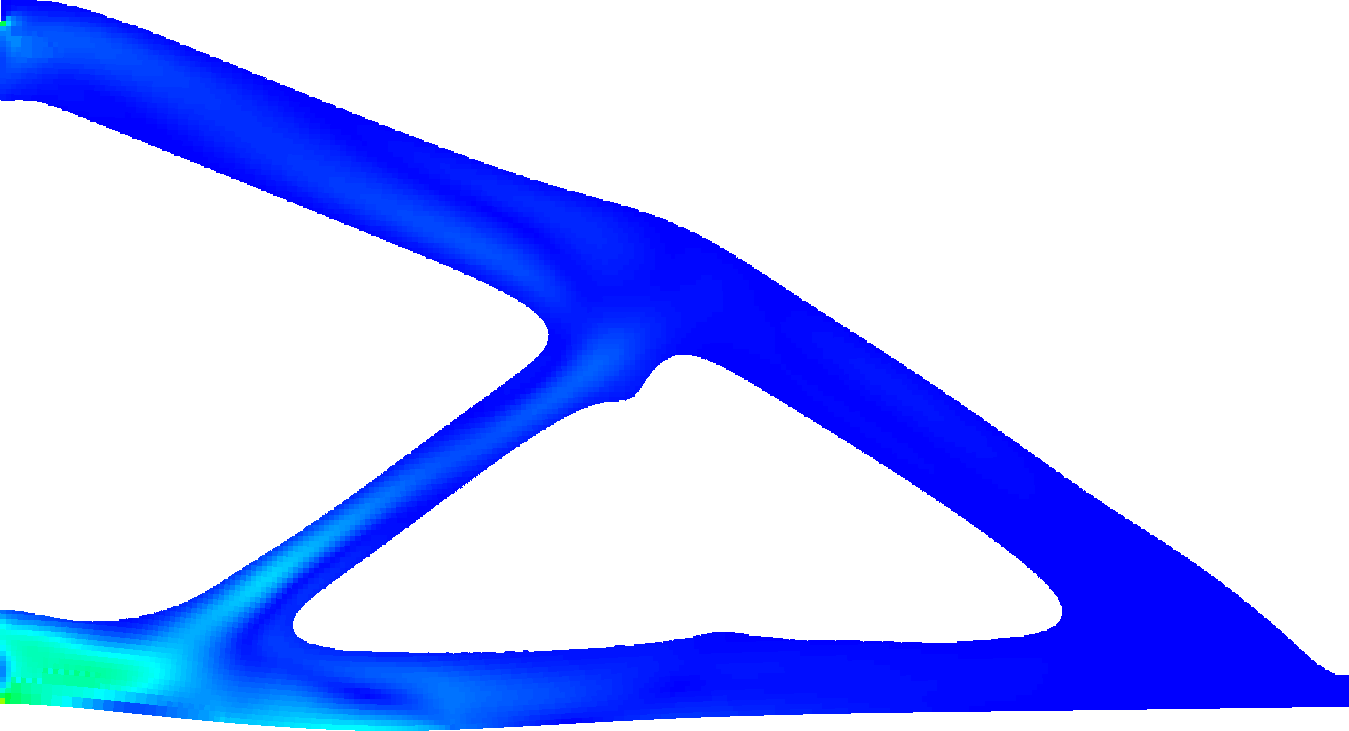}
  \end{subfigure}
  \begin{subfigure}[t]{0.11\textwidth}
    \includegraphics[width=\textwidth]{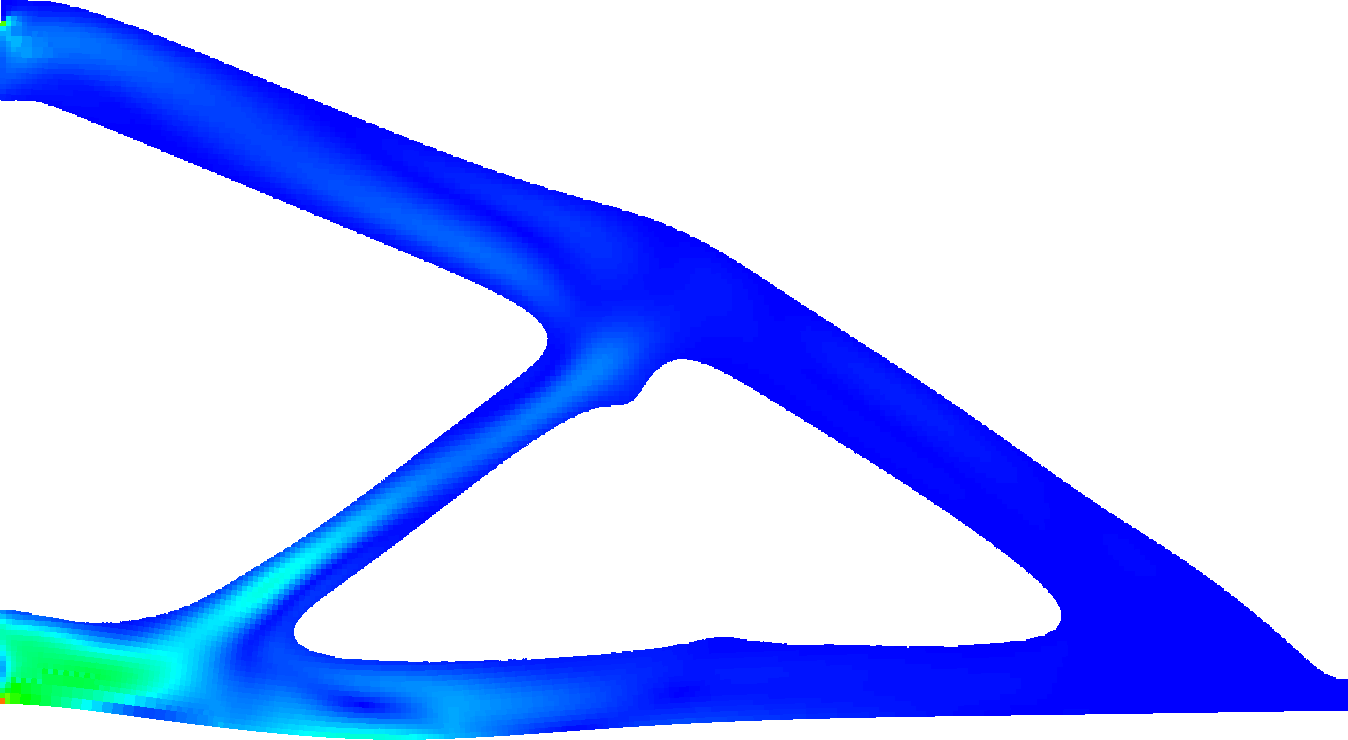}
  \end{subfigure}
  \begin{subfigure}[t]{0.11\textwidth}
    \includegraphics[width=\textwidth]{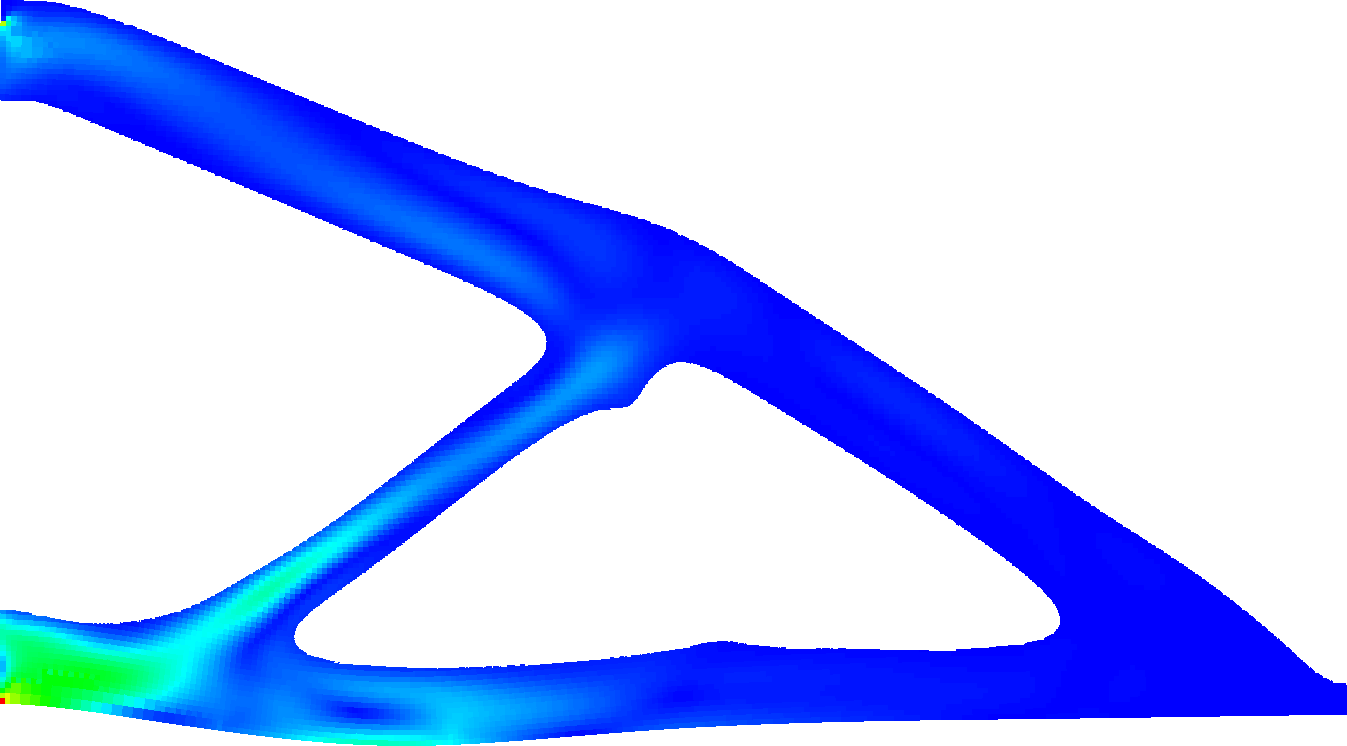}
  \end{subfigure}
  \begin{subfigure}[t]{0.11\textwidth}
    \includegraphics[width=\textwidth]{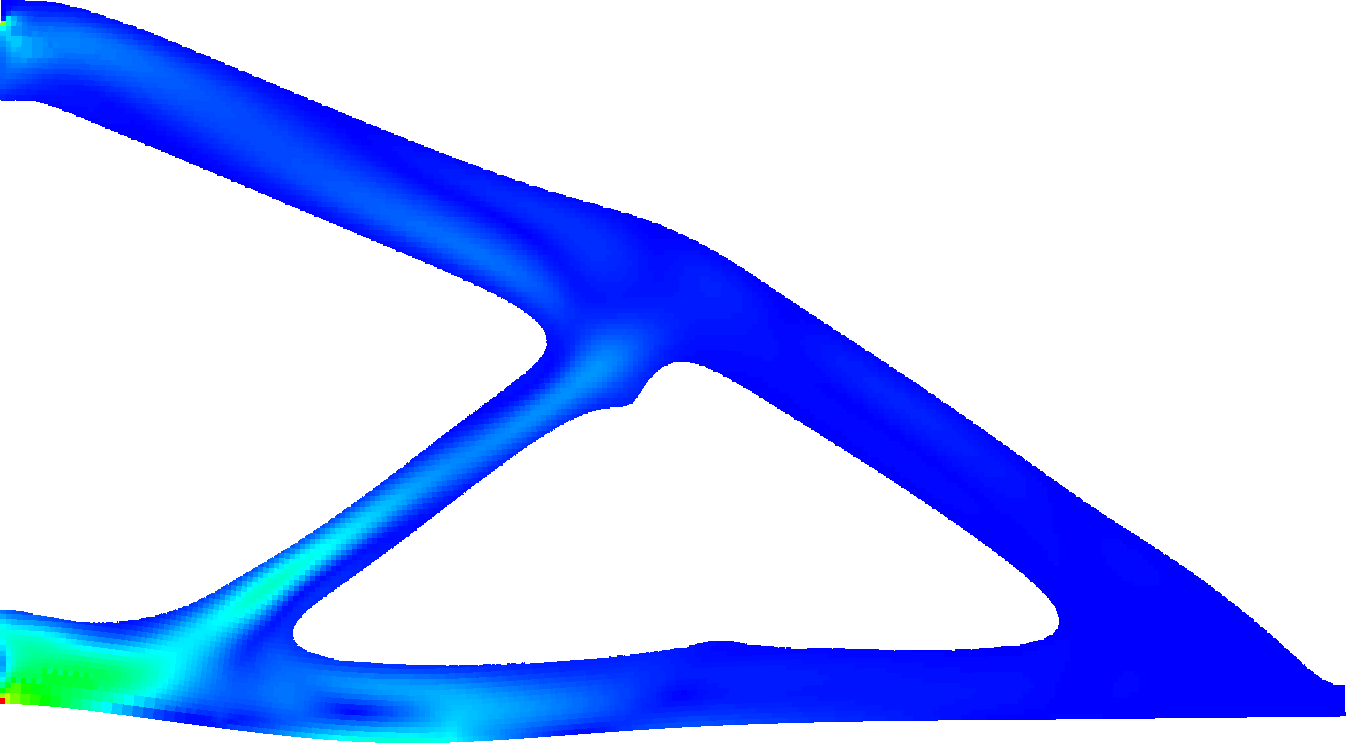}
  \end{subfigure}
  \begin{subfigure}[t]{0.11\textwidth}
    \includegraphics[width=\textwidth]{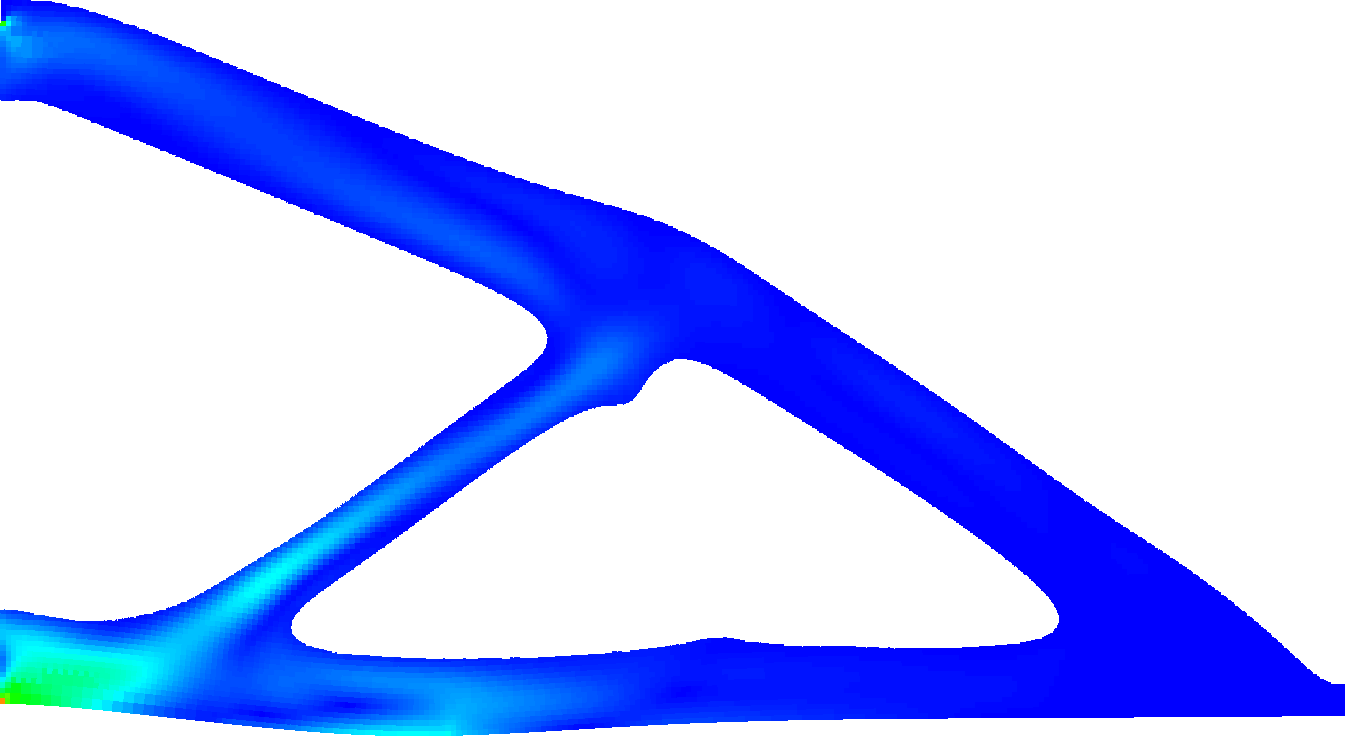}
  \end{subfigure}
  \begin{subfigure}[t]{0.11\textwidth}
    \includegraphics[width=\textwidth]{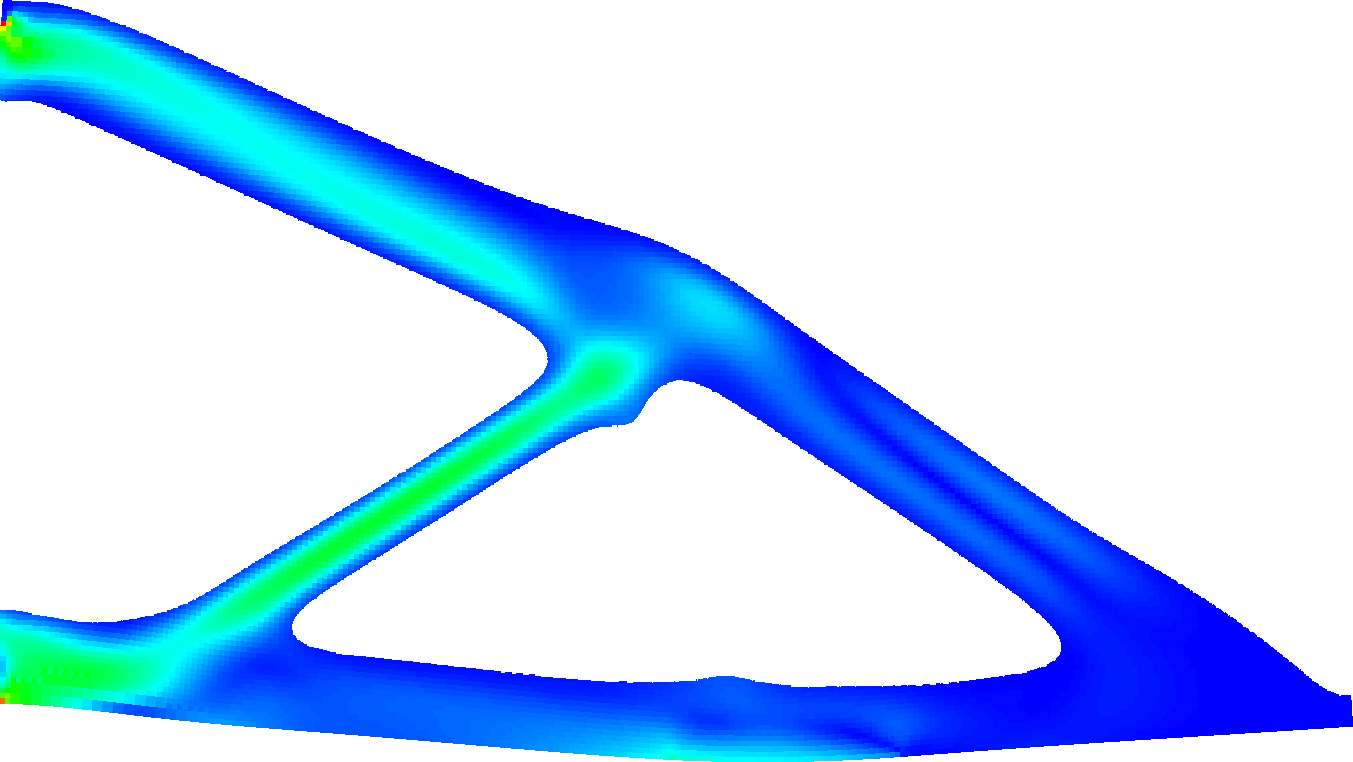}
  \end{subfigure}
  \begin{subfigure}[t]{0.11\textwidth}
    \includegraphics[width=\textwidth]{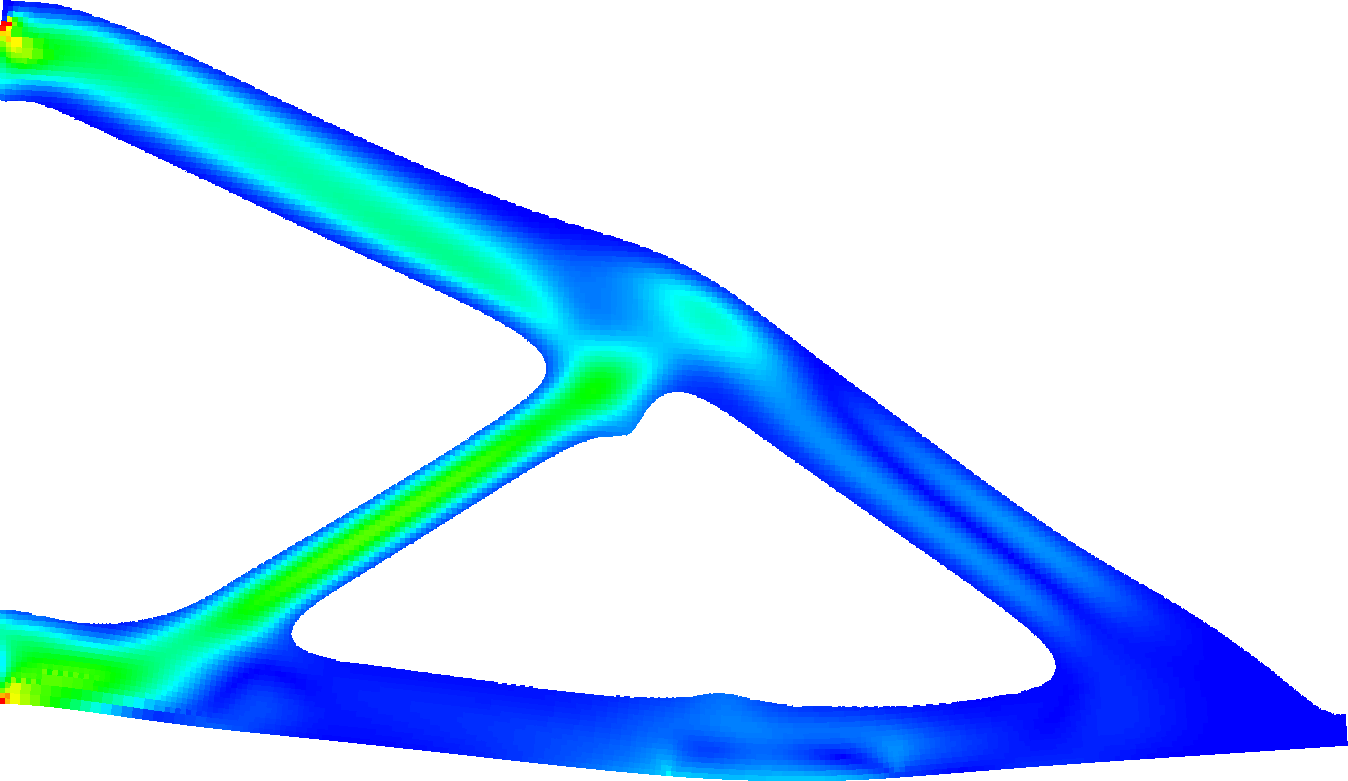}
  \end{subfigure}
  \begin{subfigure}[t]{0.11\textwidth}
    \includegraphics[width=\textwidth]{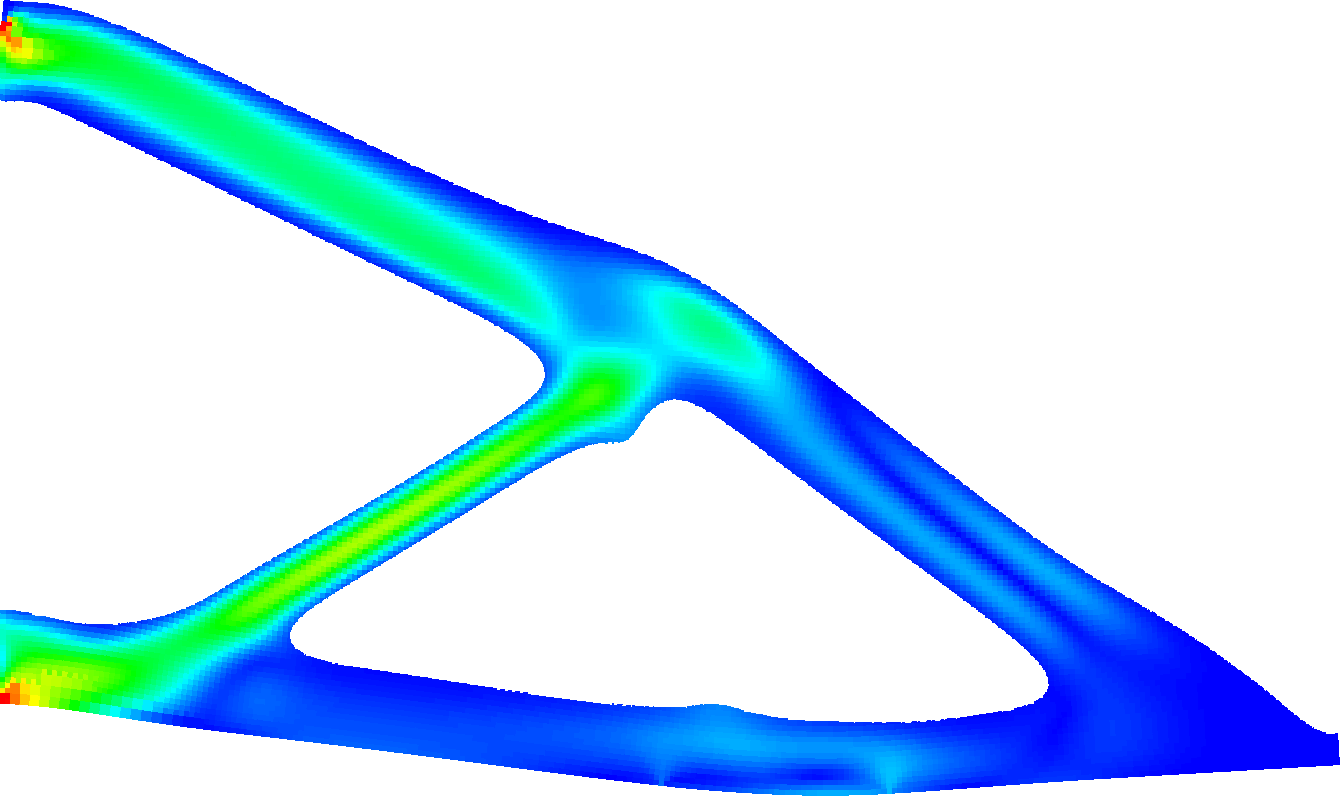}
  \end{subfigure}\\[1ex]
  \begin{subfigure}[t]{0.05\textwidth}
    \begin{tikzpicture}
      \node[rotate=90,align=left] at (0.0,0.0) {\footnotesize var. loads\\
      \footnotesize and prob.};
    \end{tikzpicture}
  \end{subfigure}
  \begin{subfigure}[t]{0.11\textwidth}
    \includegraphics[width=\textwidth]{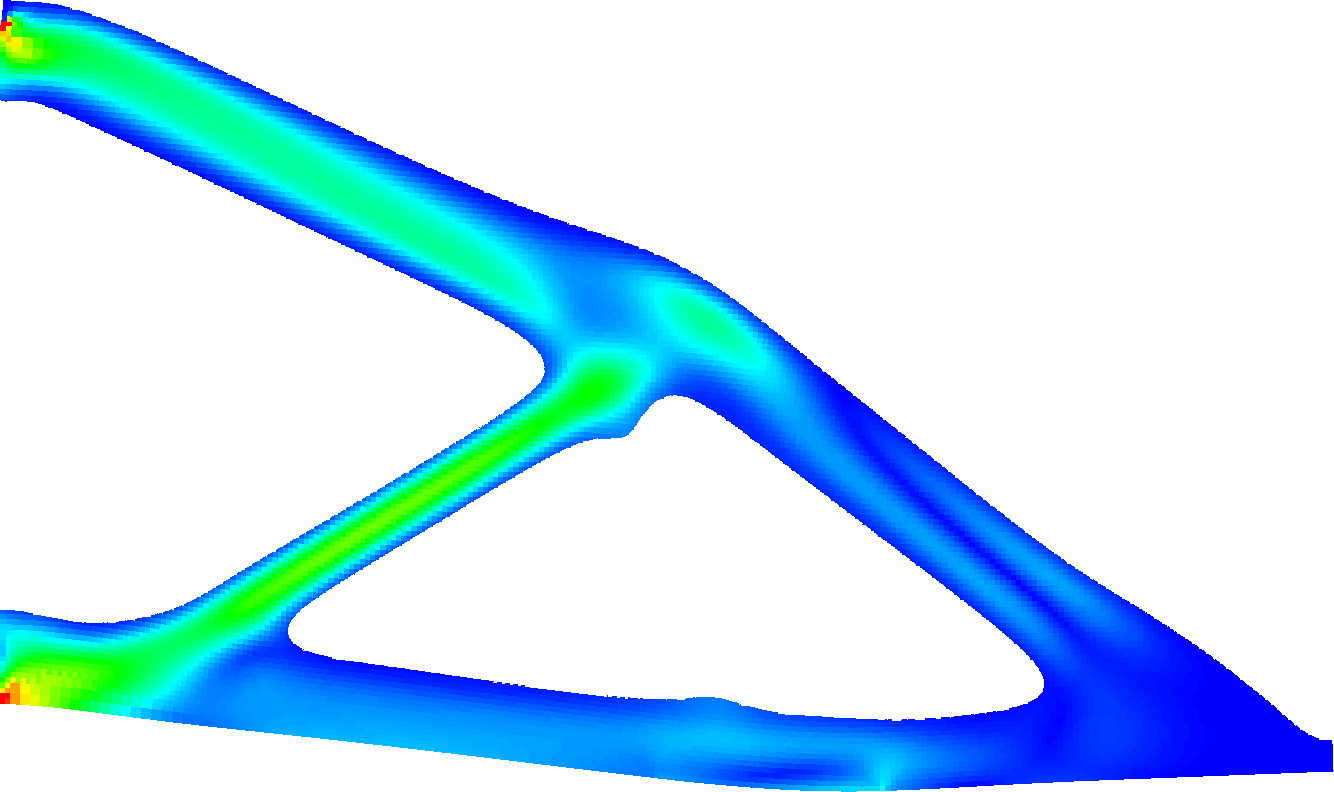}
  \end{subfigure}
  \begin{subfigure}[t]{0.11\textwidth}
    \includegraphics[width=\textwidth]{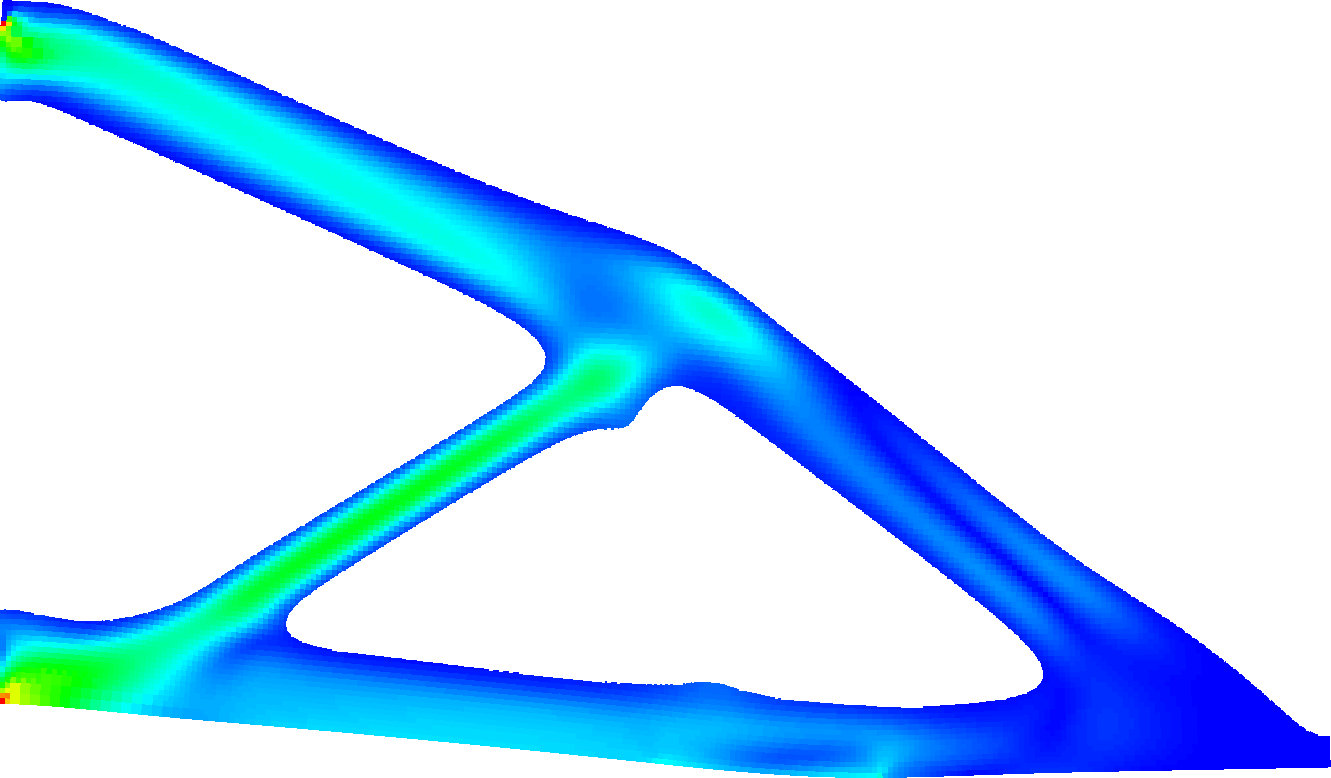}
  \end{subfigure}
  \begin{subfigure}[t]{0.11\textwidth}
    \includegraphics[width=\textwidth]{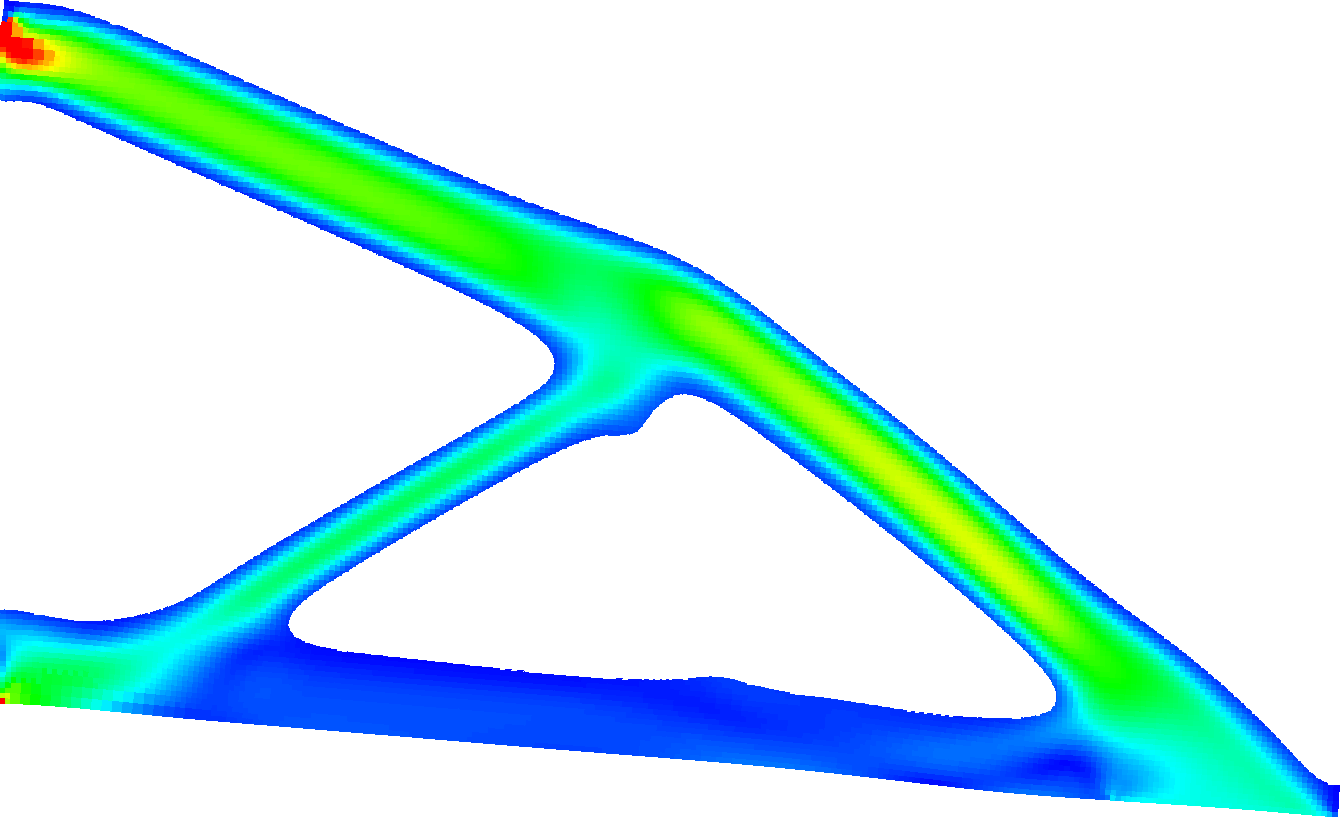}
  \end{subfigure}
  \begin{subfigure}[t]{0.11\textwidth}
    \includegraphics[width=\textwidth]{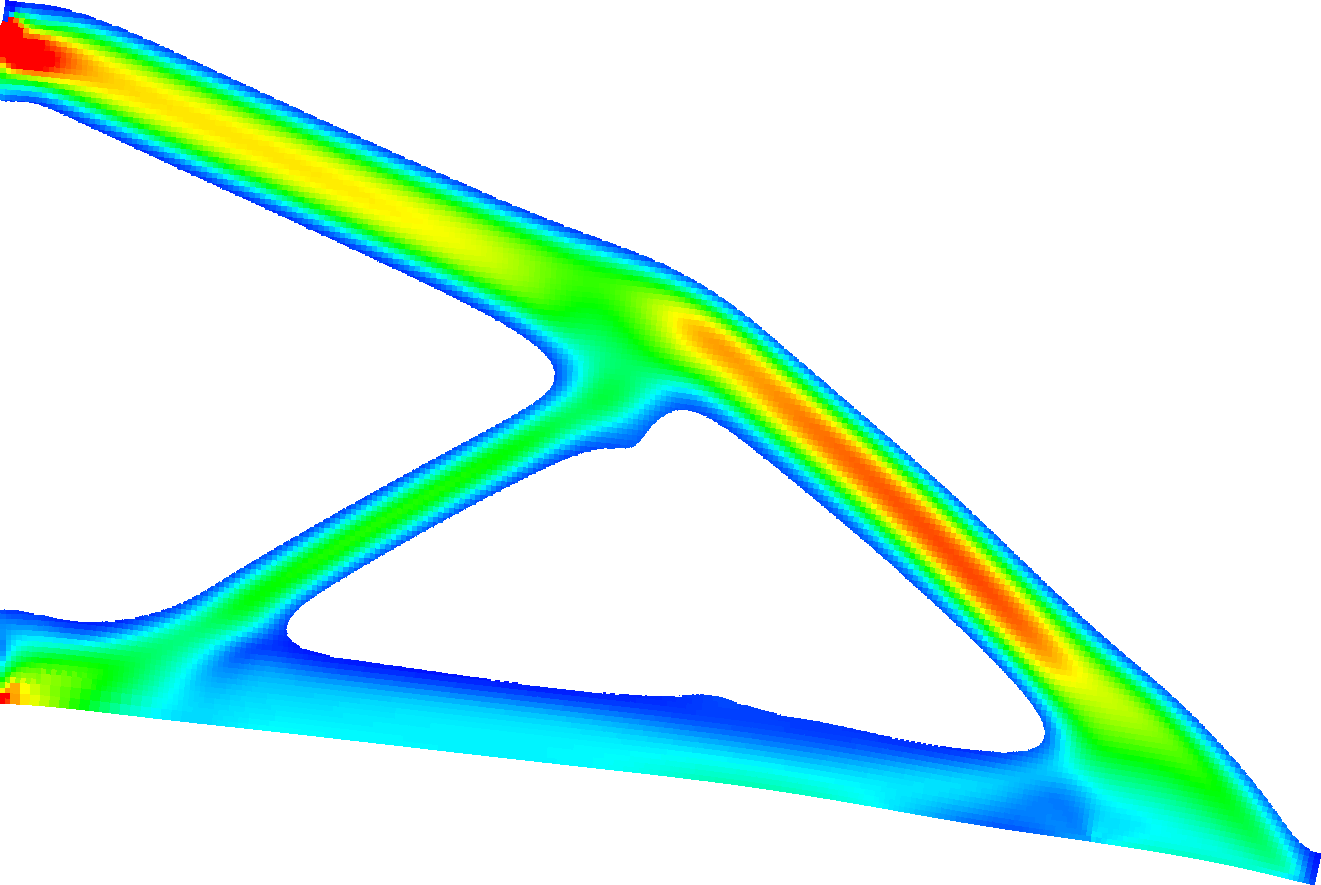}
  \end{subfigure}
  \begin{subfigure}[t]{0.11\textwidth}
    \includegraphics[width=\textwidth]{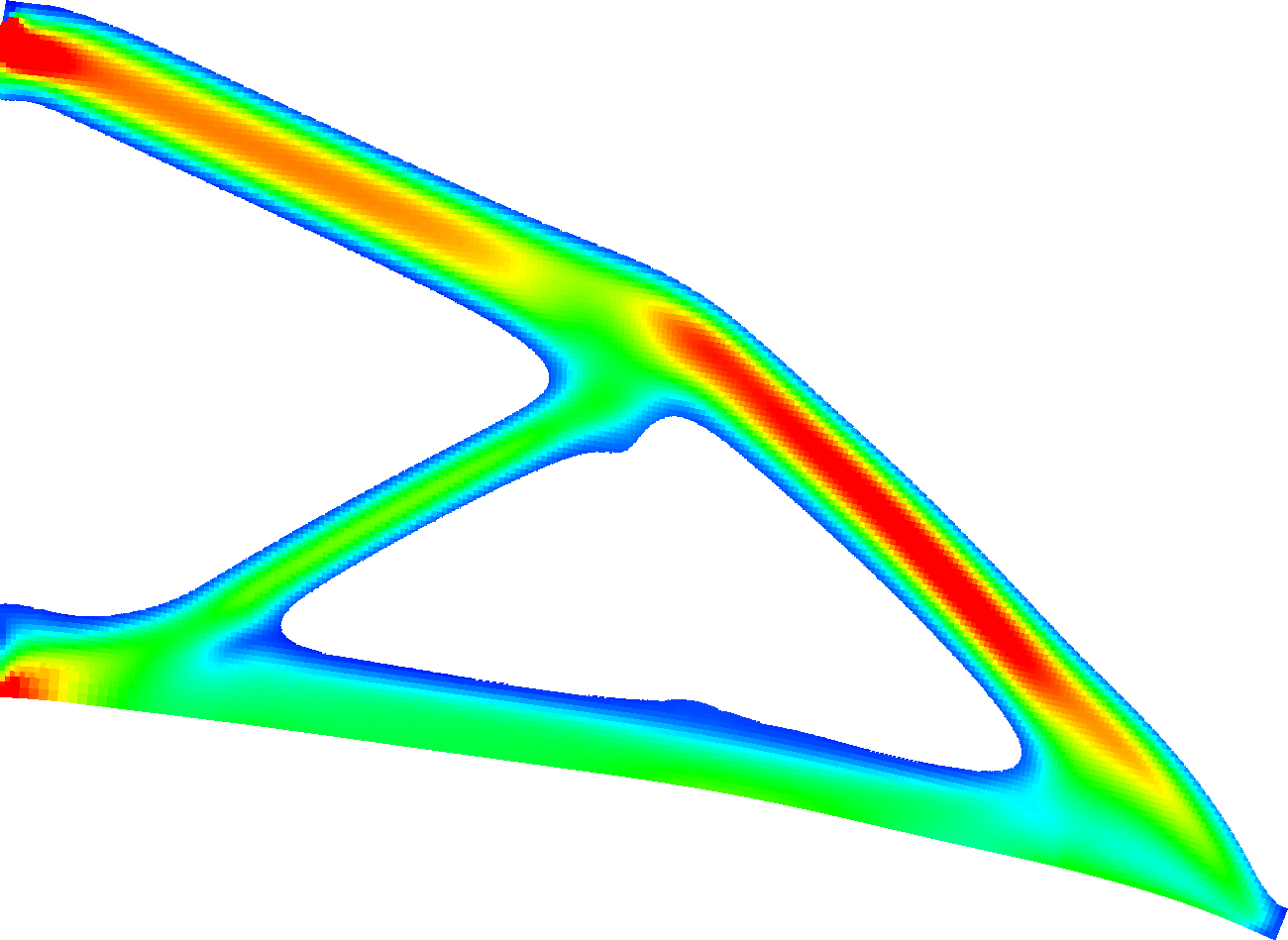}
  \end{subfigure}
  \begin{subfigure}[t]{0.11\textwidth}
    \includegraphics[width=\textwidth]{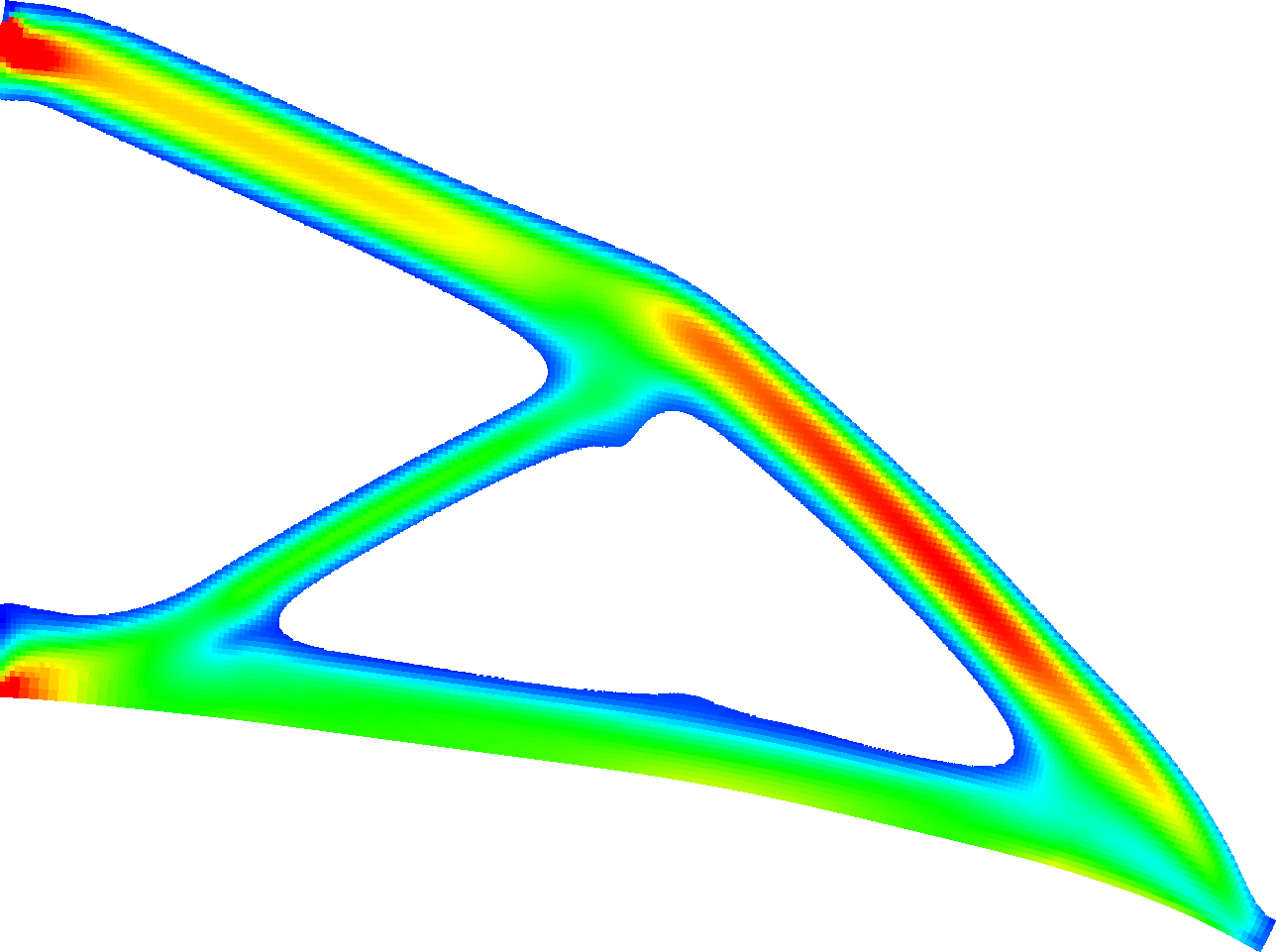}
  \end{subfigure}
  \begin{subfigure}[t]{0.11\textwidth}
    \includegraphics[width=\textwidth]{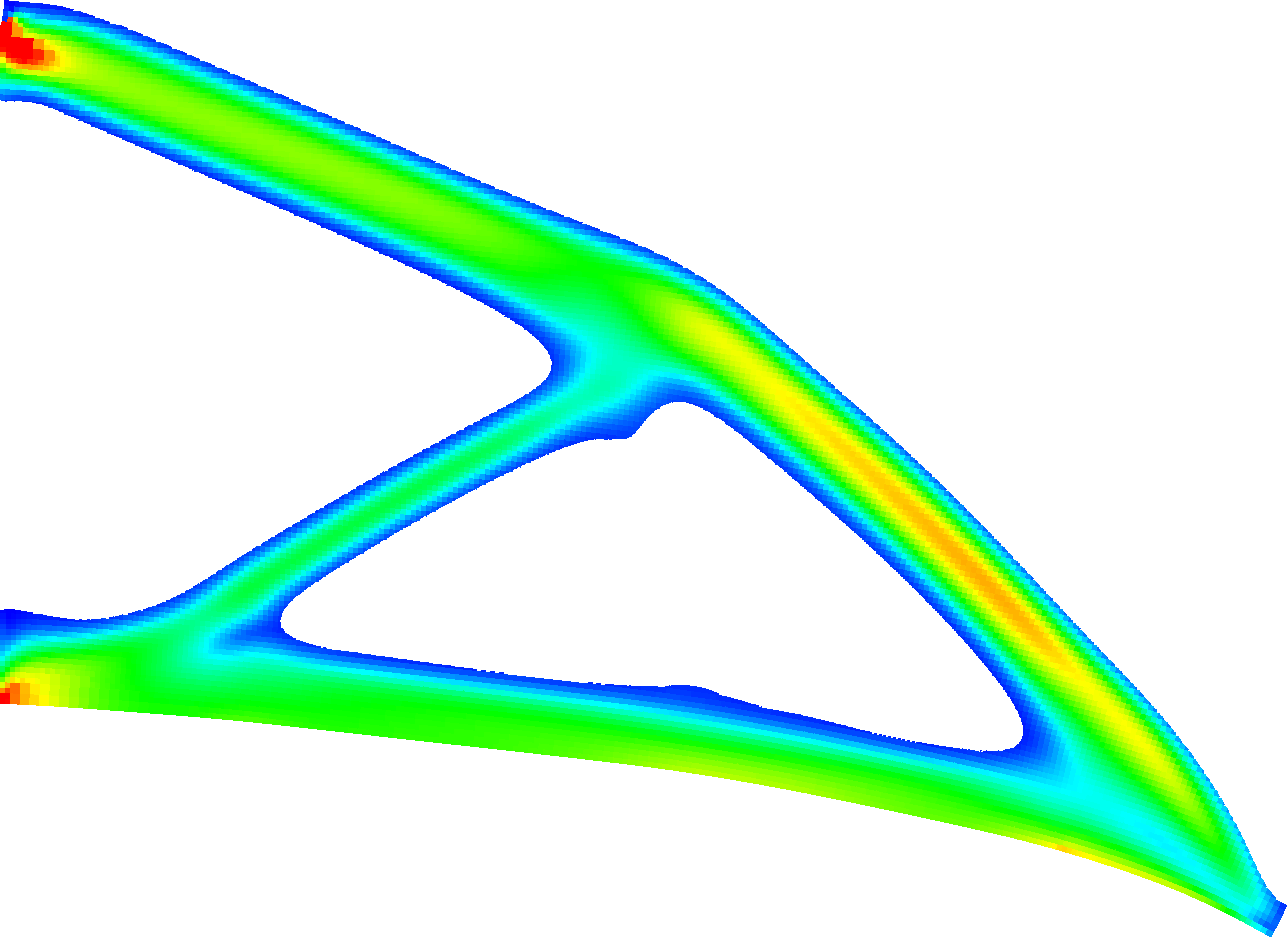}
  \end{subfigure}
  \begin{subfigure}[t]{0.11\textwidth}
    \hspace{\textwidth}
  \end{subfigure}
  \caption{Thresholded stresses for each scenario in the cantilever
    setup, $1^{st}$ and $2^{nd}$ order dominance for the equal setup (row one and two, 
    color-coded as
    $0$~\protect\includegraphics[height=1ex,width=5em]{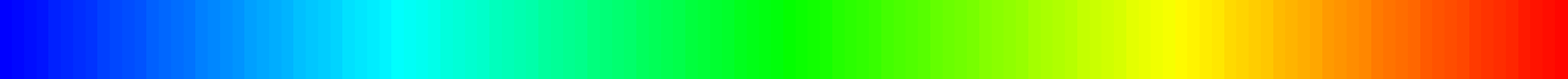}~$4.99$ and
    $0$~\protect\includegraphics[height=1ex,width=5em]{images/res/cb_stress}~$4.85$) and the
    varying load and varying probability setup (row three and four, colorcoded
    as $0$~\protect\includegraphics[height=1ex,width=5em]{images/res/cb_stress}~$9.05$
    and $0$~\protect\includegraphics[height=1ex,width=5em]{images/res/cb_stress}~$6.87$).
    Stresses $> 4.99$,$> 4.85$, $> 9.05$ and $> 6.87$ are mapped to red.}
  \label{fig:stresscanti}
\end{figure}

We compare here a setup consisting of 15 scenarios with equal absolute value of the load and equal probability (top left and first row of diagrams) and a setup 
consisting of 15 scenarios with three different absolute value  of the load (ratios $1$, $\frac23$ and $\frac13$) and different probabilities. In fact, the smaller loads are $2$ and $3$ times more probable than the 
larger ones (top right and second row of diagrams). 

The achieved volumes are listed underneath the corresponding phase field of the optimal shapes. The diagrams show the cumulative distribution function for the benchmark profile (blue)
for the first order dominance model (left) and the second order dominance model (right). In the middle the integrated survival function for the second order dominance model is shown.
Due to \eqref{eq:cdf} the benchmark values have to be smaller than or equal to the corresponding values for the shape, which is optimal with respect to first order dominance (left).
Furthermore, because of  \eqref{eq:isf} the benchmark values have to be larger than or equal to the corresponding values for the shape, 
which is optimal with respect to second order dominance (middle). Both properties are clearly reflected by the computational diagrams.
For second order dominance, values of the cumulative distribution function can be above the corresponding values for the benchmark geometry for smaller $t$, 
 as long as this is  compensated for  by larger $t$ values, in the appropriate integral sense.
This effect is substantially stronger for the configuration with varying loads and varying probabilities (diagrams in the third row).
Finally, Fig.~\ref{fig:stresscanti} displays the distribution of von Mises stresses
on the hard material phase in the 
deformed configuration for all 15 scenarios for varying and equal probability and loads,
first order and second order optimality, respectively.

\paragraph{Carrier plate.}
As a second application we consider a 2D carrier plate, where the supporting construction between a 
floor slap, whose lower boundary is assumed to be the Dirichlet boundary, and an upper plate, on which forces act, is optimized. 
In this example the stochastic loading consists of 10 scenarios consisting of 10 single loads.

\begin{figure}[htb]
  \centering
  \begin{subfigure}[t]{0.3\textwidth}
    \includegraphics[width=\textwidth]{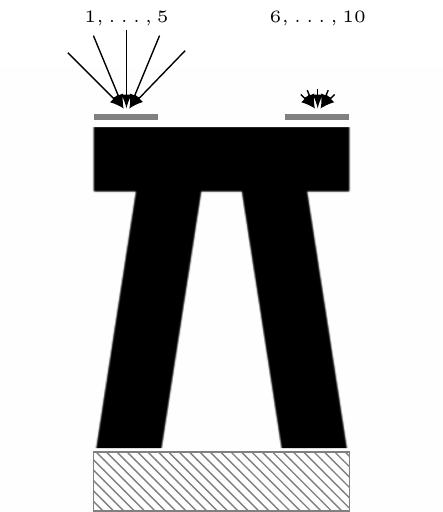}
    \caption{\color{red} Benchmark\\Vol. $0.484881$}
    \label{fig:benchcpw}
  \end{subfigure}
  \begin{subfigure}[t]{0.3\textwidth}
    \includegraphics[width=\textwidth]{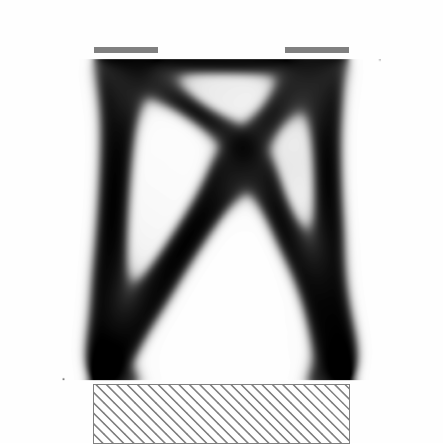}
    \caption{\color{blue!70!black} $1^{st}$ order\\Vol. $0.330986$}
    \label{fig:res1cpw}
  \end{subfigure}
  \begin{subfigure}[t]{0.3\textwidth}
    \includegraphics[width=\textwidth]{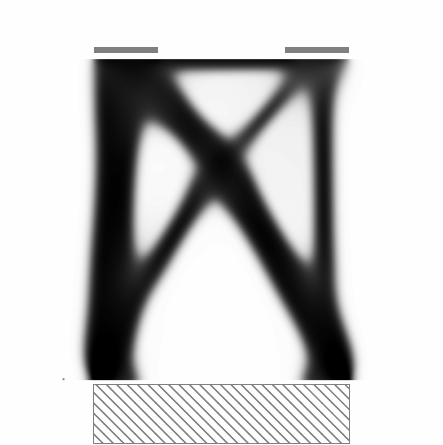}
    \caption{\color{blue!70!black} $2^{nd}$ order\\Vol. $0.12375$}
    \label{fig:res2cpw}
  \end{subfigure}
 \begin{subfigure}[t]{0.32\textwidth}
    \begin{tikzpicture}
     \begin{axis}[axis lines=middle,
        axis line style={->},
        enlargelimits=true,
        width=\textwidth,
        xmin=0, xmax=0.03,
        xlabel={$t$},
        x label style={at={(axis description cs:1.2, -0.1)}},
        xticklabel style={
          /pgf/number format/fixed,
          /pgf/number format/precision=4,
          font=\tiny },
        scaled x ticks=false,
        ylabel={$F(t)$},
        y label style={
          at={(axis description cs:-0.2, 1.3)},
          rotate=0 },
        yticklabel style={
          font=\tiny },
        ]
        \addplot[red]   table[x={x}, y={F}, mark=none] {images/res/cpw/bench_cdf1};
        \addplot[blue!70!black] table[x={x}, y={F}, mark=none] {images/res/cpw/pf_cdf1};
      \end{axis}
    \end{tikzpicture}
    \caption{CDF for $1^{st}$ order dominance}
    \label{fig:cdf1cpw}
  \end{subfigure}
  \begin{subfigure}[t]{0.32\textwidth}
    \begin{tikzpicture}
      \begin{axis}[axis lines=middle,
        axis line style={->},
        enlargelimits=true,
        width=\textwidth,
        xmin=0, xmax=0.03,
        xlabel={$t$},
        x label style={at={(axis description cs:1.2, -0.1)}},
        xticklabel style={
          /pgf/number format/fixed,
          /pgf/number format/precision=4,
          font=\tiny },
        scaled x ticks=false,
        ylabel={$\pi(t)$},
        y label style={
          at={(axis description cs:-0.2, 1.3)},
          rotate=0 },
        yticklabel style={
          font=\tiny },
        ]
        \addplot[red]   table[x={x}, y={isf}, mark=none] {images/res/cpw/bench_isf2};
        \addplot[blue!70!black] table[x={x}, y={isf}, mark=none] {images/res/cpw/pf_isf2};
      \end{axis}
    \end{tikzpicture}
    \caption{ISF for $2^{nd}$ order }
    \label{fig:isf2cpe}
  \end{subfigure}
  \begin{subfigure}[t]{0.32\textwidth}
    \begin{tikzpicture}
\begin{axis}[axis lines=middle,
        axis line style={->},
        enlargelimits=true,
        width=\textwidth,
        xmin=0, xmax=0.03,
        xlabel={$t$},
        x label style={at={(axis description cs:1.2, -0.1)}},
        xticklabel style={
          /pgf/number format/fixed,
          /pgf/number format/precision=4,
          font=\tiny },
        scaled x ticks=false,
        ylabel={$F(t)$},
        y label style={
          at={(axis description cs:-0.2, 1.3)},
          rotate=0 },
        yticklabel style={
          font=\tiny },
        ]
        \addplot[red]   table[x={x}, y={F}, mark=none] {images/res/cpw/bench_cdf2};
        \addplot[blue!70!black] table[x={x}, y={F}, mark=none] {images/res/cpw/pf_cdf2};
      \end{axis}
    \end{tikzpicture}
    \caption{CDF for $2^{nd}$ order dominance}
    \label{fig:cdf2cpw}
  \end{subfigure}
  \caption{Benchmark with load configuration plot, optimal shapes cumulative distance funcions (CDF)
    and integrated survival functions (ISF) for the pressure plate setup.}
  \label{fig:rescp}
\end{figure}
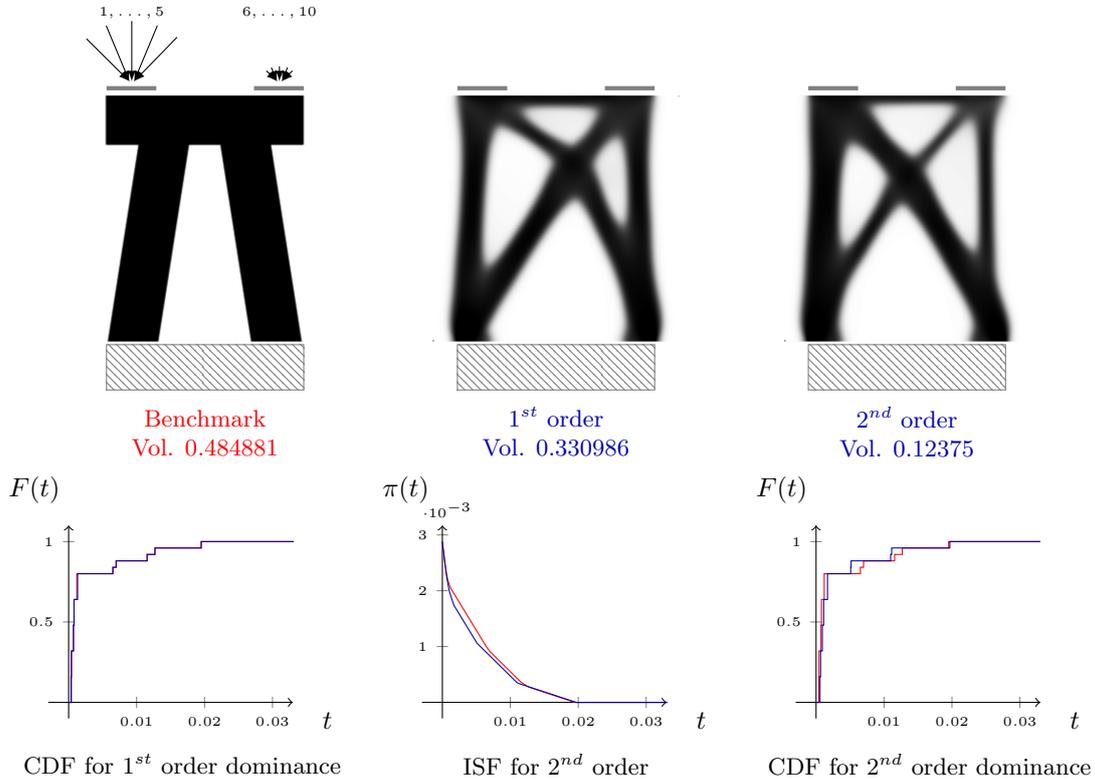

\begin{figure}[htb]
  \label{fig:benchrescpe}
  \vspace{2ex}
  \begin{subfigure}[t]{0.09\textwidth}
    \includegraphics[width=\textwidth]{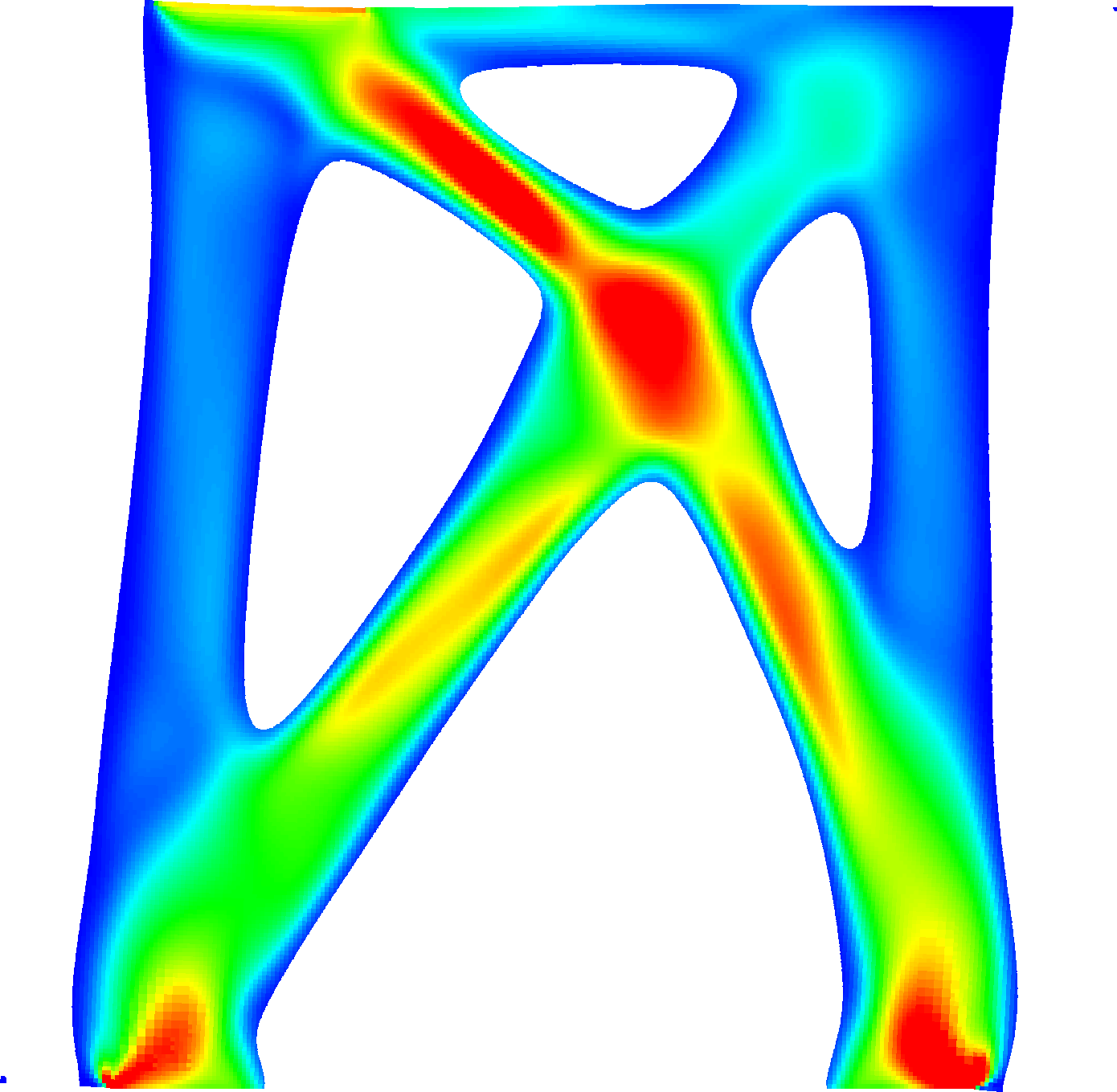}
  \end{subfigure}
  \begin{subfigure}[t]{0.09\textwidth}
    \includegraphics[width=\textwidth]{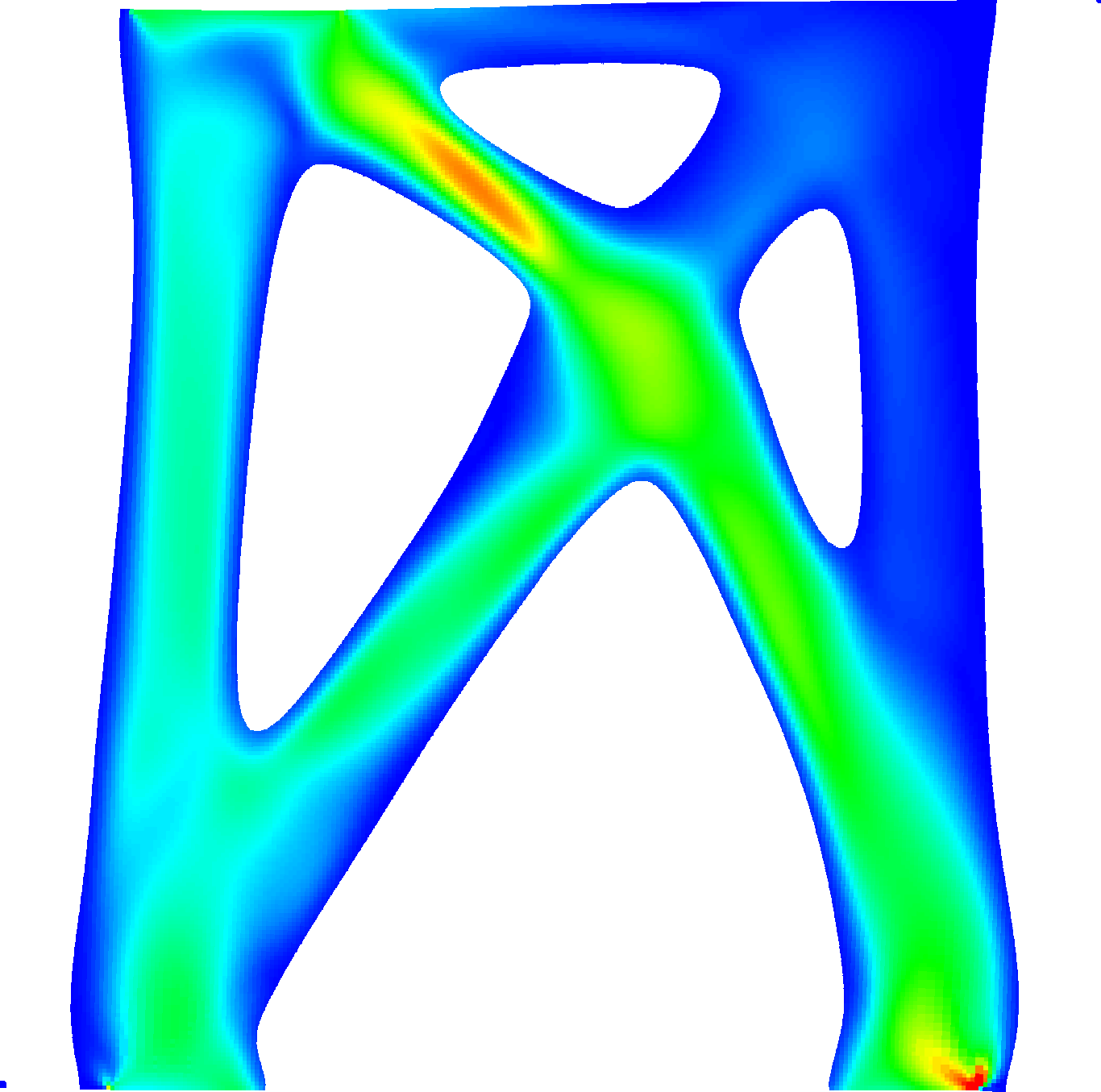}
  \end{subfigure}
  \begin{subfigure}[t]{0.09\textwidth}
    \includegraphics[width=\textwidth]{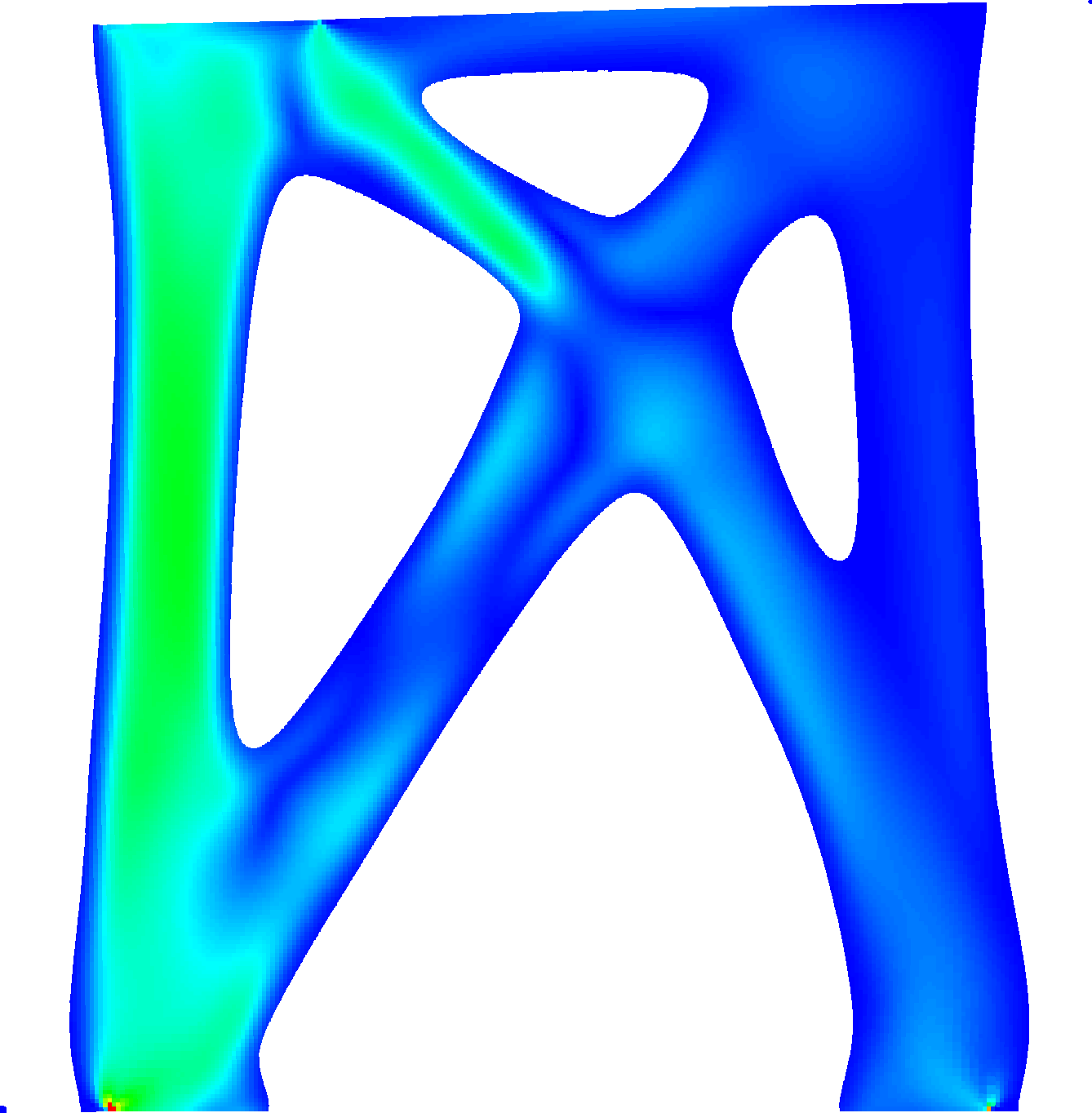}
  \end{subfigure}
  \begin{subfigure}[t]{0.09\textwidth}
    \includegraphics[width=\textwidth]{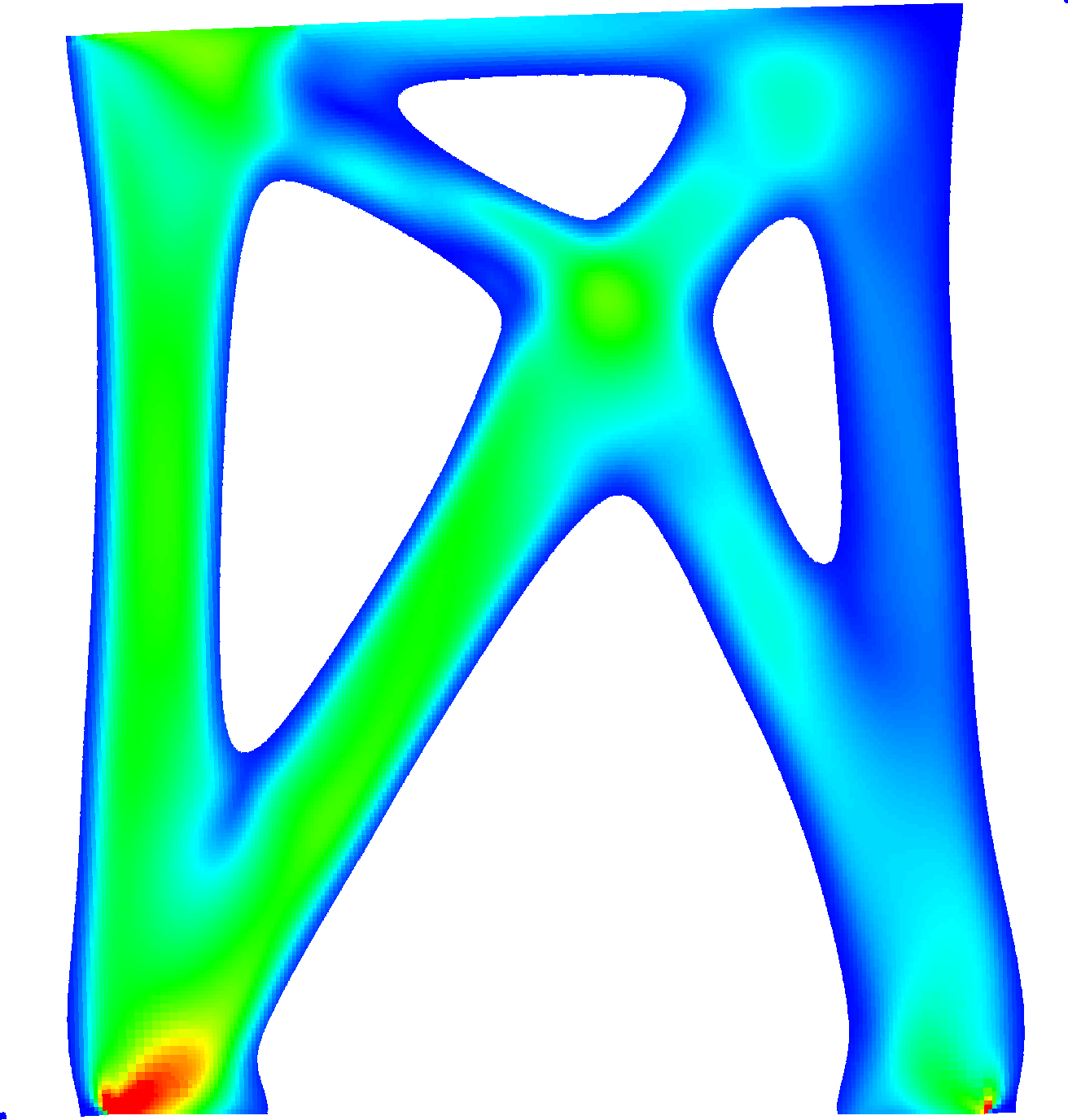}
  \end{subfigure}
  \begin{subfigure}[t]{0.09\textwidth}
    \includegraphics[width=\textwidth]{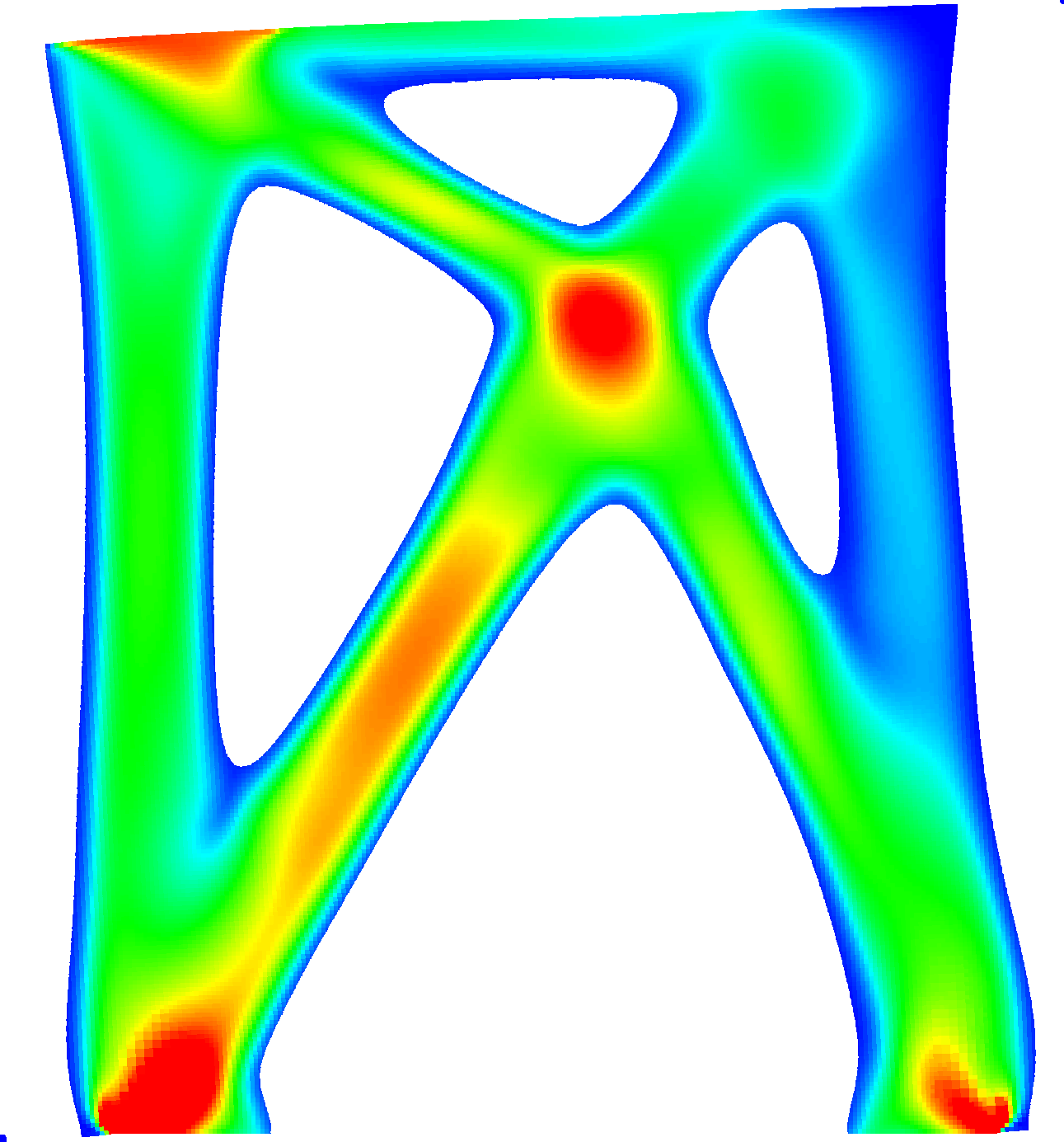}
  \end{subfigure}
  \begin{subfigure}[t]{0.09\textwidth}
    \includegraphics[width=\textwidth]{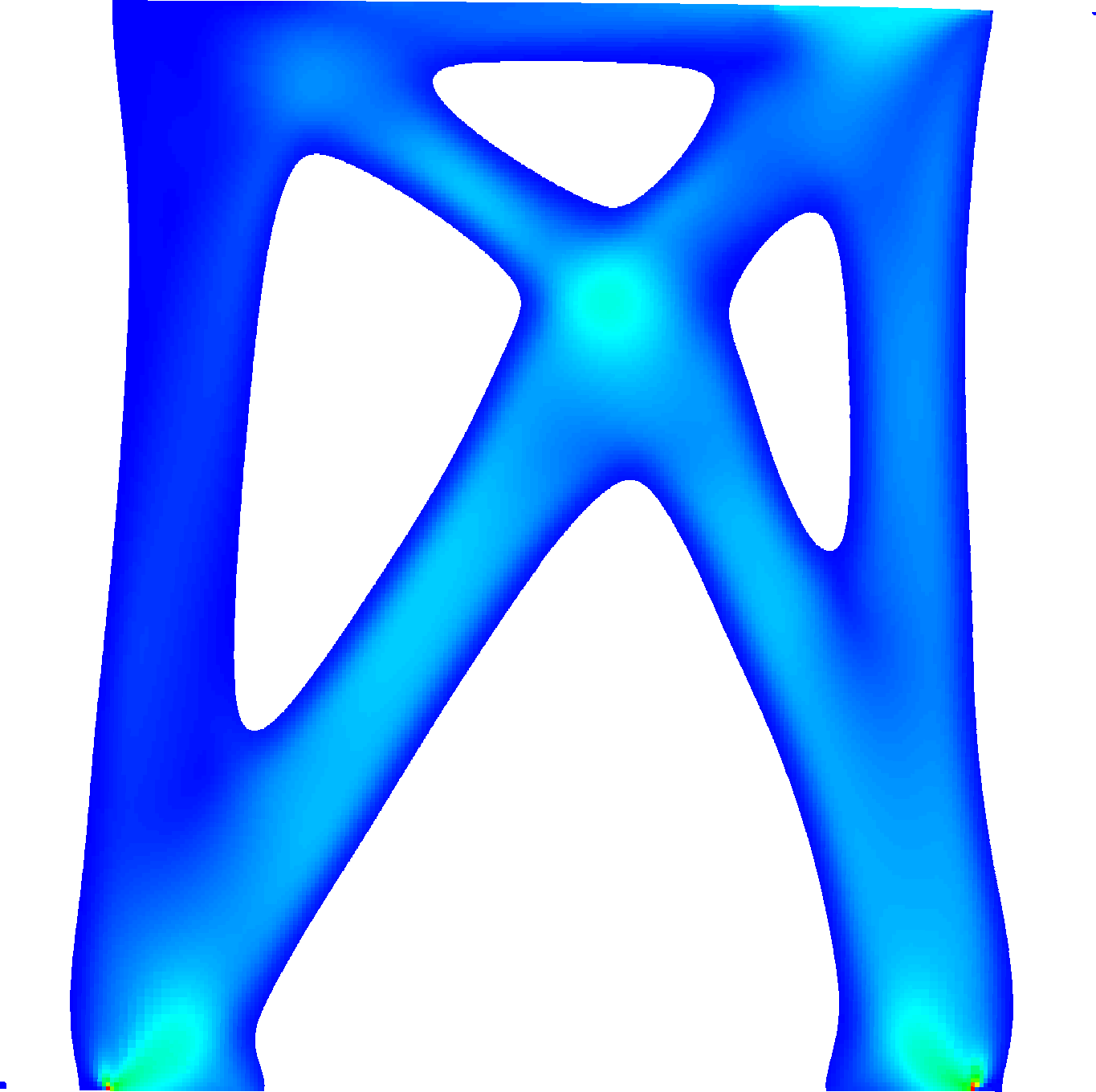}
  \end{subfigure}
  \begin{subfigure}[t]{0.09\textwidth}
    \includegraphics[width=\textwidth]{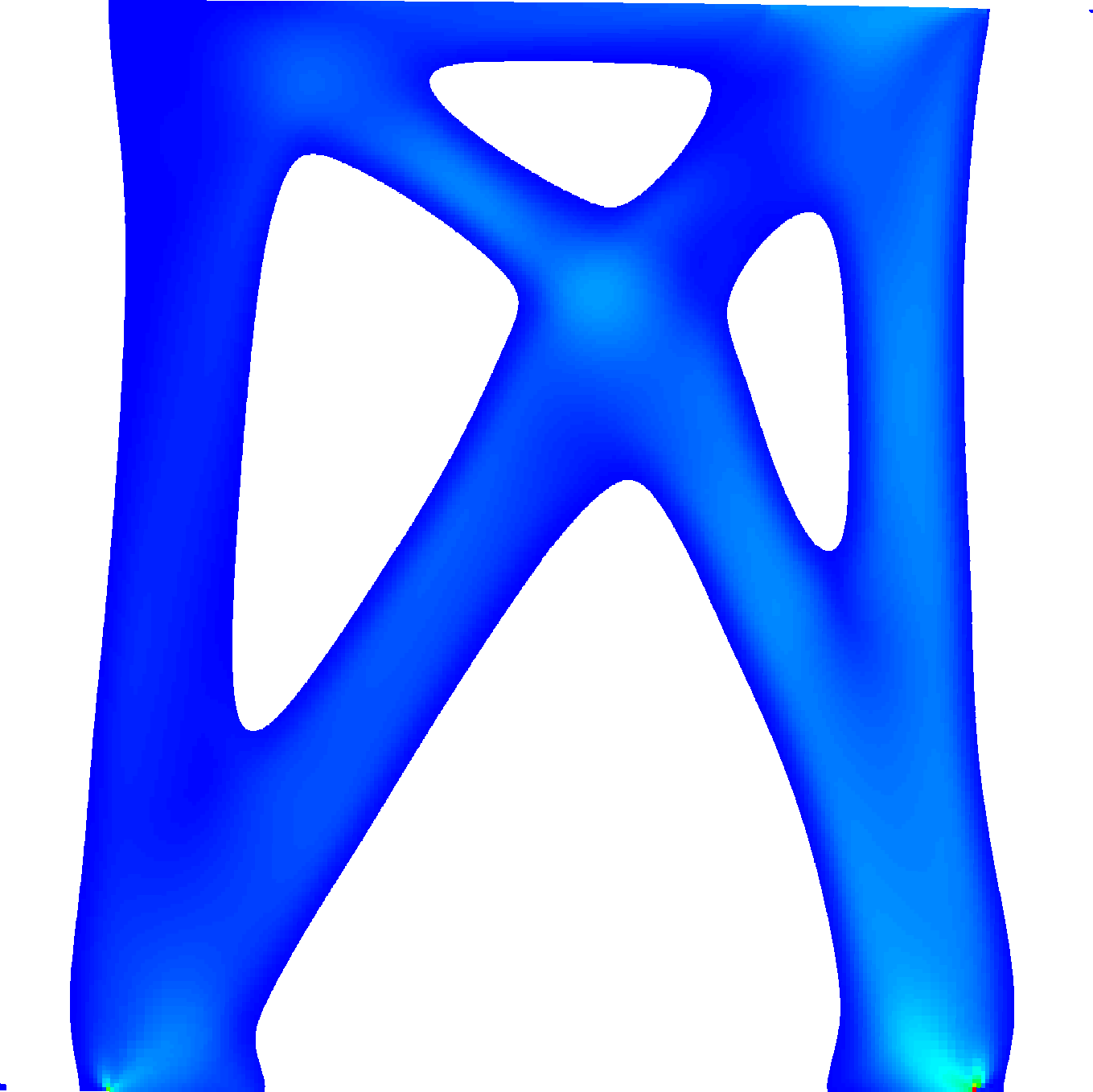}
  \end{subfigure}
  \begin{subfigure}[t]{0.09\textwidth}
    \includegraphics[width=\textwidth]{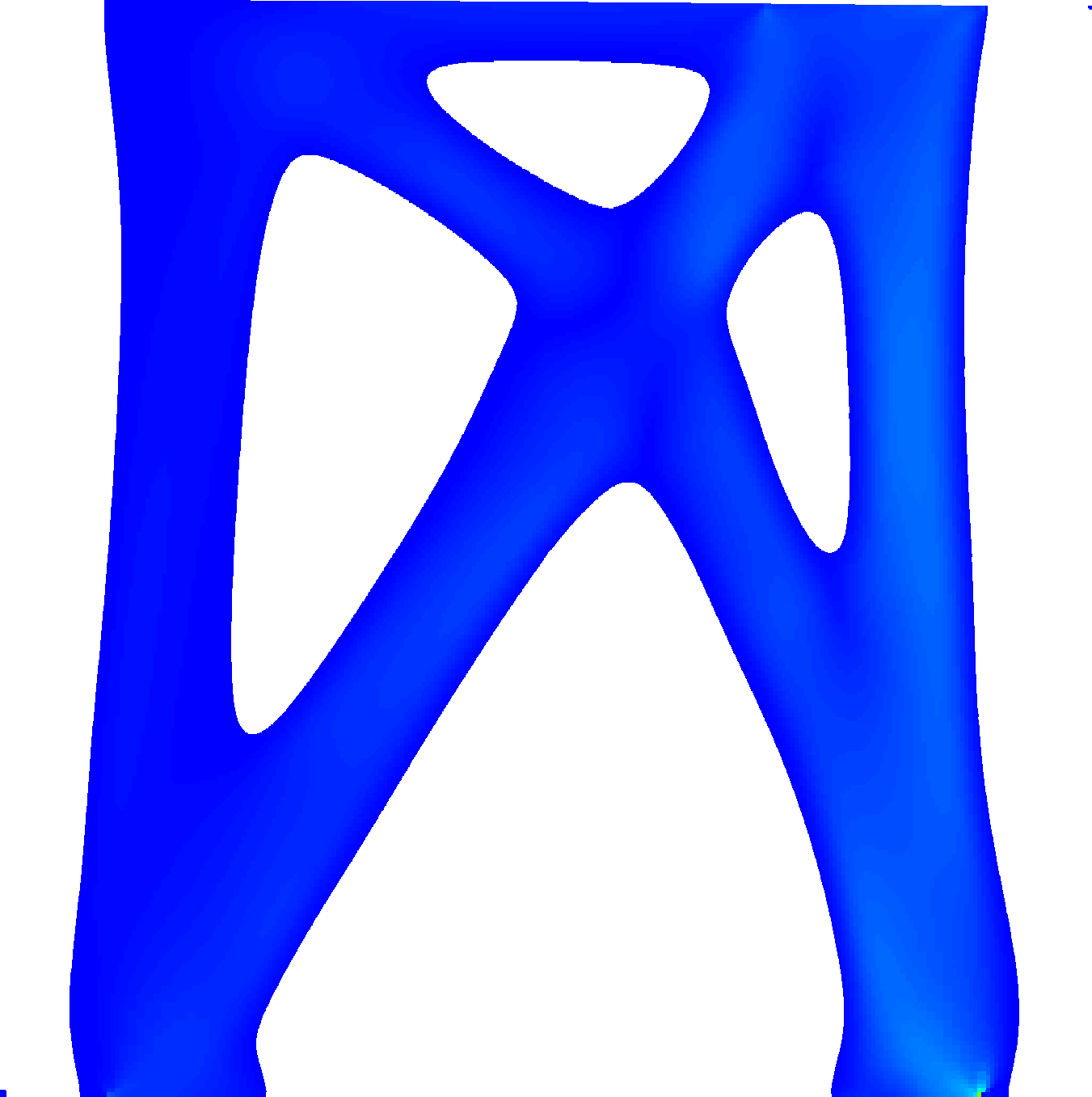}
  \end{subfigure}
  \begin{subfigure}[t]{0.09\textwidth}
    \includegraphics[width=\textwidth]{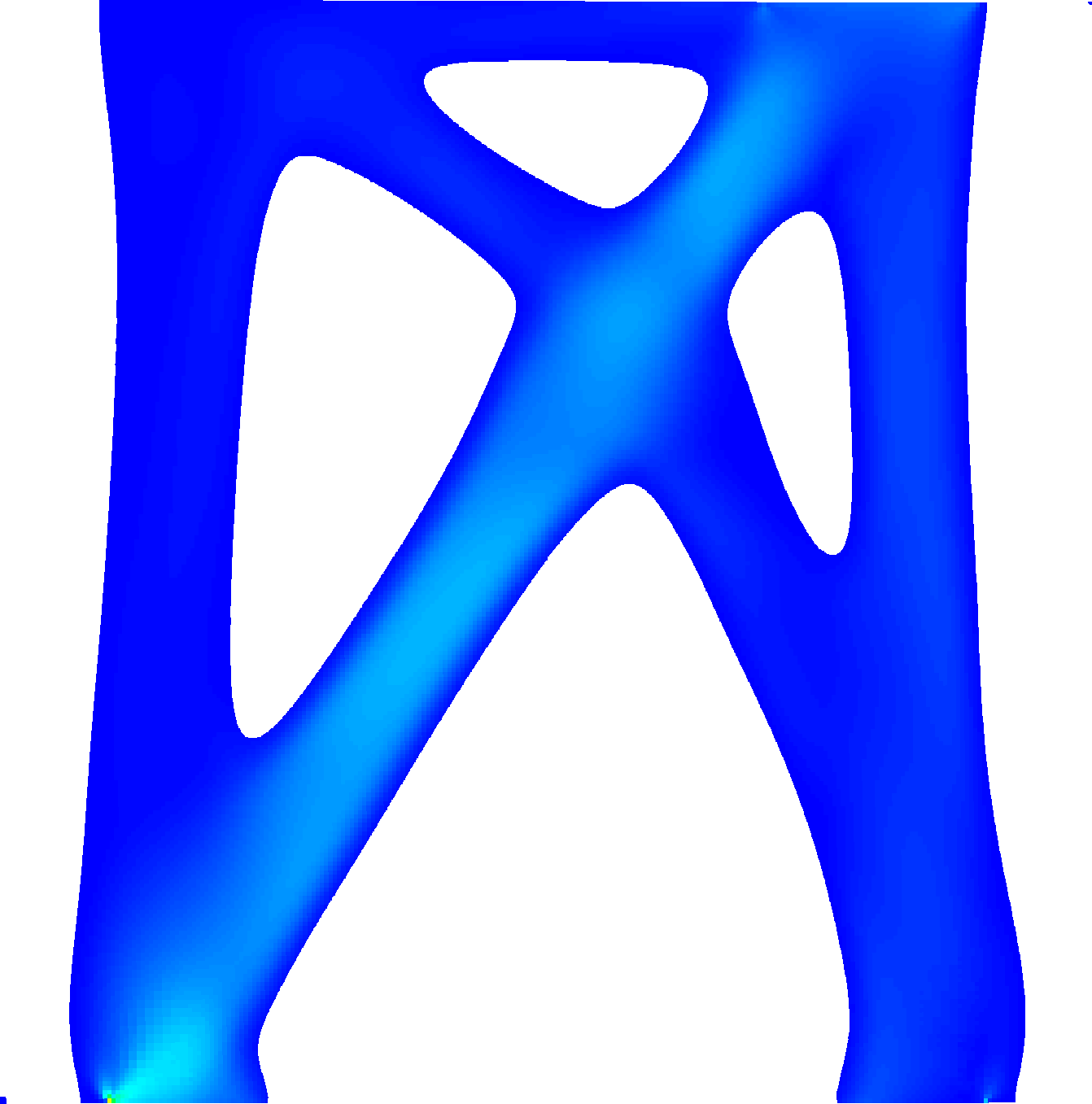}
  \end{subfigure}
  \begin{subfigure}[t]{0.09\textwidth}
    \includegraphics[width=\textwidth]{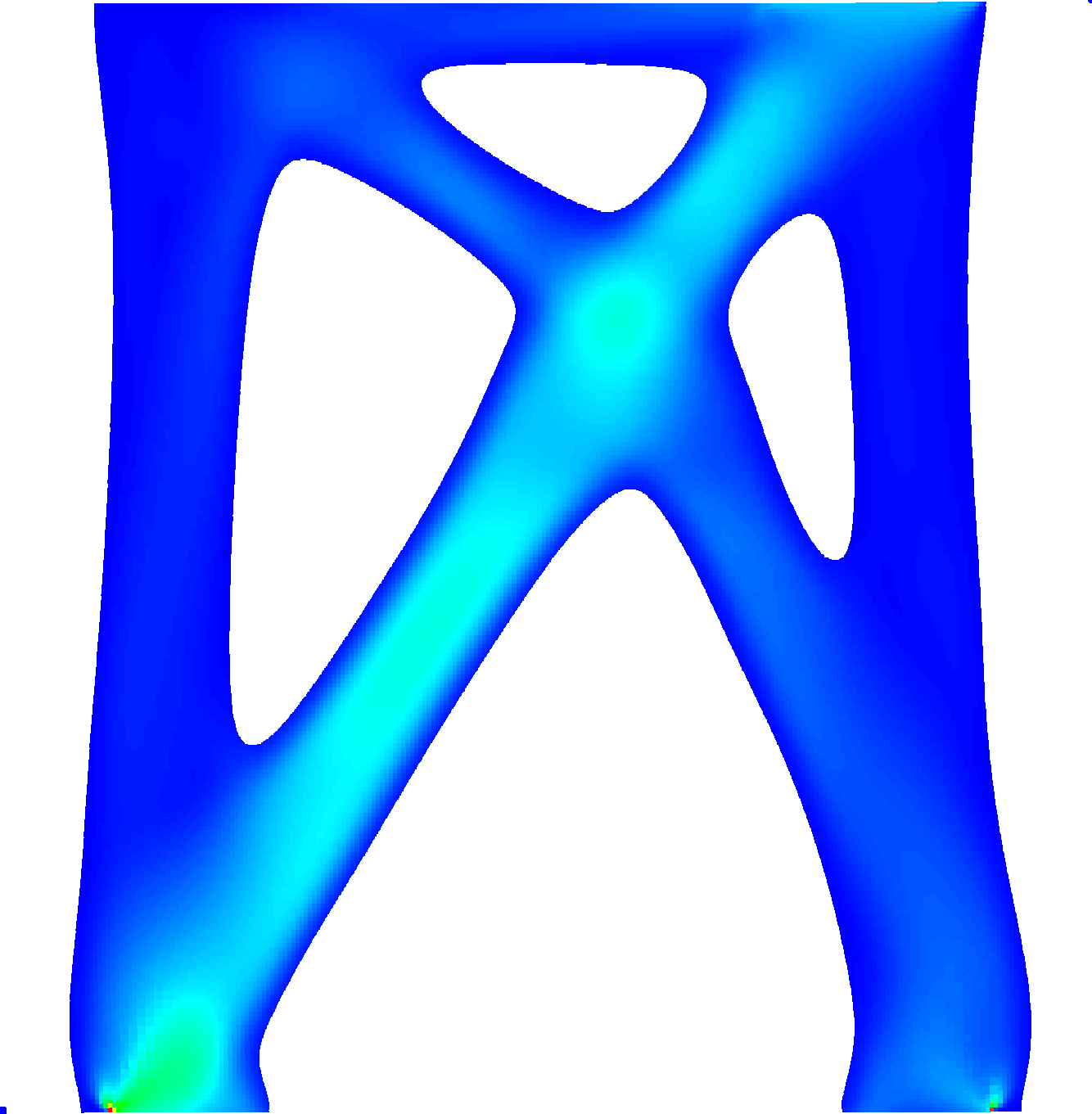}
  \end{subfigure}\\[1ex]
  \hrule \vspace{3ex}
  \begin{subfigure}[t]{0.09\textwidth}
    \includegraphics[width=\textwidth]{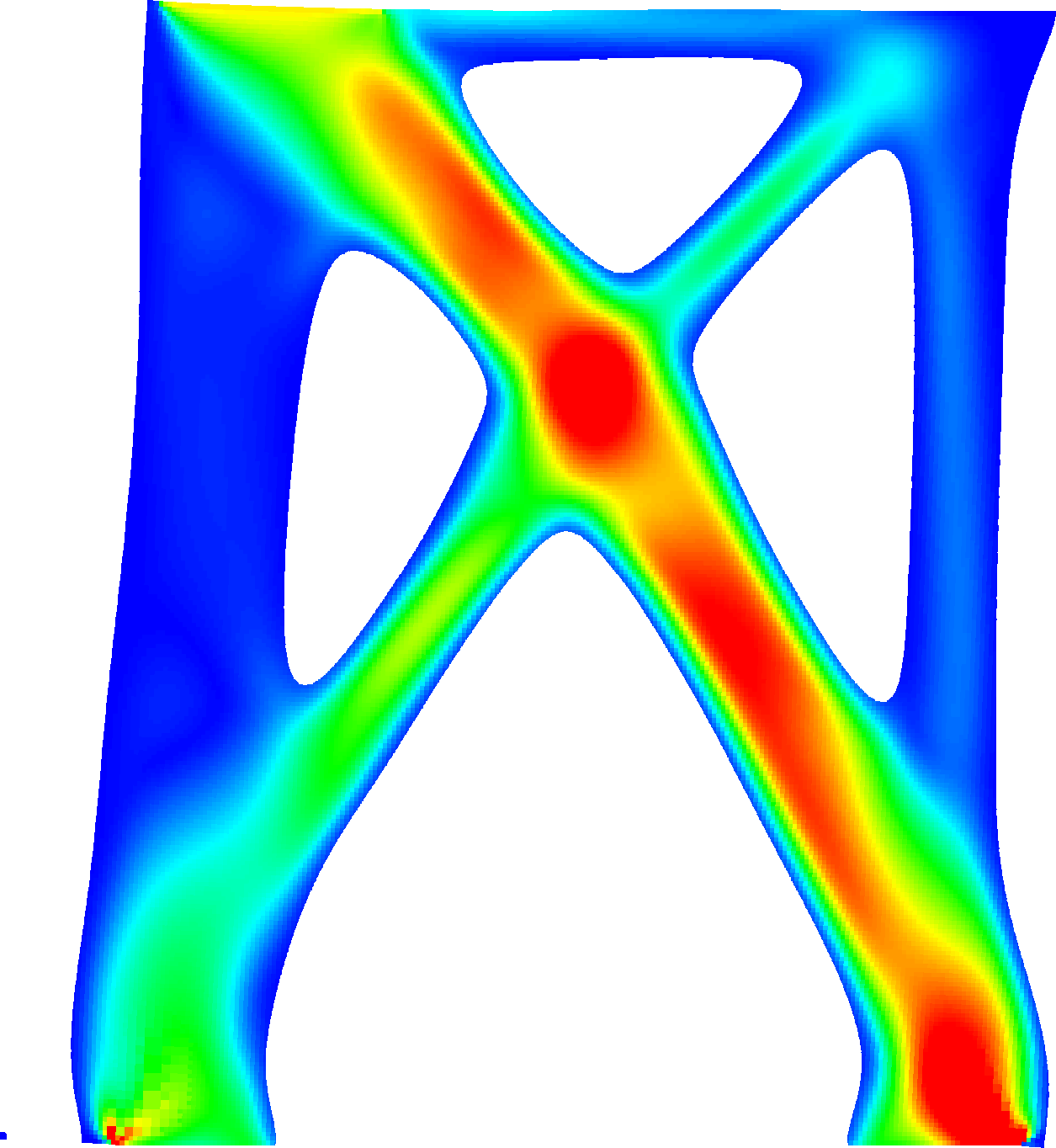}
  \end{subfigure}
  \begin{subfigure}[t]{0.09\textwidth}
    \includegraphics[width=\textwidth]{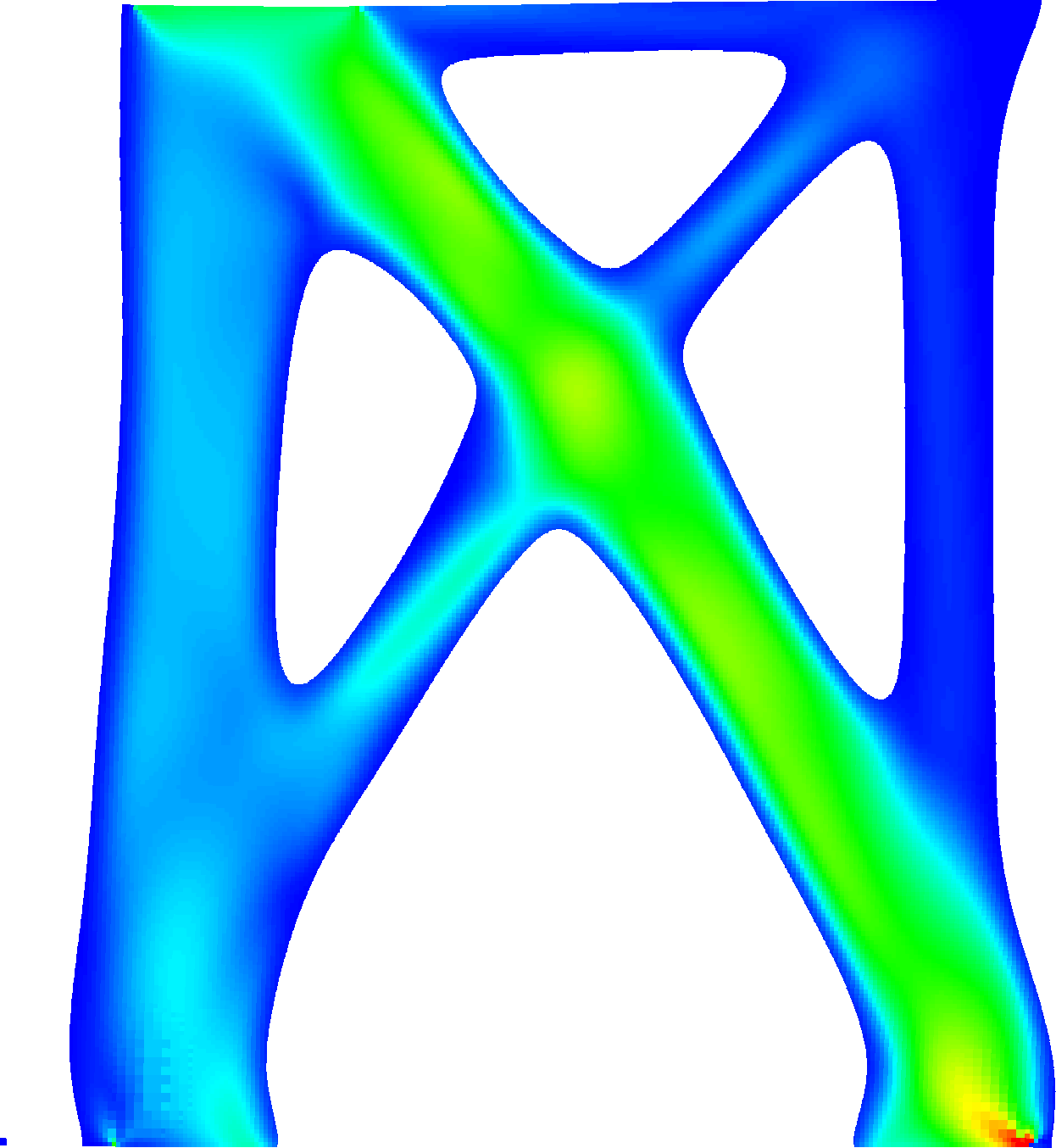}
  \end{subfigure}
  \begin{subfigure}[t]{0.09\textwidth}
    \includegraphics[width=\textwidth]{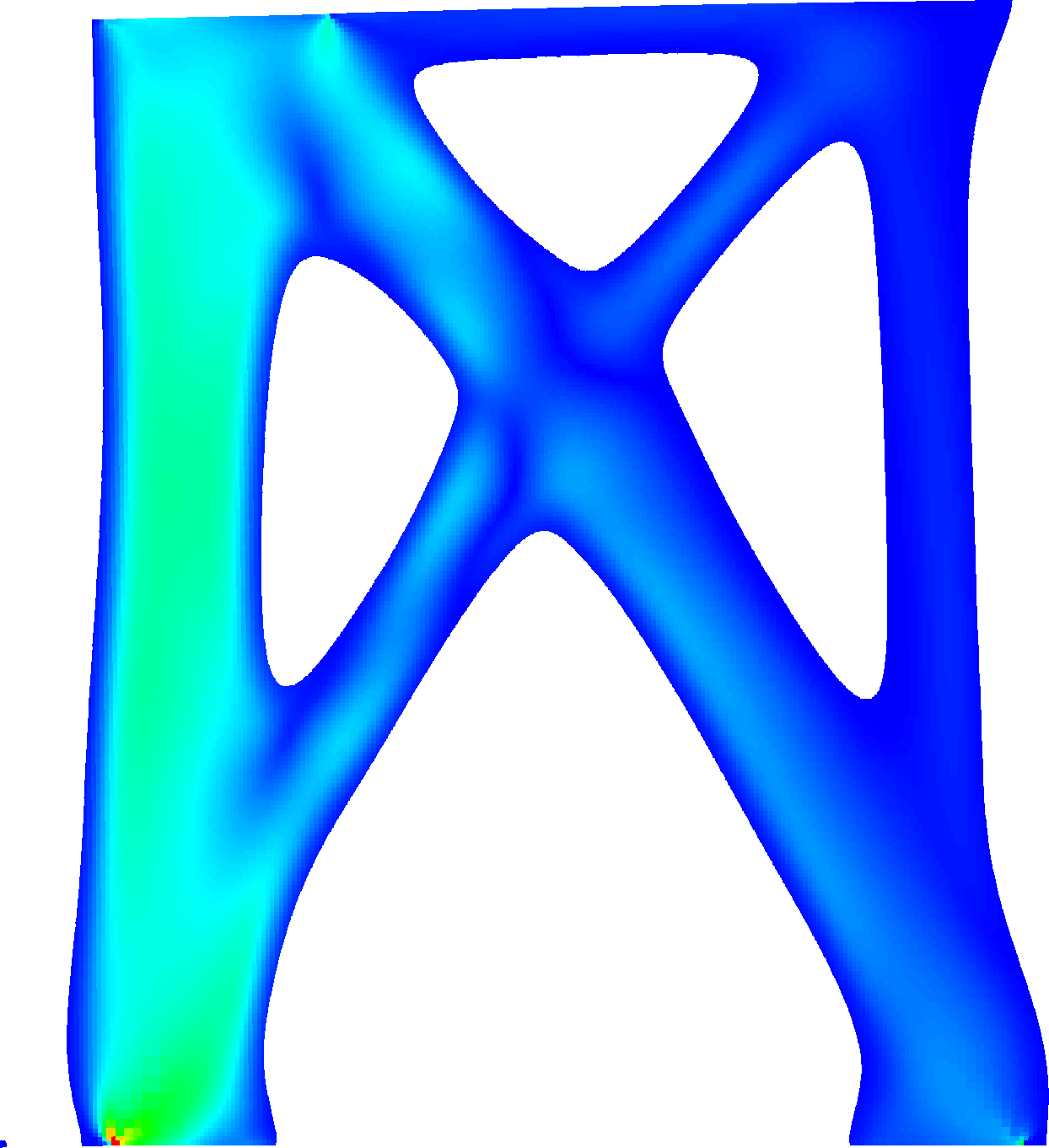}
  \end{subfigure}
  \begin{subfigure}[t]{0.09\textwidth}
    \includegraphics[width=\textwidth]{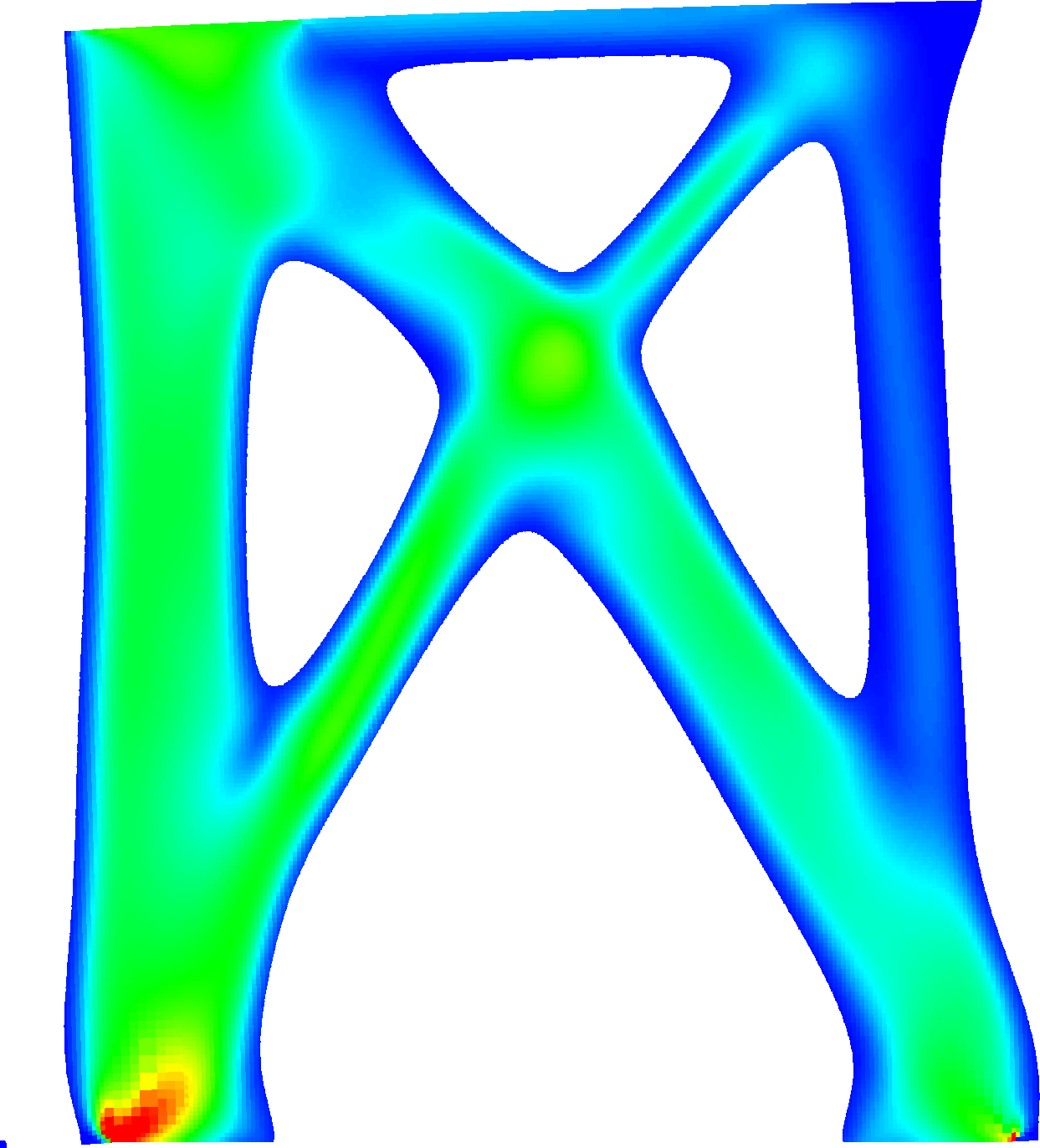}
  \end{subfigure}
  \begin{subfigure}[t]{0.09\textwidth}
    \includegraphics[width=\textwidth]{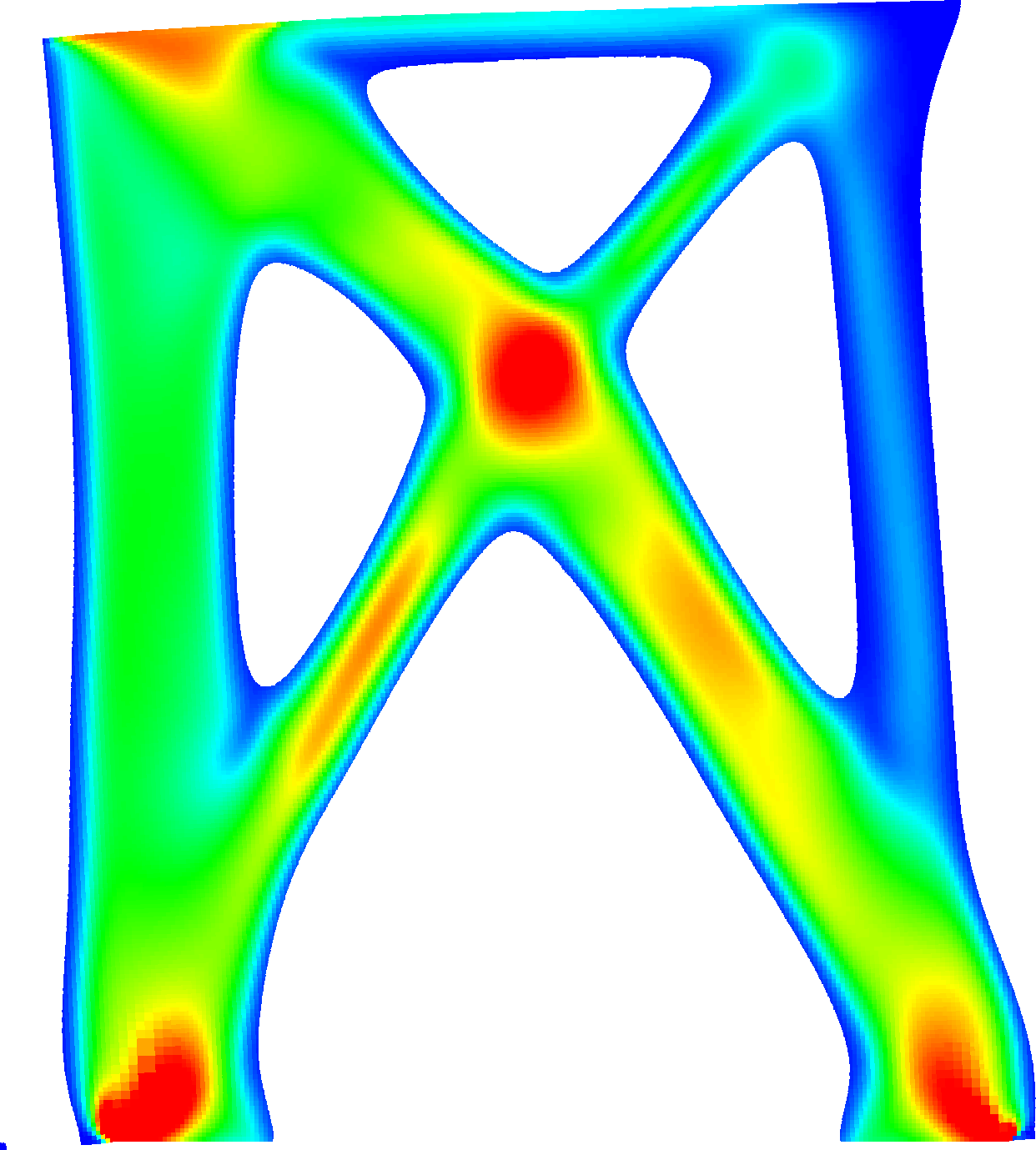}
  \end{subfigure}
  \begin{subfigure}[t]{0.09\textwidth}
    \includegraphics[width=\textwidth]{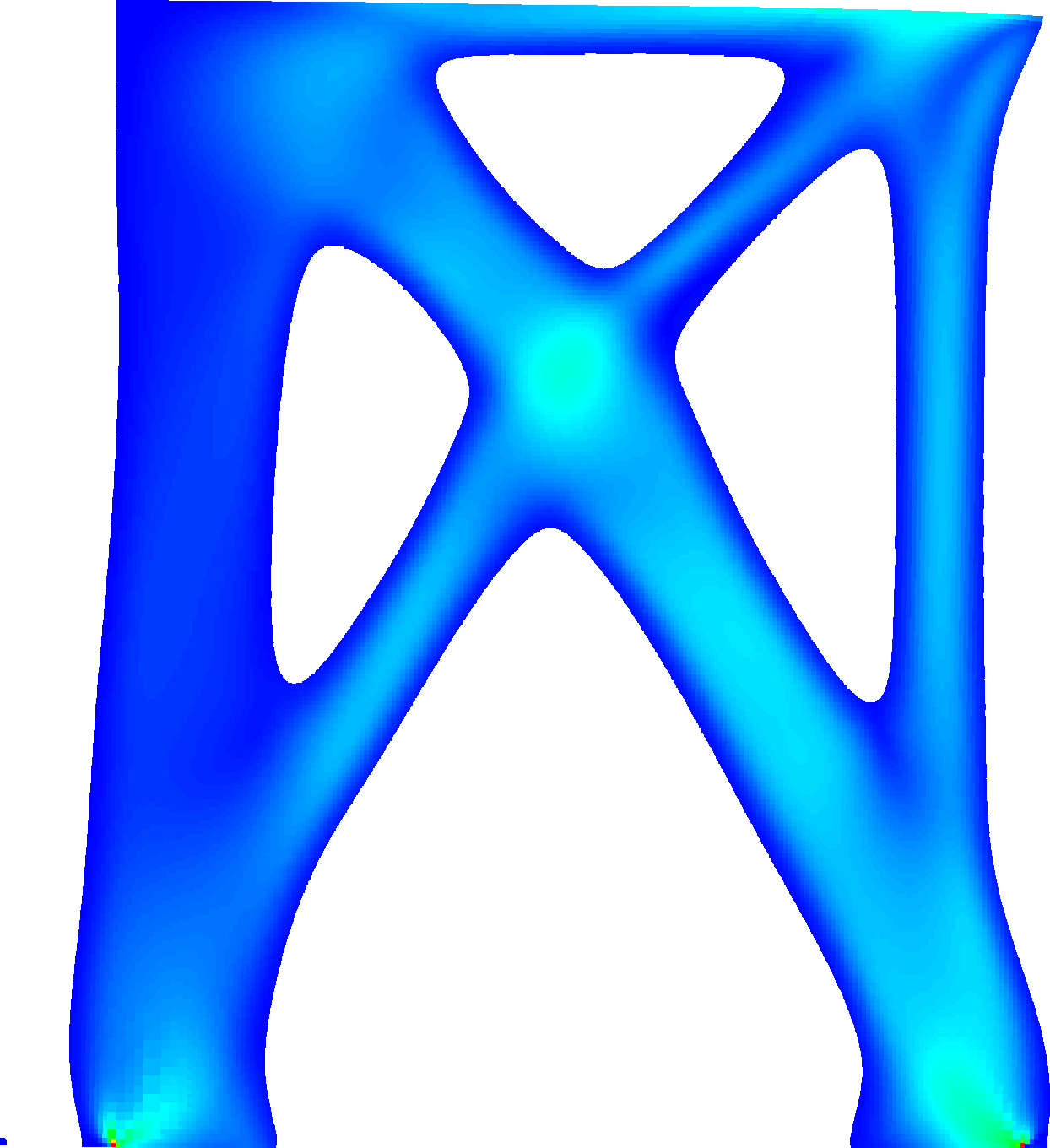}
  \end{subfigure}
  \begin{subfigure}[t]{0.09\textwidth}
    \includegraphics[width=\textwidth]{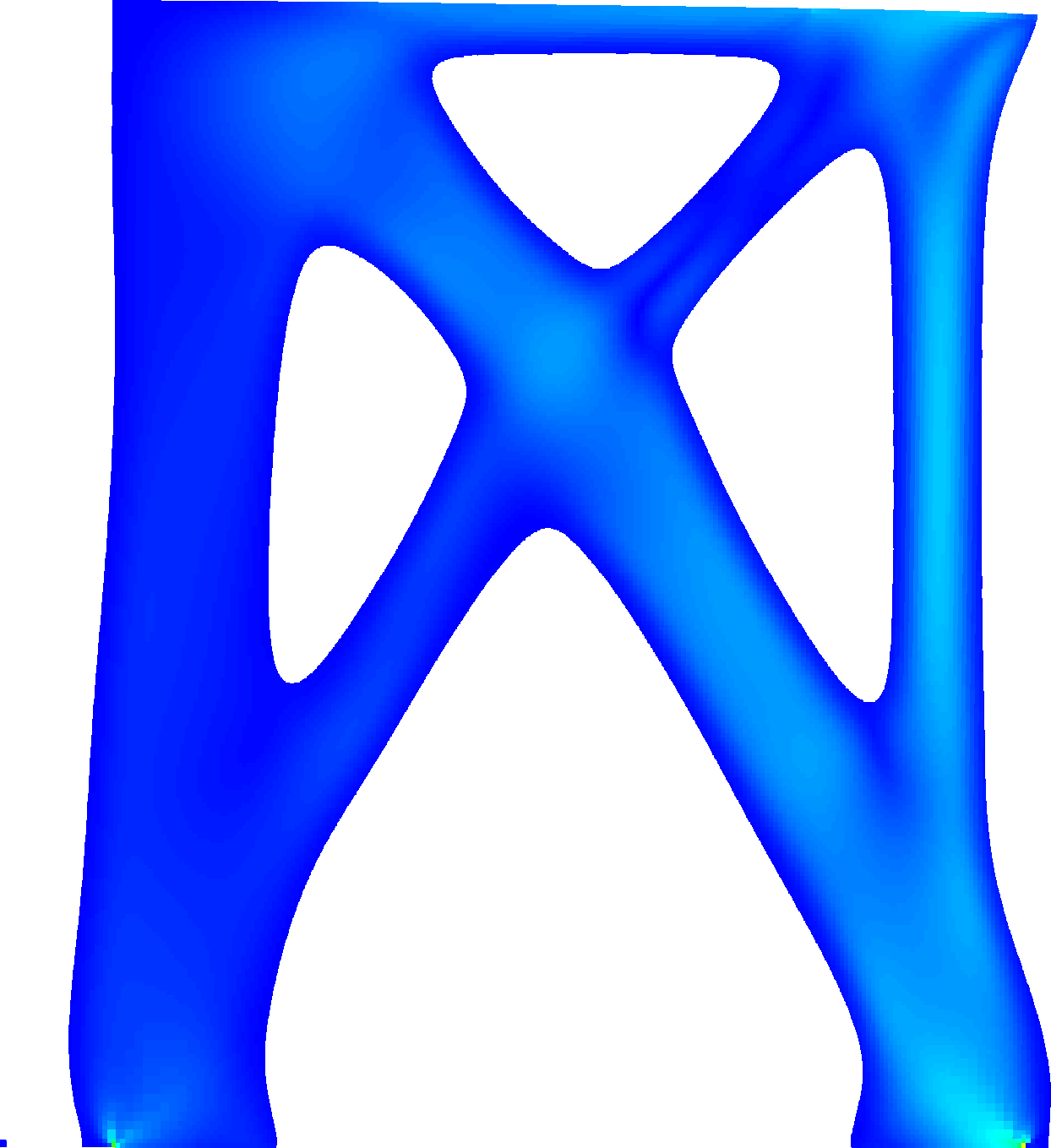}
  \end{subfigure}
  \begin{subfigure}[t]{0.09\textwidth}
    \includegraphics[width=\textwidth]{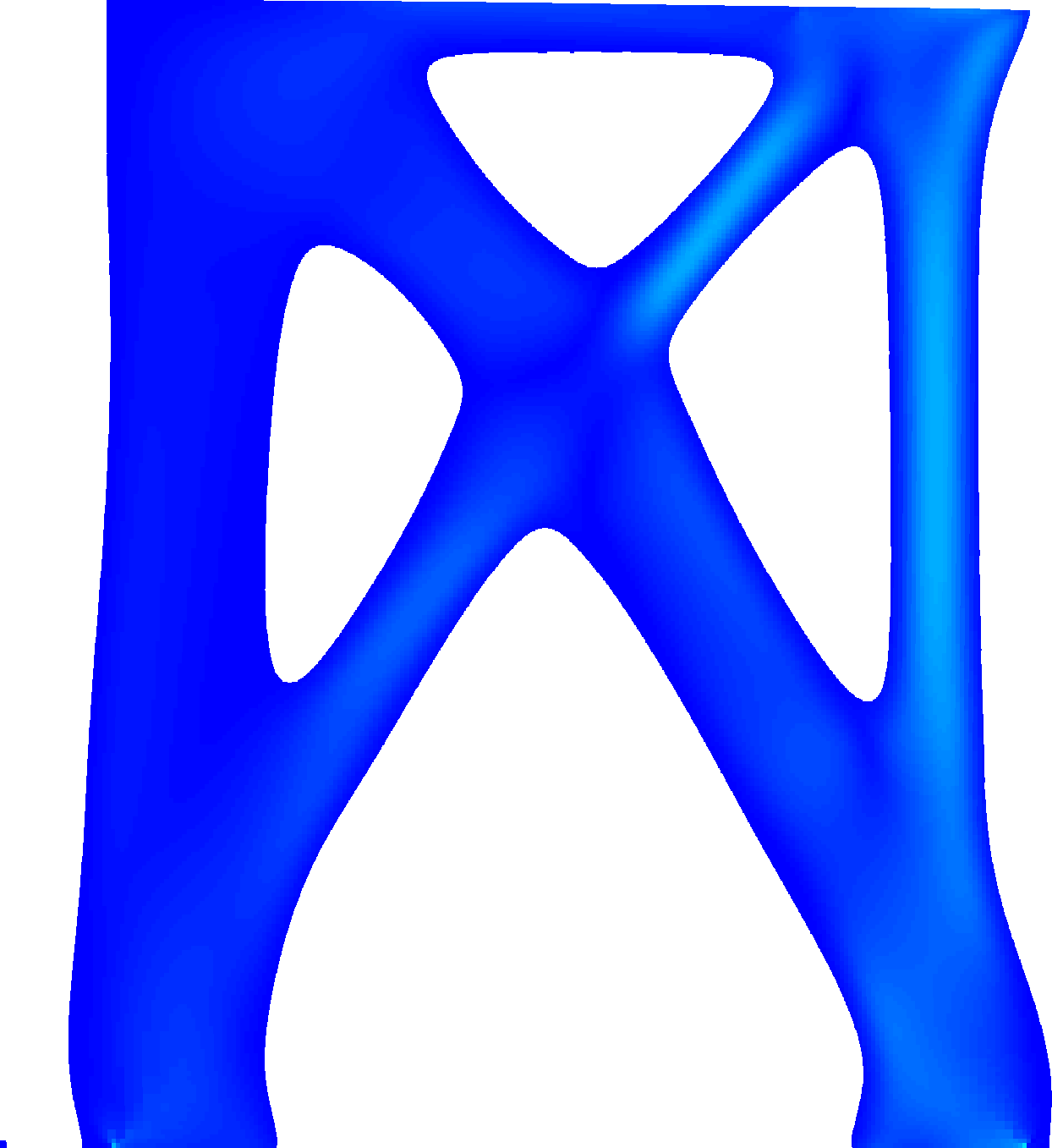}
  \end{subfigure}
  \begin{subfigure}[t]{0.09\textwidth}
    \includegraphics[width=\textwidth]{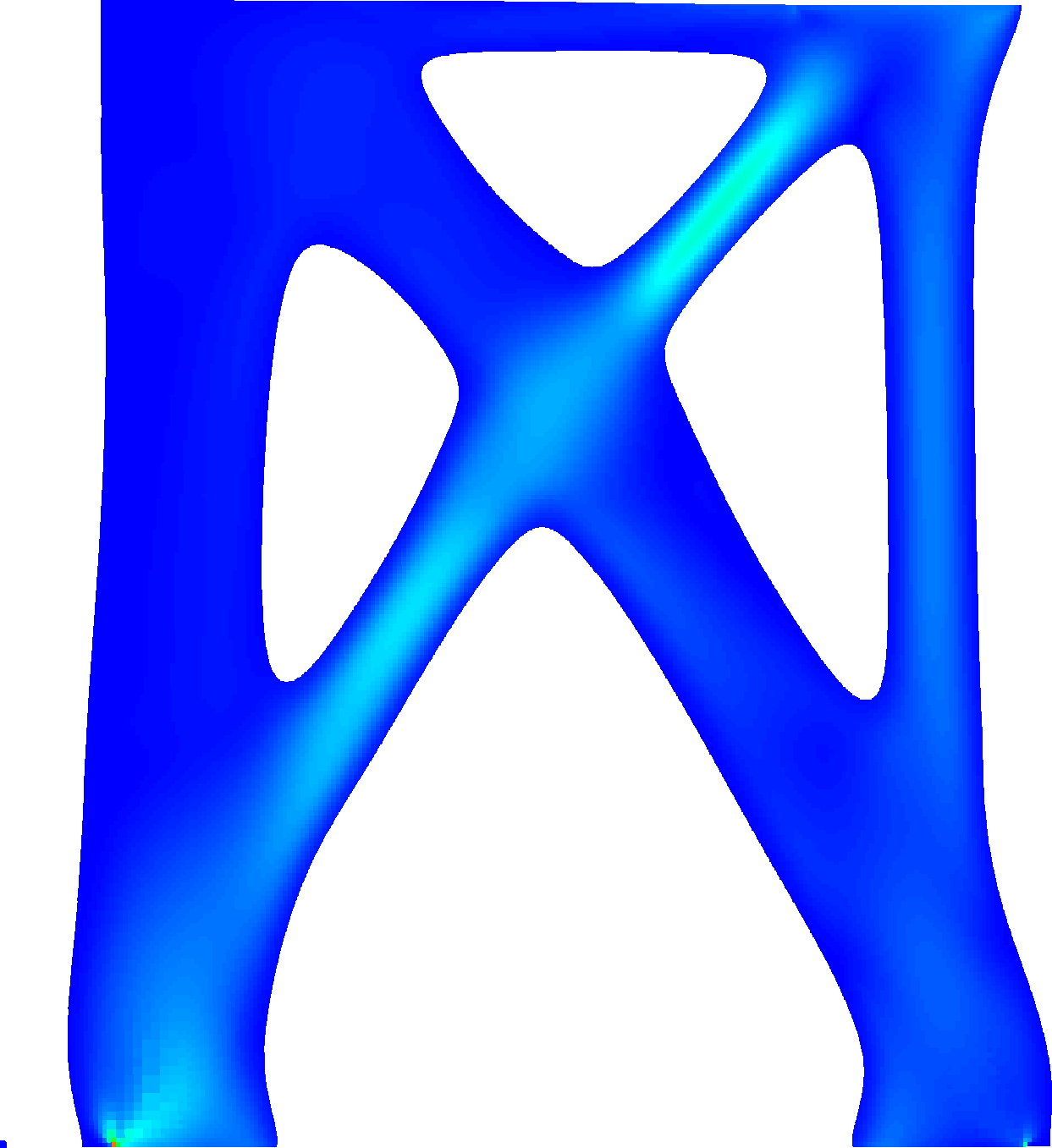}
  \end{subfigure}
  \begin{subfigure}[t]{0.09\textwidth}
    \includegraphics[width=\textwidth]{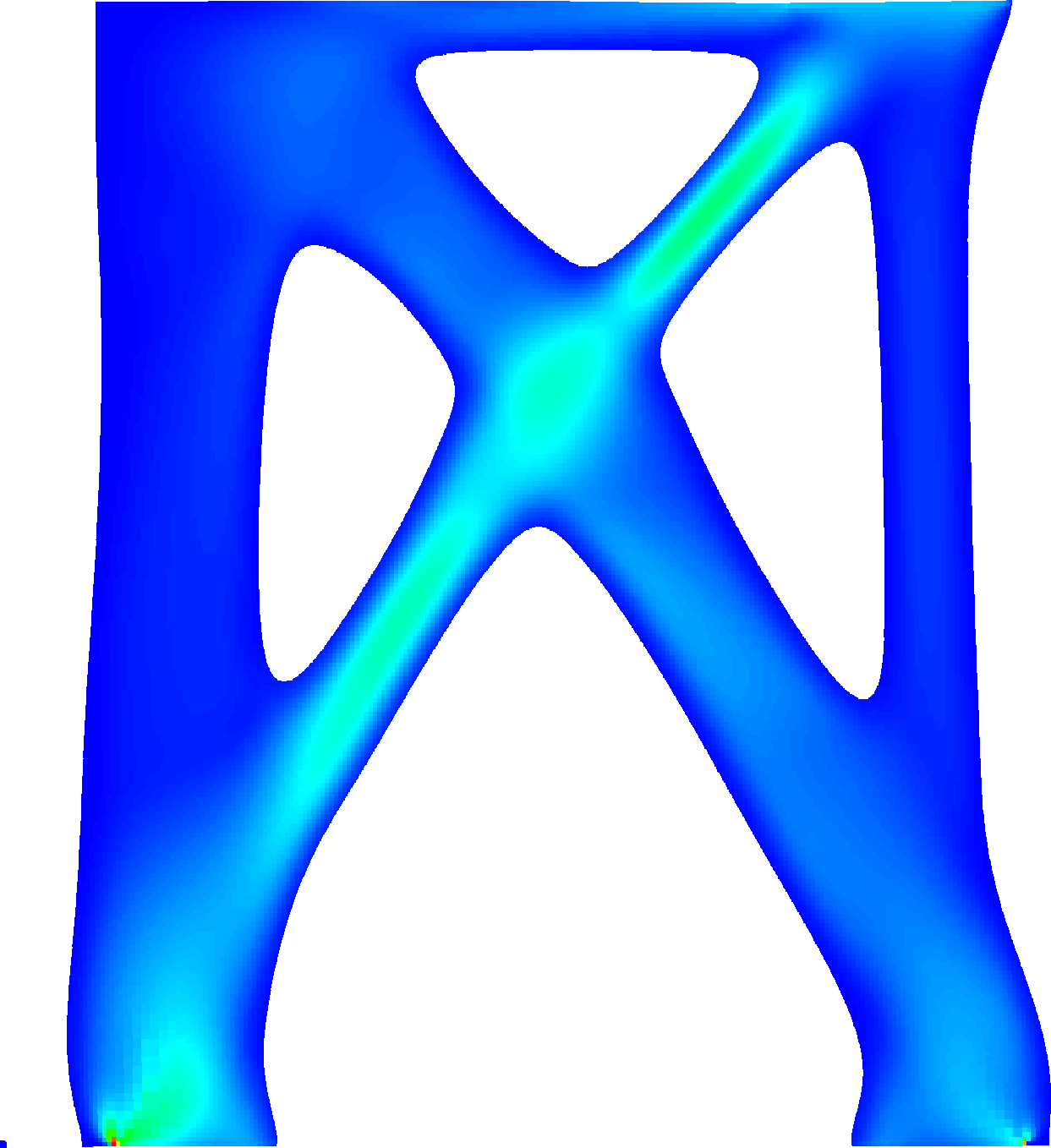}
  \end{subfigure}
  \caption{Thresholded stresses for each scenario in the weighted carrier
    plate setup, $1^{st}$ and $2^{nd}$ order dominance, colorcoded
    as $0$ \protect\includegraphics[height=1ex,width=5em]{images/res/cb_stress} $3.21$
    and $0$ \protect\includegraphics[height=1ex,width=5em]{images/res/cb_stress}
    $3.24$, stresses $> 3.21$ and $> 3.24$ are mapped to red.}
  \label{fig:stresscp}
\end{figure}

We consider a setup with varying absolute value of the applied load (ratio $\frac14$) and varying
probability of all scenarios, where high-probability loads come with a small load value.
Here, a phase field parameter of $3.13 \cdot 10^{-2}$ was used. After each
grid refinement, $\epsilon$ was multiplied by $0.75$.
The smallest $\epsilon$ was  $1.32 \cdot 10^{-2}$.
The value of $\heavreg$ decreased from $1$
to  $1.95 \cdot 10^{-3}$ and $7.81 \cdot 10^{-3}$
for the first and second order optimizations, respectively.
The different load configurations and the benchmark shape are displayed in Fig. \ref{fig:rescp} together with the cumulative distribution functions 
and the integrated survival function.
The corresponding cumulative distribution functions and integrated survival functions show
that the constraints are fulfilled for the first and second order dominance model.
Fig. \ref{fig:stresscp} shows the distribution of von Mises stresses on the hard material for all
10 scenarios.

\subsection*{Acknowledgment}
This work was supported by the Deutsche Forschungsgemeinschaft through the
Collaborative Research Center 1060 {\em  Mathematics of Emergent Effects}.
Furthermore, the   third author acknowledges support by the Deutsche Forschungsgemeinschaft in the  
 Collaborative Research Center TRR~154 {\em Mathematical Modelling, Simulation and Optimization Using the Example of Gas Networks}.
\renewcommand{\refname}{References}
\bibliographystyle{acm}
\bibliography{bibtex/all,bibtex/own}

\end{document}